\documentclass[12pt]{article}
\usepackage{amsfonts, amsmath, amsthm, amssymb, stmaryrd}
\usepackage[all]{xy}

\newtheorem{theorem}{Theorem}[subsection]
\newtheorem{lemma}[theorem]{Lemma}

\newtheorem{cor}[theorem]{Corollary}
\newtheorem{prop}[theorem]{Proposition}
\theoremstyle{definition}
\newtheorem{defn}[theorem]{Definition}
\newtheorem{remark}[theorem]{Remark}
\newtheorem{hypo}[theorem]{Hypothesis}
\numberwithin{equation}{theorem}

\newcommand\AAA{\mathbb{A}}

\newcommand\NN{\mathbb{N}}
\newcommand\PP{\mathbb{P}}
\newcommand\QQ{\mathbb{Q}}
\newcommand\RR{\mathbb{R}}
\newcommand\ZZ{\mathbb{Z}}
\newcommand\calC{\mathcal{C}}
\newcommand\calD{\mathcal{D}}
\newcommand\calE{\mathcal{E}}
\newcommand\calF{\mathcal{F}}
\newcommand\calG{\mathcal{G}}

\newcommand\calO{\mathcal{O}}
\newcommand\calR{\mathcal{R}}
\newcommand\calS{\mathcal{S}}

\newcommand\ida{\mathfrak{a}}
\newcommand\idm{\mathfrak{m}}
\newcommand\be{\mathbf{e}}
\newcommand\bv{\mathbf{v}}
\newcommand\bw{\mathbf{w}}
\newcommand\bx{\mathbf{x}}
\newcommand\hattimes{\widehat{\otimes}}

\newcommand\dual{\vee}
\newcommand\del{\partial}
\newcommand{\tilP}{\tilde{P}}

\DeclareMathOperator{\alg}{alg}
\DeclareMathOperator{\Aut}{Aut}
\DeclareMathOperator{\coker}{coker}

\DeclareMathOperator{\Frac}{Frac}
\DeclareMathOperator{\Gal}{Gal}
\DeclareMathOperator{\Hom}{Hom}
\DeclareMathOperator{\im}{im}
\DeclareMathOperator{\inte}{int}
\DeclareMathOperator{\loc}{loc}
\DeclareMathOperator{\Maxspec}{Maxspec}

\DeclareMathOperator{\ord}{ord}
\DeclareMathOperator{\pr}{pr}
\DeclareMathOperator{\qu}{qu}
\DeclareMathOperator{\rig}{rig}

\DeclareMathOperator{\spe}{sp}
\DeclareMathOperator{\Spec}{Spec}
\DeclareMathOperator{\Spf}{Spf}
\DeclareMathOperator{\Tor}{Tor}

\DeclareMathOperator{\Tr}{Tr}
\DeclareMathOperator{\Trace}{Trace}

\newcounter{fixmectr}

\begin{document}

\title{Finiteness of rigid cohomology with coefficients}
\author{Kiran S. Kedlaya \\
Department of Mathematics, Room 2-165 \\ Massachusetts Institute
of Technology \\ 77 Massachusetts Avenue \\
Cambridge, MA 02139 \\
\texttt{kedlaya@math.mit.edu}}
\date{November 2, 2005}

\maketitle

\begin{abstract}
We prove that for any field $k$ of characteristic $p>0$,
any separated scheme $X$ of finite type over $k$, and any overconvergent
$F$-isocrystal $\calE$ over $X$, the rigid cohomology 
$H^i_{\rig}(X, \calE)$ and rigid cohomology with compact supports
$H^i_{c,\rig}(X, \calE)$ are finite dimensional vector spaces over
an appropriate $p$-adic field. We also
establish Poincar\'e duality and the K\"unneth formula with coefficients.
The arguments use a pushforward construction in relative dimension $1$, based
on a relative version of Crew's conjecture on the quasi-unipotence
of certain $p$-adic differential equations.
\end{abstract}

\tableofcontents

\section{Introduction}

\subsection{Rigid cohomology: an overview}

The purpose of this paper is to prove some basic structural properties
about Berthelot's rigid cohomology for algebraic varieties
over fields of characteristic $p>0$.
This theory unifies most of the existing $p$-adic
cohomology theories, including crystalline cohomology (developed by
Berthelot, Ogus et al.; see \cite{bib:ber0} and \cite{bib:bo}),
Monsky-Washnitzer cohomology (see \cite{bib:mw} and \cite{bib:vdp}), 
Serre's Witt vector
cohomology (and its descendant, the de Rham-Witt complex of Deligne
and Illusie; see \cite{bib:illusie}), Lubkin's
bounded Witt vector construction (see \cite{bib:lubkin}), et al. 

Rigid cohomology can be viewed as a version of
algebraic de Rham cohomology in positive characteristic; in that theory,
one defines cohomology of a variety by embedding it into a smooth variety
and working on its formal neighborhood. To get a good theory (in which
the divisions by $p$ inherent in integration make sense), Berthelot
passes to rigid analytic spaces and replaces the formal neighborhoods
by certain analytic neighborhoods (the ``strict neighborhoods'').

Grothendieck's proof that algebraic de Rham cohomology (of a variety
over a field of characteristic zero) is finite dimensional uses a reduction
to the case of a smooth proper variety via Hironaka's resolution of 
singularities. Similarly, using de Jong's alterations theorem to stand in for
resolution of singularities, Berthelot \cite{bib:ber2}, \cite{bib:ber6}
succeeded in proving finite dimensionality of rigid cohomology with
constant coefficients, along with Poincar\'e duality and the K\"unneth
formula.

However, one can also work with nonconstant coefficients, analogous
to local systems in de Rham cohomology; the most convenient coefficient
objects are the ``overconvergent $F$-isocrystals''.
In this paper, we generalize Berthelot's results to
the cohomology of these coefficient objects.

\subsection{The main results}

Here are our results on the rigid cohomology of varieties over
a field $k$ of characteristic $p>0$. Let $K$ be a complete
discrete valuation ring of mixed characteristics $(0,p)$
with residue field $k$.
\begin{theorem} \label{thm:finite}
Let $\calE/K$ be an overconvergent $F$-isocrystal on a separated
scheme $X$ of finite type over $k$. Then
the rigid cohomology spaces
$H^i_{\rig}(X/K, \calE)$ are finite
dimensional $K$-vector spaces for all $i$.
\end{theorem}
\begin{theorem} \label{thm:finite2}
Let $\calE/K$ be an overconvergent $F$-isocrystal on a separated
scheme $X$ of finite type over $k$. Then
the rigid cohomology spaces with compact supports
$H^i_{c,\rig}(X/K, \calE)$ are finite dimensional $K$-vector spaces
for all $i$.
\end{theorem}
\begin{theorem}[Poincar\'e duality] \label{thm:poincare}
Let $\calE/K$ be an overconvergent $F$-isocrystal on a smooth
$k$-scheme $X$ of pure dimension $d$. Then
for any closed subscheme $Z$ of $X$, there are natural perfect pairings
\begin{equation} \label{eq:pairing}
H^i_{Z,\rig}(X/K, \calE) \otimes_K H^{2d-i}_{c,\rig}(Z/K, \calE^\dual) 
\to K.
\end{equation}
\end{theorem}
\begin{theorem}[K\"unneth formula] \label{thm:kunneth}
For $i=1,2$,
let $\calE_i/K$ be an overconvergent $F$-isocrystal on a
separated scheme $X_i$ of finite type over $k$.
Put $X = X_1 \times_k X_2$,
let $\pr_i: X_1 \times X_2 \to X_i$ be the
natural projections, and put
$\calE = \pr_1^* \calE_1 \otimes \pr_2^* \calE_2$.
Then there are natural isomorphisms
\begin{equation} \label{eq:kun1}
\bigoplus_{j+l=i} H^j_{c,\rig}(X_1/K, \calE_1) \otimes_K
H^l_{c,\rig}(X_2/K, \calE_2) \to H^i_{c,\rig}(X/K, \calE).
\end{equation}
Moreover, if $X_1$ and $X_2$ are smooth over $k$, then for
any closed subschemes $Z_i$ of $X_i$,
there are also natural isomorphisms
\begin{equation} \label{eq:kun2}
\bigoplus_{j+l=i} H^j_{Z_1, \rig}(X_1/K, \calE_1) \otimes_K
H^l_{Z_2, \rig}(X_2/K, \calE_2) \to H^i_{Z_1 \times_k Z_2/K, \rig}(X, \calE).
\end{equation}
\end{theorem}

\subsection{Previously known cases}

Numerous special cases of these results have
been established previously, modulo Berthelot's comparison theorems
\cite[Propositions~1.9 and 1.10]{bib:ber2}
between rigid, crystalline and Monsky-Washnitzer (MW) cohomology;
we discuss the various strategies behind these methods in the next
section.
For now, we simply list the previously known special cases.
\begin{itemize}
\item
For $\calE$ trivial and $X$ smooth proper, all four results were
proved in crystalline cohomology by Berthelot \cite{bib:ber0}.
Berthelot's result also applies to the case where $\calE$ is obtained
from an $F$-crystal of locally free coherent sheaves,
 via the comparison theorem between
crystalline and rigid cohomologies with such coefficients
\cite[Theorem~0.7.7]{bib:ogus}.
\item
For $\calE$ trivial and $X$ a smooth affine curve, Theorem~\ref{thm:finite}
was proved in MW cohomology by Monsky \cite{bib:monsky}.
\item
For $\calE$ trivial and $X$ smooth affine, Theorem~\ref{thm:finite}
was proved in MW cohomology by Mebkhout \cite{bib:meb1},
using results of Christol and Mebkhout \cite{bib:cm3}.
\item
For $\calE$ trivial and $X$ smooth, Theorems~\ref{thm:finite} 
and~\ref{thm:finite2} were proved by Berthelot in \cite{bib:ber2},
and Theorems~\ref{thm:poincare} and~\ref{thm:kunneth} in \cite{bib:ber6}.
\item
For $\calE$ trivial, Theorems~\ref{thm:finite} and~\ref{thm:finite2}
were proved
by Grosse-Kl\"onne \cite{bib:gk2} and by Tsuzuki \cite{bib:tsu6}.
\item
For $X$ a smooth affine curve, the first three results were proved
conditionally by Crew \cite{bib:crew2} (see below).
\item
For $\calE$ unit-root and $X$ smooth, all four results were proved by
Tsuzuki \cite{bib:tsu5}.
\end{itemize}

The finiteness result of Crew, for $X$ a curve but $\calE$ arbitrary,
relies on the hypothesis that $\calE$ is quasi-unipotent at each point
of a compactification of $X$.
We will discuss this hypothesis further in the next section.

\subsection{Approaches to $p$-adic cohomology}
\label{subsec:philosophy}

Before giving a description of what is actually in the paper,
we take a moment to fit our methods into the already rich context
of the existing literature on $p$-adic cohomology. This discussion
focus on the past and present; for thoughts on the future,
see Section~\ref{subsec:further}.

One of the main strategies for approaching $p$-adic cohomology is
to extract information from the case of smooth proper varieties.
That is, one attempts to reduce consideration to crystalline cohomology.
This is the strategy adopted by Berthelot, Tsuzuki, Shiho, et al.;
the work of Grosse-Kl\"onne is in a similar vein.

A second approach, practiced by Christol and Mebkhout,
is perhaps the most direct heir to the legacy of Dwork and Robba. It
is more analytic in nature, using strong structural
properties of $p$-adic differential equations to control 
a suitable ``Dworkian'' cohomology theory; the most geometric of these
is the ``formal cohomology'' introduced by Monsky and Washnitzer. It
can be viewed as a close analogue of ``na\"\i ve'' algebraic de Rham
cohomology of smooth affine schemes in characteristic zero.

A third approach, introduced by Crew, attempts to mediate the two
extremes. It combines geometric arguments
with arguments from $p$-adic functional analysis. As in the 
Dworkian approach, the aim is to prove assertions directly
on affine schemes, then build general results from them, rather
than have to deduce results on affines from results about proper schemes.

In all three points of view, the same obstruction arises
to treating nonconstant coefficients: one needs an analogue for
$p$-adic differential equations of basic monodromy properties of
algebraic differential equations, or if one prefers, of Grothendieck's
local monodromy theorem in \'etale cohomology. This obstruction was isolated by
Crew, who posited the hypothesis that became known as ``Crew's conjecture''.
(It is also sometimes called the Crew-Tsuzuki conjecture, since the
``correct'' formulation first appears in \cite{bib:tsu3}.)
Briefly put, Crew's conjecture states that every $p$-adic differential
equation on an annulus with Frobenius structure has quasi-unipotent monodromy.
(We give the precise statement in Chapter~\ref{sec:rel}.)

Happily, this obstruction has recently been lifted.
Crew's conjecture has been proved independently by
Andr\'e \cite{bib:andre}, Mebkhout \cite{bib:mebkhout}, and the author
\cite{bib:me7}. (The unit-root case was proved earlier by Tsuzuki
\cite{bib:tsu1}; that case was a key ingredient to the aforementioned
results \cite{bib:tsu5}. Incidentally, it also plays a pivotal role in 
\cite{bib:me7}.) This circumstance makes it feasible to prove strong
structural properties of rigid cohomology with coefficients; however,
the ``raw'' form of Crew's conjecture is only sufficient for handling
curves. One must refine the result into a more global form for it to
be of cohomological value; one method of doing so was suggested by Shiho
(see Section~\ref{subsec:further}), but our method seems to be somewhat simpler.
We formulate a version of Crew's conjecture for a differential equation
on a family of annuli, and reduce it to the original result.

With a global version of Crew's conjecture in hand, we set off
to obtain cohomological results. However, we do not particularly adhere
to any of the approaches; our methods are something of a hodgepodge.
For most of our calculations, we work with Monsky-Washnitzer cohomology;
this is consonant with the approach
of Mebkhout, although we tend to prefer algebraic and geometric constructions
over $p$-adic differential equations. 
In particular, we construct ``generic'' higher direct images for families
of affine curves, use a Leray spectral sequence to relate the cohomology
of the total space to that of the direct images, then perform
d\'evissage.
The construction of the direct images, though algebraic in character,
uses the formalism of Crew's functional analytic approach.
However, to treat cohomology with compact supports and Poincar\'e
duality, we 
pass to the functional analytic approach of Crew, and to manipulate
cohomology on nonsmooth schemes, we invoke cohomological descent results
of Chiarellotto and Tsuzuki. (It may also be possible to use Grosse-Kl\"onne's
methods to treat the nonsmooth case.)

\subsection{Structure of the paper}

We now give an overview of the content of the various sections of the paper.

In Chapter~\ref{sec:prelim}, we introduce the basic rings with which
we compute: affinoid algebras (the building blocks of rigid analytic
geometry), dagger algebras (the building blocks of Monsky-Washnitzer
cohomology and of Gr\"osse-Klonne's theory of ``dagger spaces''),
the classical Robba ring, and our relative generalization thereof.

In Chapter~\ref{sec:frobdiff}, we construct auxiliary structures
on the rings of Chapter~\ref{sec:prelim} for use in rigid cohomology:
modules of differential forms, and Frobenius lifts. These are
used to define $(\sigma, \nabla)$-modules, the coefficient objects
in the rigid cohomology of smooth affine schemes. (These may be thought
of as ``rigidified local systems'': the differential forms provide
the structure of a connection, and the Frobenius rigidifies the connection
so as to make up for the disconnectedness of $p$-adic topologies,
specifically by filling the role played by analytic continuation
in complex analysis. On top of this, the Frobenius plays a role in controlling
local monodromy analogous
to that of a variation of Hodge structure on a punctured disc.)

In Chapter~\ref{sec:rigid}, we introduce Berthelot's rigid cohomology
and its coefficient objects, the overconvergent $F$-isocrystals. These
are closely related to the $(\sigma, \nabla)$-modules of the previous
chapter, which are essentially ``local models'' of overconvergent 
$F$-isocrystals. We also recall (mostly without proof) some of the
necessary formal properties of Berthelot's construction.

In Chapter~\ref{sec:rel}, we state and prove the relative form of
the $p$-adic local monodromy theorem. This states essentially that 
a family of $p$-adic differential equations on annuli, equipped with
Frobenius structure, which by the original $p$-adic local monodromy theorem
is quasi-unipotent over the generic point, is also quasi-unipotent
after restricting the base to an open dense subset. The proof is a 
reduction to
the original form of Crew's conjecture, using a solution to a differential
equation over the generic point as a guide for solving the relative equation.

In Chapter~\ref{sec:cohcurves}, we apply the
(nonrelative) $p$-adic local monodromy theorem to prove,
or rather to reprove, the finite dimensionality of rigid cohomology
with coefficients on an affine curve. As mentioned above, assuming the
$p$-adic local monodromy theorem, this has already been shown by Crew
in \cite{bib:crew2}. However, our approach is very different: we use
algebraic rather than topological consequences of the quasi-unipotence
of differential equations. In any case, this chapter can be omitted at
the preference of the reader.

In Chapter~\ref{sec:push}, we construct ``relative Monsky-Washnitzer
cohomology'' for simple families of curves, namely products of a base
with the affine line. (The basic strategy is enough to construct
pushforwards more generally, at least in relative dimension 1,
but this topic deserves a separate paper in itself.)
We also show that the construction commutes with flat base change.

In Chapter~\ref{sec:push2}, we compute the cohomology with compact
supports of overconvergent $F$-isocrystals on the affine line,
and establish Poincar\'e duality there. For the latter, we must bring
in the topological machinery of \cite{bib:crew2}.

Finally, in Chapter~\ref{sec:final}, we prove our main results using
d\'evissage (i.e., induction on dimension), 
plus some other formal properties of rigid
cohomology that do not occur elsewhere in the paper.
As our methods are only adapted to smooth varieties, we must invoke
cohomological descent results of Chiarellotto and Tsuzuki to handle
nonsmooth varieties.

\subsection{Further remarks}
\label{subsec:further}

We conclude this introduction by pointing out some consequences,
both actual and potential, of the present paper.

As an application of the finiteness of rigid cohomology with
coefficients, one
can now give a proof of the full 
Weil
conjectures within the framework of $p$-adic cohomology,
independent of the theory
of \'etale cohomology.
Moreover, one can give a $p$-adic analogue of ``Weil II'',
the main theorem of Deligne's \cite{bib:deligne}.
We carry this out in a separate paper \cite{bib:me10}.
(The reader is invited to imagine other results in \'etale cohomology
that may have $p$-adic analogues; we will leave such speculation
for another occasion.)

Our analogue of Weil II is not perfect, however, because our coefficient
objects only mimic the category
of lisse $\ell$-adic sheaves. The analogues
of constructible sheaves should be holonomic arithmetic $\mathcal{D}$-modules
(see \cite{bib:ber4} or \cite{bib:ber7}), but it remains to prove that they admit
Grothendieck's six operations. This fault may eventually be circumvented by the 
notion of an ``overholonomic'' arithmetic $\mathcal{D}$, as introduced by
Caro \cite{caro}.

It may also be possible to establish results on $\mathcal{D}$-modules
by proving a hypothesis of Shiho \cite[Conjecture~3.1.8]{bib:shiho},
which roughly states that
every overconvergent $F$-isocrystal on a smooth $k$-scheme $X$ admits
an extension to a log-$F$-isocrystal on some proper $k$-scheme containing
an alteration of $X$ (in the sense of de~Jong \cite{bib:dej0}).
Indeed, under this hypothesis, Shiho establishes the results of this
paper by imitating the arguments of \cite{bib:ber2}, \cite{bib:ber6}
with crystalline cohomology replaced by its logarithmic analogue.

For isocrystals on a curve, Shiho's conjecture
is de~Jong's formulation of Crew's conjecture \cite{bib:dej3},
which has been deduced from the usual formulation
by the author \cite{bib:me8}. In higher dimensions, 
it should be viewed as a globalized form of Crew's conjecture.
Shiho's conjecture is known in the unit-root case by a result of
Tsuzuki \cite{bib:tsu4}; indeed, the results of \cite{bib:tsu5} rely
on that case the same way that Shiho hopes to use the general assertion.
See \cite{bib:me14} for the first paper in a projected series addressing
Shiho's conjecture.

\subsection*{Acknowledgments}
Part of this work was carried out during a visit to the
Universit\'e de Rennes in May--June~2002; the author thanks Pierre 
Berthelot and Bas Edixhoven
for their hospitality. The author also thanks Vladimir Berkovich,
Johan de~Jong and Arthur Ogus for helpful discussions,
and the referees for numerous useful suggestions.
The author was supported by an NSF Postdoctoral Fellowship
and by NSF grant DMS-0400727.

\section{The rings of rigid cohomology}
\label{sec:prelim}

In this section, we construct the rings that occur in our
study of rigid cohomology. First, we adopt some notational conventions
for the whole paper.

For the rest of the paper, $q$ will denote a fixed power of a prime $p$,
and ``Frobenius'' will always mean the $q$-power Frobenius. Thus our
notion of an $F$-isocrystal will coincide with what other authors call
an $F^a$-isocrystal, where $a = \log_p q$.

Let $k$ be a field of characteristic $p>0$,
and let $\calO$ be a
complete discrete valuation ring of characteristic 0 with residue
field $k$ and fraction field $K$.
Let $\Gamma^*$ be the divisible closure of the value group of $K$.
Let $\pi$ be a uniformizer of $\calO$. Let $v_p$ denote
the valuation on $\calO$ normalized so that $v_p(p) = 1$, and
let $|\cdot|$ denote the corresponding norm, that is, $|x| = p^{-v_p(x)}$.
We will always assume extant and
chosen a ring endomorphism $\sigma_K$ of $K$ preserving $\calO$ and
lifting the $q$-power Frobenius on $k$. For instance, if $k$ is perfect and
$K = \Frac W(k)$,
$\sigma_K$ exists, is unique, and coincides with a power of the
Witt vector Frobenius.

Recall that the \emph{direct (inductive) limit topology}
on the union $U$ of a direct system of topological spaces
is the topology
in which a subset of $U$ is open if and only if its preimage in each $U_i$
is open. In other words, it is the finest topology for which the
maps $U_i \to U$ are continuous.
If $U$ is Hausdorff, then any compact subset of $U$ is contained
in some finite union of the $U_i$ 
(exercise, or see \cite[Lemma~15.10]{bib:gray}); in particular,
a sequence $\{x_n\}$ converges to $x$ if and only
if it does so within some finite union of the $U_i$.

\subsection{Affinoid algebras}

We recall the basic properties of affinoid algebras following 
\cite[Chapter~6]{bib:bgr}.
For an $n$-tuple $I = (i_1, \dots, i_n)$ of nonnegative
integers, write $x^I$ for $x_1^{i_1}
\cdots x_n^{i_n}$ and $|I| = i_1 + \cdots + i_n$.

\begin{defn} \label{defn:tn}
The \emph{Tate algebra} (over $K$) in the variables $x_1, \dots, x_n$
is defined as
\[
T_n = \left\{\sum_I a_I x^I: a_I \in K, \lim_{|I| \to \infty} |a_I| = 0
\right\};
\]
this ring is also denoted more explicitly $K \langle x_1, \dots, x_n 
\rangle$. It is a noetherian ring \cite[Theorem~5.2.6/1]{bib:bgr}.
The elements of $T_n$ are also referred to as \emph{strictly convergent
power series} in the $x_i$. 
\end{defn}

\begin{defn}
The \emph{Gauss norm} on $T_n$ is given by the formula
\[
\left|\sum_I a_I x^I \right| = \max_I \{ |a_I|\}.
\]
Under the Gauss norm, $T_n$ is complete, 
and every ideal of $T_n$ is closed \cite[Corollary~5.2.7/2]{bib:bgr}.
\end{defn}

\begin{defn}
To any quotient $T_n/\ida$ of $T_n$ we can associate the \emph{residue norm}:
\[
|f + \ida| = \inf\{|f+a|: a \in \ida\},
\]
which induces the quotient topology on $T_n/\ida$.
Beware that this norm need not be power-multiplicative, e.g., if
$T_n/\ida$ has nilpotent elements.
\end{defn}

\begin{defn}
Any topological $K$-algebra isomorphic to a 
quotient of a Tate
algebra is called an \emph{affinoid algebra} (over $K$).
It turns out that every $K$-algebra homomorphism between affinoid algebras
is continuous \cite[Theorem~6.1.3/1]{bib:bgr}, so the topology of
an affinoid algebra is uniquely determined by the underlying $K$-algebra.
We call that topology the \emph{intrinsic topology}.
Similarly, every finitely generated module over an affinoid algebra $A$
admits a unique topology, and all $A$-module homomorphisms between
finitely generated modules are continuous \cite[Proposition~3.7.3/3]{bib:bgr}.
\end{defn}

\begin{defn}
For $A$ an affinoid algebra, let $\Maxspec A$ be the set of maximal ideals of
$A$. For any $\idm \in \Maxspec A$, 
the field $A/\idm$ is a finite extension of $K$
\cite[Corollary~6.1.2/3]{bib:bgr}, 
so there is a unique norm $|\cdot|_{\idm}$ on $A/\idm$ extending
the norm on $K$.
\end{defn}

\begin{defn}
For $A$ an affinoid algebra, we say $A$ is \emph{reduced/normal/regular} at 
a point $\idm \in \Maxspec A$
if the local ring of $A$ at $\idm$ is reduced/normal/regular. We may drop
``at $\idm$'' if the property holds at all points, as in that case $A$
is reduced/normal/regular in the ring-theoretic sense 
\cite[Corollary~7.3.2/9]{bib:bgr}.
\end{defn}

\begin{defn}
For $A$ an affinoid algebra,
we define the \emph{spectral seminorm} $|\cdot|_{\sup,A}$ on $A$ by
\[
|x|_{\sup,A} = \sup_{\idm \in \Maxspec A} |x|_{\idm}
\]
and the \emph{spectral valuation} $v_A$ by 
\[
v_A(x) = -\log_p \left(|x|_{\sup,A}\right).
\]
We drop $A$ from the notation when it is clear from context.
As catalogued in \cite[Section~6.2]{bib:bgr},
the basic properties of the spectral seminorm include:
\begin{itemize}
\item
The spectral seminorm is a $K$-algebra seminorm:
that is, $|x+y|_{\sup} \leq \max\{|x|_{\sup}, |y|_{\sup}\}$, $|xy|_{\sup}
\leq |x|_{\sup} |y|_{\sup}$, and $|\lambda|_{\sup} = |\lambda|$ for all
$\lambda \in K$.
\item
The spectral seminorm is power-multiplicative: $|x|_{\sup}^n = |x^n|_{\sup}$ for
$n \in \NN$.
\item
Given any residue norm $|\cdot|$ on $A$, the spectral seminorm can be computed
as
\[
|x|_{\sup} = \lim_{n \to \infty} |x^n|^{1/n}.
\]
\item
The maximum modulus principle holds for $|\cdot|_{\sup}$, that is, for any
$x \in A$, there
exists $\idm \in \Maxspec A$ such that $|x|_{\idm} = |x|_{\sup}$. In particular,
the values of $|\cdot|_{\sup}$ all lie in $\Gamma^*$.
\item
If $f: A \to B$ is a $K$-algebra homomorphism between affinoid algebras,
then $|f(a)|_{\sup,B} \leq |a|_{\sup,A}$ for all $a \in A$.
\item
The spectral seminorm is a norm (that is, $|x|_{\sup} = 0$ if and only if $x=0$)
if and only if $A$ is reduced. In this case we also refer to it as the
\emph{spectral norm}.
\item
The spectral norm on $T_n$ coincides with the Gauss norm.
\end{itemize}
\end{defn}

\begin{defn}
For $A$ an affinoid algebra, we say that $x \in A$ is 
\emph{integral} if $|x|_{\sup,A} \leq 1$; the integral elements
of $A$ form a subring of $A$, denoted $A^{\inte}$.
We say that $x \in A$ is \emph{topologically nilpotent} if 
$|x|_{\sup,A} < 1$, or equivalently
(by \cite[Proposition~6.2.3/2]{bib:bgr}) if $x^n \to 0$
as $n \to \infty$.
The topologically nilpotent elements form an ideal of the integral
subring; the quotient by this ideal is called the
\emph{reduction} of $A$. The reduction of $A$ is a finitely
generated $k$-algebra \cite[Corollary~6.3.4/3]{bib:bgr};
its $\Spec$ is called the \emph{special fibre} of $A$.
Note that $A^{\inte}/\pi A^{\inte}$ is reduced if and only if it
coincides with the reduction of $A$.
\end{defn}

\begin{defn}
For $A$ an affinoid algebra and $f \in A$, the
\emph{localization} of $A$ at $f$ is a minimal affinoid algebra
containing $A$ in which $f$ is invertible. The localization is unique
up to canonical isomorphism; for $A = T_n/\ida$, it is given by
\[
A \langle f^{-1} \rangle = T_{n+1}/(\ida T_{n+1} + (x_{n+1} f - 1)).
\]
\end{defn}

\begin{defn}
For $A_1$ and $A_2$ affinoid algebras, let
$A_1 \hattimes A_2$ denote the completion of the usual tensor product.
Then $A_1 \hattimes A_2$ is again an affinoid algebra, and it is
in fact a product of $A_1$ and $A_2$ in the
category of affinoid algebras and $K$-algebra morphisms.
\end{defn}

\subsection{Dagger algebras}

Next, we recall the basic properties of dagger algebras, as described
in \cite[Section~1]{bib:gk}.
\begin{defn} \label{defn:wn}
For each $\rho > 0$, define
\[
T_{n,\rho} 
= \left\{\sum_I a_I x^I: a_I \in K, \lim_{|I| \to \infty} |a_I|\rho^{|I|} 
= 0 \right\};
\]
then $T_{n,\rho}$ is an affinoid algebra for any $\rho \in \Gamma^*$
\cite[Theorem~6.1.5/4]{bib:bgr}.
Define the \emph{Monsky-Washnitzer algebra} (over $K$) in $x_1, \dots, x_n$
as the union of the $T_{n,\rho}$ over all $\rho > 1$; this ring is
denoted $W_n$ or more explicitly
$K \langle x_1, \dots, x_n \rangle^\dagger$. 
The elements of $W_n$
are also referred to as \emph{overconvergent power series} in the $x_i$.
\end{defn}

By topologizing $W_n$ using the Gauss norm, we may view it as a dense
subring of $T_n$.
The ring $W_n$ is noetherian by a theorem of Fulton
\cite{bib:fulton},
and every ideal of $W_n$ is 
closed by the same proof as for $T_n$ (as pointed out in
\cite[Section~1.4]{bib:gk}).

\begin{defn}
Define a \emph{dagger algebra} to be any topological $K$-algebra
isomorphic to a quotient of some $W_n$.
Again, any $K$-algebra homomorphism between dagger algebras is
automatically continuous \cite[Proposition~1.6]{bib:gk}, so the topology
is uniquely determined by the underlying $K$-algebra.
We call this topology the \emph{affinoid topology}, in order to distinguish
it from a second topology on dagger algebras to be introduced in the next
section.
\end{defn}

All of the definitions of the previous section can be carried over to the
category of dagger algebras. For instance, the \emph{localization} of a
dagger algebra $A = W_n/\ida$ at an element $f$ is given by
\[
A \langle f^{-1} \rangle^\dagger = W_{n+1}/(\ida W_{n+1} + (x_{n+1} f - 1)).
\]

For $A$ a dagger algebra, the completion $\widehat{A}$ of 
$A$ with respect to the
ideal of topologically nilpotent elements is an affinoid algebra: if
$f: W_n \to A$ is a surjection with kernel $\ida$, 
then the completion of $A$
is isomorphic to $T_n/\ida T_n$.
Write $A^{\inte} = A \cap (\widehat{A})^{\inte}$.
From \cite[Theorem~1.7]{bib:gk}, we have the following relationships
between a dagger algebra $A$ and its completion $\widehat{A}$.
\begin{itemize}
\item $\widehat{A}$ is faithfully flat over $A$.
\item $\Maxspec A = \Maxspec \widehat{A}$.
\item If $A$ is reduced (resp.\ normal, regular), so is $\widehat{A}$.
\item If $A$ is reduced, its affinoid topology is induced by the
spectral norm on $\widehat{A}$.
\end{itemize}
Thanks to the faithful flatness, every finitely generated $A$-module $M$
acquires a canonical ``affinoid'' topology, namely the subspace topology 
induced by
the inclusion of $M$ into $M \otimes_A \widehat{A}$.
Again, all homomorphisms between finitely
generated $A$-modules are continuous for this topology.

\begin{lemma} \label{lem:flat}
Let $A$ and $B$ be dagger algebras such that
$A^{\inte}/\pi A^{\inte}$ and $B^{\inte}/\pi B^{\inte}$ are reduced.
Let $f: A \to B$ be a $K$-algebra homomorphism of dagger
algebras. If the induced map $R \to S$ on reductions is flat, then
$f$ is flat, as is the induced map $A^{\inte} \to B^{\inte}$.
\end{lemma}
\begin{proof}
Since $A^{\inte}/\pi A^{\inte}$ and $B^{\inte}/\pi B^{\inte}$ are reduced,
we have
$R = A^{\inte}/\pi A^{\inte}$ and $S = B^{\inte}/\pi B^{\inte}$,
and $A^{\inte}$ and $B^{\inte}$ are both $\pi$-adically separated.
Since the exact sequence $0 \to \pi A^{\inte} \to A^{\inte} \to R \to 0$
remains exact upon tensoring with $B^{\inte}$ over $A^{\inte}$, 
$\Tor_1^{A^{\inte}}(R, B^{\inte})$
vanishes. Thus \cite[Lemma~2.1]{bib:mw} implies that the map
$A^{\inte} \to B^{\inte}$ is flat, which implies in turn that $f$ is flat.
\end{proof}

\begin{remark}
Although a dagger algebra $A$ is not complete with respect to the affinoid
topology, it is ``weakly complete'' 
(in the terminology of \cite{bib:mw})
in the sense that for any $y_1, \dots, y_n \in A^{\inte}$,
the homomorphism
$K[x_1, \dots, x_n] \to A$ sending $x_i$ to $y_i$ extends to a continuous
homomorphism $W_n \to A$ \cite[1.9]{bib:gk}.
For example, the series $(1-a)^{-1} = \sum_{n=0}^\infty a^n$ converges whenever
$a$ is topologically nilpotent; hence any element of $A^{\inte}$
whose reduction is 1 is a unit.
\end{remark}

\subsection{Fringe algebras}

\begin{defn}
For $A$ a dagger algebra, we define a \emph{fringe algebra} to be any
algebra $B$ of the form $f(T_{n,\rho})$ for some 
surjection $f: W_n \to A$ and some $\rho \in \Gamma^*$ with $\rho > 1$.
(The name can be justified by visualizing $\Maxspec B$ as being $\Maxspec A$ 
plus some inessential ``fringe'' along its boundary.)
\end{defn}

Because of the restriction on $\rho$, $B$ is an affinoid 
algebra, and so has an intrinsic topology. The intrinsic topology
is finer than that induced from the affinoid topology of $A$; in 
particular, the latter is not complete.

The key technical tool for handling fringe algebras is the following
lemma, for whose proof see \cite[Lemma~1.8]{bib:gk}.
\begin{lemma} \label{lem:fringe}
Let $g: A_1 \to A_2$ be a morphism of dagger algebras, and choose
surjections $f_i: W_{n_i} \to A_i$. Then for any $\rho_1 >1$,
there exists $\rho_2 > 1$ such that 
$g(f_1(T_{n_1,\rho_1})) \subseteq f_2(T_{n_2,\rho_2})$.
\end{lemma}
\begin{cor} \label{cor:fringecon}
Let $g: A_1 \to A_2$ be a $K$-algebra homomorphism between dagger algebras.
Then any fringe algebra of $A_1$ maps under $g$ into a fringe algebra of 
$A_2$.
\end{cor}
\begin{cor} \label{cor:fringe2}
For any fringe algebras $B_1, \dots, B_m$ of $A$, there is a fringe
algebra $B$ containing $B_1, \dots, B_m$.
\end{cor}
\begin{proof}
Choose a surjection $f: W_n \to A$. Then Lemma~\ref{lem:fringe} implies
that $B_i \subseteq f(T_{n,\rho})$ for $\rho>1$ sufficiently small.
\end{proof}
\begin{cor} \label{cor:fringe3}
Let $A$ be a dagger algebra and let $x \in A$ be a topologically nilpotent
element for the affinoid topology. Then $x$ is topologically nilpotent in
some fringe algebra $B$.
\end{cor}
\begin{proof}
Since $x$ is topologically nilpotent, we can find an integer $n$ and
an element $\lambda$ of $K$ with $|\lambda| > 1$ such that
$x^n \lambda \in A^{\inte}$. By the weak completeness
of dagger algebras, there is a $K$-algebra homomorphism from
$W_1$ to $A$ taking $x_1$ to $x^n \lambda$. By
Corollary~\ref{cor:fringecon}, for any $\rho < |\lambda|$,
the fringe algebra $T_1(\rho)$ maps into a fringe algebra $B$ of $A$;
since $x_1 \lambda$ is topologically nilpotent in $T_1(\rho)$,
its image $x^n$ is topologically nilpotent in $B$, as then is $x$.
\end{proof}

\begin{remark}
The converse of Corollary~\ref{cor:fringe3} does not hold. For instance,
if $A = W_n$, then $x_1$ is topologically nilpotent in each fringe algebra
of $A$ but not in the affinoid topology of $A$.
\end{remark}

\begin{defn}
For $A$ a dagger algebra, define the \emph{fringe topology} of $A$ to
be the direct limit of the intrinsic topologies on its fringe algebras.
(This topology is the one called the ``direct limit topology''
in \cite[4.2]{bib:gk}.) Note
that the fringe topology is Hausdorff, since it is finer than the affinoid
topology. Also, by Corollary~\ref{cor:fringecon},
any $K$-algebra homomorphism between dagger algebras is
automatically continuous for the fringe topologies (as well as for the affinoid
topologies). Moreover, any finitely generated $A$-module also acquires
a fringe topology (by base extension from a suitable fringe algebra),
and homomorphisms between finitely generated $A$-modules are automatically
continuous for these fringe topologies (as well as for the affinoid
topologies).
\end{defn}

\begin{defn}
For $A_1, A_2$ dagger algebras,
we define the completed tensor product $A_1 \hattimes A_2$
as the direct limit of $B_1 \hattimes B_2$ as each $B_i$ runs over all
fringe algebras of $A_i$. 
Then $A_1 \hattimes A_2$ is again a
dagger algebra, and (thanks to Lemma~\ref{lem:fringe} and its corollaries)
is a product of $A$ and $B$ in the category
of dagger algebras and (continuous) $K$-algebra morphisms. For example,
\[
K \langle x_1, \dots, x_n \rangle^\dagger
\hattimes K \langle y_1, \dots, y_n \rangle^\dagger
\cong K \langle x_1, \dots, x_n, y_1, \dots, y_n \rangle^\dagger.
\]
For $A$ a dagger algebra, we write
\[
A \langle x_1, \dots, x_n \rangle^\dagger
= A \hattimes K \langle x_1, \dots, x_n \rangle^\dagger.
\]
\end{defn}
\begin{remark} \label{rem:flatprod}
If $A_1^{\inte}/\pi A_1^{\inte}$ and $A_2^{\inte}/\pi A_2^{\inte}$ are reduced,
then $A_1 \hattimes A_2$ is flat over $A_1$ and $A_2$ by
Lemma~\ref{lem:flat}. Namely, if $R_1$ and $R_2$ are the corresponding
reductions, then $R_1$ is vacuously flat over $k$ and so $R_1 \otimes_k
R_2$ is flat over $R_2$, and vice versa.
\end{remark}

\subsection{Interlude: a Gr\"obner basis calculation}

The following lemma allows us to relate the spectral norm of a fringe
algebra with that of its underlying dagger algebra.
The proof is essentially a construction of
a Gr\"obner basis; see \cite{bib:eisenbud} for an introduction to this
circle of ideas.
\begin{lemma} \label{lem:grobner}
Let $\ida$ be an ideal of $W_n$. Then there exists $\rho_0 > 1$ such that
whenever $\rho \in (1, \rho_0] \cap \Gamma^*$, 
and $y \in W_n$ and  $z \in T_{n,\rho}$ satisfy $y - z \in \ida$,
we can find $u \in T_{n,\rho}$ with
\[
u-z \in \ida, \quad
|u| \leq |y|, \quad
|u|_{\rho} \leq
|z|_{\rho},
\] 
where $|\cdot|_\rho$ is the spectral norm on $T_{n,\rho}$.
\end{lemma}
We note in passing that the argument also yields
an alternate proof that $W_n$ and its integral subring are noetherian.
\begin{proof}
Define the partial ordering $\preceq$ on $\ZZ^n_{\geq 0}$ so that
$I \preceq J$ when each component of $I$ is less than or equal to
the corresponding component of $J$. Then $\preceq$ is a well partial
ordering (wpo), i.e., 
any infinite sequence contains an nondecreasing subsequence.
(To see this, first pass to a subsequence nondecreasing in the first component,
then also in the second component and so on.)

Choose a total ordering $\leq$ of $\ZZ^n_{\geq 0}$ extending both $\preceq$ 
and the quasi-ordering by total degree, which is also 
compatible with addition (that is, $I\leq J$ if and only if $I+K \leq J+K$).
For example, one may choose the ``deglex'' ordering, in which $I \leq J$
if either $I = J$, $|I| < |J|$, or $|I| = |J|$ and 
in the first position where $I$ and $J$ differ, $I$ has
the lesser entry. Since $\preceq$ is a well partial ordering,
$\leq$ is a well ordering.
For $a = \sum a_I x^I \in T_{n,\rho}$,
define the \emph{$\rho$-leading term} of $a$ to be the expression
$a_I x^I$, for $I$ the largest tuple under $\leq$
which maximizes $|a_I x^I|_\rho = |a_I| \rho^{|I|}$. (There is a largest such $I$
because there are only finitely many such $I$.)

We claim that for each $a \in W_n$, 
the $1$-leading term of $a$ coincides with the
$\rho$-leading term of $a$ for each sufficiently small $\rho > 1$ (where
the sense of ``sufficiently small'' depends on $a$).
To see this, let $a_I x^I$ be the 1-leading
 term of $a$.
For each tuple $J = (j_1, \dots, j_n)$, we then have either 
\begin{enumerate}
\item[(a)] $|a_J| < |a_I|$, or 
\item[(b)]
$|a_J| = |a_I|$ and $J \leq I$; in this case we have $|J| \leq |I|$ since
$\leq$ refines the quasi-ordering by total degree.
\end{enumerate}
Pick $\eta>1$ for which $a \in T_{n,\eta}$.
If $|a_J x^J|_\eta \leq |a_I x^I|_\eta$, then in case (a),
we have $|a_J x^J|_\rho < |a_I x^I|_\rho$ for all
$\rho \in [1, \eta) \cap \Gamma^*$;
in case (b), we have $|a_J x^J|_\rho \leq |a_I x^I|_\rho$ and $J \leq I$.
So these terms are all okay for any $\rho$; in fact, because $a \in T_{n,\eta}$,
there are only finitely many tuples $J$ with
$|a_J x^J|_\eta > |a_I x^I|_\eta$.
For each such $J$, we must be in case (a), so $|y_J x^J|_\rho 
< |y_I x^I|_\rho$ for $\rho \in (1,\eta]$
sufficiently small. This yields the claim.

Define elements $a_1, a_2, \dots$ of $\ida$ as follows.
Given $a_1, \dots, a_{i-1}$, choose $a_i$ if possible to be an element
of $\ida$ whose $1$-leading term is not a multiple of the
1-leading term of $a_j$ for any $j<i$,
otherwise stop. Since $\preceq$ is a well partial ordering,
this process must eventually stop; at that point, every 1-leading
term of every element of $\ida$ occurs as a multiple of the
1-leading term for some $a_i$.

Let $A$ be the (finite) set of the $a_i$.
By the argument above, we can choose $\rho_0 > 1$ such that
for each $a \in A$, $a \in T_{n,\rho_0}$ and 
the $\rho_0$-leading term and 1-leading term of $a$ are equal.
We claim that this choice of $\rho_0$ works.
To see this, given $\rho \in (1, \rho_0] \cap \Gamma^*$, and
$y \in W_n$ and $z \in T_{n,\rho}$ with $y - z \in \ida$
(and $y \neq 0$, or else $u=0$ trivially works),
we construct a sequence $\{c_j\}$ of elements of $T_{n,\rho_0}$
and a sequence $\{d_j\}$ of elements of $A$,
as follows. Given the sequences up to $c_j$ and $d_j$, 
put $z_j = z - c_1 d_1 - \cdots - c_j d_j$ (where if $j=0$ we take
$z_0 = z$).
If $|z_j| \leq |y|$, then stop. Otherwise,
let $e_I x^I$
be the leading term of $z_j - y$.
By the construction of the $a_i$, we can find a monomial $c_{j+1}$
and an element $d_{j+1} \in A$ so that $c_{j+1} d_{j+1}$ has
1-leading term, and hence $\rho$-leading term, equal to $e_I x^I$.

From the construction, we clearly have $|z_j|_\rho \leq |z|_\rho$.
On the other hand, if the process were never to terminate, we would have
$|z_j| \to 0$ as $j \to \infty$: if $|z_j| = |z_{j+1}|$, the 1-leading term
of $z_{j+1}$ must be smaller under $\leq$ than that of $z_j$, and
well-orderedness says this cannot continue forever, so eventually
$|z_j|$ drops (by a discrete factor). Since $y \neq 0$, this yields
a contradiction. Thus the process terminates at some $z_j$, and we
may take $u=z_j$.
\end{proof}

Our principal application of Lemma~\ref{lem:grobner} is the following
result, reminiscent of the Hadamard three circles theorem.
\begin{prop} \label{prop:hadamard}
Let $A$ be a dagger algebra, and let $B$ be a fringe
algebra of $A$. Then
for any rational $\epsilon \in (0,1]$, there exists a fringe algebra $C$ of $A$
containing $B$ such that
\[
|x|_{\sup,C} \leq |x|_{\sup,A}^{1-\epsilon} |x|_{\sup,B}^{\epsilon}
\]
for all $x \in B$.
\end{prop}
\begin{proof}
Choose a surjection $f: W_n \to A$, and let $|\cdot|_A$ denote the
induced residue norm on $A$. By Corollary~\ref{cor:fringe2},
we can choose $\rho>1$ such that
$f(T_{n,\rho})$ contains $B$; by shrinking $\rho$ further, we can
ensure that $\rho \leq \rho_0$ for some $\rho_0$ as in 
Lemma~\ref{lem:grobner}.
Put $B' = f(T_{n,\rho})$, and let $|\cdot|_{B'}$ denote the 
residue norm on $B'$ induced by the spectral norm on $T_{n,\rho}$.

Let $|\cdot|_\eta$ denote the spectral norm on $T_{n,\eta}$.
For $y \in T_{n,\rho}$, we have the inequality
\[
|y|_{\rho^\epsilon} \leq |y|^{1-\epsilon} |y|_{\rho}^{\epsilon}
\]
since it is actually an equality for monomials.

Given $x \in B$, choose $y \in W_n$ and $z \in T_{n,\rho}$,
both reducing to $x$ in $A$, with $|y| = |x|_A$ and $|z|_\rho = |x|_B$. By
Lemma~\ref{lem:grobner}, we can find $u \in T_{n,\rho}$,
also reducing to $x$ in $A$,
with $|u| \leq |y|$ and $|u|_\rho \leq |z|_\rho$.
We then have
\[
|u|_{\rho^\epsilon} \leq |u|^{1-\epsilon} |u|_{\rho}^{\epsilon}.
\]
Put $C = f(T_{n,\rho^\epsilon})$; 
then the residue norm $|\cdot|_C$ induced by the spectral norm on
$T_{n,\rho^\epsilon}$ satisfies
\[
|x|_C \leq |x|_A^{1-\epsilon} |x|_{B'}^{\epsilon}
\]
for all $x \in B$.
Replacing $x$ by $x^i$, taking $i$-th roots and taking limits as
$i \to \infty$ yields
\[
|x|_{\sup,C} \leq |x|_{\sup,A}^{1-\epsilon} |x|_{\sup,B'}^{\epsilon}.
\]
Since $|x|_{\sup,B'} \leq |x|_{\sup,B}$, this yields the desired result.
\end{proof}

\subsection{Robba rings}
\label{subsec:robba}

The Robba ring is typically used to study $p$-adic differential
equations over a field. We will generalize its definition so that it
can be constructed not just over a field, but over a reduced 
affinoid or dagger
algebra. (It may be possible to avoid the reducedness hypothesis, but 
we have no need to do so here.)

\begin{defn}
Let $A$ be a reduced 
affinoid algebra.
For each $r>0$, let $\calR_{A,r}$ denote the set of
formal (bidirectional)
power series $x = \sum_{i \in \ZZ} x_i t^i$, with $x_i \in A$,
such that
\[
\lim_{i \to -\infty} v_A(x_i) + ri = \infty,
\qquad
\liminf_{i \to +\infty} \frac{v_A(x_i)}{i} \geq 0;
\]
an equivalent condition is that $\lim_{i \to \pm \infty} v_A(x_i)+si = \infty$
for $0 < s \leq r$.
This set forms a ring under the usual multiplication law for power series.
For $0 < s \leq r$, we define the valuation $w_{A,s}$ on $\calR_{A,r}$ by the
formula
\[
w_{A,s}(x) = \min_i \{ v_A(x_i) + si \};
\]
then $\calR_{A,r}$ is complete for the Fr\'echet topology defined by these
valuations (that is, any sequence which is Cauchy with respect to each 
$w_{A,s}$ is convergent).
The union $\calR_A = \cup_{r>0} \calR_{A,r}$ is called the
\emph{Robba ring} of $A$; we equip it with the direct limit of the Fr\'echet
topologies on the $\calR_{A,r}$.
\end{defn}

\begin{remark} \label{rem:notnoeth}
Robba rings are not as well-behaved as the other rings we have considered
so far; for instance, the Robba ring over $K$ is not noetherian.
By a theorem of Lazard \cite{bib:laz},
it does however have the B\'ezout property: 
every finitely generated ideal over $\calR_K$ (or $\calR_{K,r}$) is principal.
Over a B\'ezout ring, any finitely presented projective module, or any
finitely generated locally free module, is actually free.
\end{remark}

Extending this definition to the dagger case is a bit tricky; the
most straightforward extension is not the one that we really need
(see Remark~\ref{rem:notrobba}).

\begin{defn} \label{defn:dagrobba}
Let $A$ be a reduced dagger algebra. 
Define the \emph{Robba ring} $\calR_A$ over $A$
as the ring of formal series $x = \sum_{i \in \ZZ} x_i t^i$ over $A$
with the following property: 
for each sufficiently small $r>0$,
there exists a fringe algebra $B$ such that the $x_i$ belong to $B$ and
\[
\lim_{i \to \pm \infty} v_B(x_i) + ri = \infty,
\]
or equivalently, $x_i p^{\lfloor ri \rfloor}$ converges to zero as
$i \to \pm \infty$ in the fringe topology of $A$.
\end{defn}
In practice, it is sometimes easier to check the following condition for
membership in the Robba ring of a dagger algebra.
\begin{prop} \label{prop:robcrit}
Let $A$ be a reduced dagger algebra, and let $\widehat{A}$ be the completion
of $A$ for the affinoid topology.
Suppose $x = \sum_i x_i t^i \in \calR_{\widehat{A},r}$ satisfies $x_i \in A$
for all $i$. If the condition that $x_i p^{\lfloor si \rfloor} \to 0$
as $i \to \pm \infty$
in the fringe topology of $A$ holds for $s=r$, it holds also when $0<s < r$.
\end{prop}
\begin{proof}
Let $B$ be a fringe algebra in which $x_i p^{\lfloor ri \rfloor} \to
0$ as $i \to \pm \infty$.
Then we also have $x_i p^{\lfloor si \rfloor} \to 0$ in $B$ as $i \to -\infty$.
On the other hand, by Proposition~\ref{prop:hadamard},
for any rational $\epsilon \in (0,1]$, 
we can find a fringe algebra $C$ such that
\[
|x_i|_{\sup,C} \leq |x_i|^{1-\epsilon}_{\sup,A} |x_i|^{\epsilon}_{\sup,B}
\qquad \mbox{for all $i$.}
\]
If we choose $\epsilon < s/r$, we may then write
\[
p^{- \lfloor si \rfloor}  |x_i|_{\sup,C} 
\leq
\left( p^{-\lfloor (s-\epsilon r)i \rfloor} 
|x_i|^{1-\epsilon}_{\sup,A} 
\right)
\left(
 p^{-\lfloor \epsilon ri \rfloor} |x_i|^{\epsilon}_{\sup,B}
 \right);
\]
the right side tends to 0 as $i \to \infty$, so the left side does as well.
This yields the desired result.
\end{proof}
\begin{cor}
Let $A$ be a reduced dagger algebra, and let $\widehat{A}$ be the completion
of $A$ for the affinoid topology.
Let $x = \sum_{i \in \ZZ} x_i t^i$ be a formal series with coefficients in $A$.
Then $x \in \calR_A$ if and only if $x \in \calR_{\widehat{A}}$ and
for some $r>0$, $x_i p^{\lfloor ri \rfloor} \to 0$ as $i \to \pm \infty$
in the fringe topology of $A$.
\end{cor}
This means checking convergence in a Robba ring over a dagger algebra
is typically a two-step process: one first checks convergence
over the completion,
then one checks the condition on the fringe topology.

\begin{remark}
We mention an example in passing (that will not be used later).
For $A = K \langle x \rangle^\dagger$, it can be shown that
a series $\sum_{i \in \ZZ}
a_i t^i$ belongs to $\calR_A$ if and only if there exist some $\delta>1$
and $\epsilon>0$ such that the series converges for 
$|x| \leq \delta$, $\delta^{-1} \leq |t| < 1$ and $|x^\epsilon t| < 1$.
\end{remark}

\begin{remark} \label{rem:notrobba}
Do not confuse $\calR_A$ with the
the direct limit of $\calR_B$ over the
fringe algebras $B$ of $A$, where $\calR_B$ is defined by viewing $B$
as an affinoid. (The latter is obtained by reversing
the order of the quantifiers
``for $r>0$ sufficiently small'' and ``there exists a fringe algebra $B$''
in Definition~\ref{defn:dagrobba}.)
This ring is strictly smaller than $\calR_A$: for instance,
for $A = K \langle x \rangle^\dagger$,
$1-xt$ is a unit in $\calR_A$ but not in $\calR_B$ for any fringe algebra $B$.
\end{remark}

\begin{defn} \label{defn:robbatop}
For $A$ a reduced dagger algebra,
we define $w_r$ on $\calR_A$ (or rather, on the subring where it makes sense)
by restriction from $\calR_{\widehat{A}}$.
For $B$ a fringe algebra of $A$,
we also define $w_{B,r}$ on the subring where it makes sense.
We topologize $\calR_A$ as follows: the sequence $\{x_n\}$
converges if and only if the
$w_{B,r}(x_n)$ are defined and tend to $\infty$ for some $B$ and $r$,
and for sufficiently small $s$ (not depending on $n$), the $w_s(x_n)$
are defined and tend to $\infty$. Equivalently (as in the proof of
Proposition~\ref{prop:robcrit}), for each sufficiently small $r$,
there exists a fringe algebra $B$ such that the $w_{B,r}(x_n)$
are defined and tend to $\infty$.
\end{defn}

\begin{defn}
In the notation of a Robba ring, we use the superscript ``$+$''
to denote the subring consisting
of power series with only nonnegative powers of $t$.
We call this subring the \emph{plus part}.
\end{defn}

\begin{defn}
In the notation of a Robba ring (or plus part), we use the superscript
``$\inte$'' to denote the subring consisting of series with
integral coefficients.
We call this subring the \emph{integral subring}.
Note that this subring can be topologically characterized, as the set of
$x$ for which the sequence $\{(cx^a)^n\}_{n=1}^\infty$ converges to
zero for any $c \in K$ with $|c| < 1$ and any positive integer $a$.
Consequently, any continuous map between Robba rings induces a continuous
map on integral subrings.
Again, the subset of the integral subring consisting of topologically
nilpotent elements is an ideal; we call the quotient by this ideal the
\emph{reduction} of the Robba ring. It is isomorphic to the Laurent series
ring $R((t))$, where $R$ is the reduction of $A$.
\end{defn}

We will also need a multidimensional version of this construction.
\begin{defn}
Let $T = \{t_1, \dots, t_n\}$ 
be a finite set and $A$ a reduced affinoid algebra.
We then define $\calR^T_{A,r}$ as the set of formal sums
$\sum_I x_I t^I$, for $I$ ranging over $\ZZ^n$, such that
for $0 < s_1, \dots, s_n \leq r$,
\[
\lim_I v_A(x_I) + s_1i_1 + \cdots + s_n i_n = \infty.
\]
(That is, the quantity on the left takes values below any given cutoff
for only finitely many $I$.)
We define $\calR^T_A$ as the union of the $\calR^T_{A,r}$. 
(If $T$ is empty, we put $\calR^T_A = A$.) Again, each $\calR^T_{A,r}$
has a natural Fr\'echet topology, and we topologize $\calR^T_A$ as the direct
limit thereof.
To define $\calR^T_A$ for $A$ a reduced dagger algebra, we require that for
each tuple of $s_1, \dots, s_n>0$ 
with $\max_j\{s_j\}$ sufficiently small, we have
\[
\lim_I v_B(x_I) + s_1i_1 + \cdots + s_n i_n = \infty
\]
for some fringe algebra $B$; again, for series in $\calR^T_{\widehat{A}}$
it is enough to check this for a single tuple. We topologize as in
Definition~\ref{defn:robbatop}.
\end{defn}

One immediate benefit of the multidimensional construction 
(besides its ultimate value in computing cohomology with compact supports
in Chapter~\ref{sec:push2}) is that it allows us to
interpret products of Robba rings. That is, the category of topological rings 
isomorphic to $\calR^S_A$ for some finite set $S$ and some reduced affinoid
(resp.\ dagger) algebra $A$, with continuous $K$-algebra morphisms,
has pairwise products given by
\[
\calR^{S}_A \hattimes \calR^{T}_B = \calR^{S \cup T}_{A \hattimes B}.
\]

\begin{remark}
Although Robba rings are never noetherian (see Remark~\ref{rem:notnoeth}), 
as noted above, their integral
subrings are better behaved. For instance, $\calR^{\inte}_K$ 
is a discrete valuation
ring with residue field $k((t))$. More generally, we expect that 
$\calR^{T,\inte}_A$ is always noetherian for $A$ a reduced affinoid or
dagger algebra, though this assertion is a bit subtle in case $T$
has more than one element.
For instance, if $T = \{t,u\}$, then
$\calR^{T,\inte}_{K}/\pi \calR^{T,\inte}_K$ consists of formal Laurent
series $\sum_{i,j \in \ZZ} c_{i,j} t^i u^j$ such that for each $a,b > 0$,
there exist only finitely many pairs $(i,j)$ with $c_{i,j} \neq 0$
and $i + aj < b$. This ring includes the usual double Laurent series ring
$k \llbracket t,u \rrbracket[t^{-1},u^{-1}]$
but is somewhat bigger, and its noetherianness is not immediately evident.
\end{remark}

\begin{lemma} \label{lem:nilpotent}
Let $A$ be a reduced dagger algebra. Suppose that $x \in \calR_A \cap
\calR_{B,r}$ for some fringe algebra $B$ and some $r>0$, and that
$w_r(x) > 0$. Then there exists another fringe algebra $C$ of $A$ such that
$w_{C,r}(x) > 0$.
\end{lemma}
\begin{proof}
Write $x = \sum_i x_i t^i$; then already for all but finitely many
$i$, we have $w_{B,r}(x_i t^i) > 0$, and so $w_{C,r}(x_i t^i) > 0$ for
any fringe algebra $C$ containing $B$. To get the inequality for each
of the finitely many remaining $i$, apply Proposition~\ref{prop:hadamard}.
\end{proof}

\section{Frobenius and differential structures}
\label{sec:frobdiff}

Overconvergent $F$-isocrystals on smooth $k$-schemes are essentially
modules over certain dagger algebras equipped with ``Frobenius''
and ``differential'' structures.
In this chapter, we review the notions that make it possible to
describe these extra structures, culminating in the notion of
$(\sigma, \nabla)$-modules and their cohomology.

\subsection{Differentials}

\begin{defn} \label{defn:moddiff}
In order to make uniform statements, we adopt the following \emph{ad 
hoc} definition. An \emph{eligible category} is the category of one of
the following types of objects.
\begin{itemize}
\item Affinoid algebras and $K$-algebra morphisms.
\item Dagger algebras and $K$-algebra morphisms.
\item Multidimensional Robba rings over reduced affinoid algebras,
and continuous $K$-algebra morphisms.
\item Multidimensional Robba rings over reduced dagger algebras,
and continuous $K$-algebra morphisms.
\item Plus parts of the rings in the previous two categories.
\end{itemize}
\end{defn}

\begin{defn} 
Let $\calC$ be an eligible category.
Let $R$ be an object in $\calC$,
and let $J$ be the kernel of the multiplication map
$R \hattimes R \to R$ (with the subspace topology). 
Define the \emph{module of continuous differentials}
of $R$ as 
\[
\Omega^1_{R/K} = J/J^2,
\]
equipped with the quotient topology and the continuous $K$-linear derivation 
$d: R \to \Omega^1_{R/K}$ sending $x$ to $1 \otimes x - x \otimes 1$.
Then $\Omega^1_{R/K}$ is 
universal for continuous $K$-linear
derivations from $A$ into separated topological $A$-modules. 
As usual, we write
$\Omega^i_{R/K}$ for the exterior power $\wedge^i \Omega^1_{R/K}$
over $R$, and $d$ for the map $\Omega^i_{R/K} \to
\Omega^{i+1}_{R/K}$ induced by the canonical derivation $d: R \to
\Omega^1_{R/K}$.
\end{defn}

It is easily checked that $\Omega^1_{T_n/K}$ and $\Omega^1_{W_n/K}$
are freely generated by $dx_1, \dots, dx_n$. 
It follows that $\Omega^1_{A/K}$ is finitely generated whenever $A$
is an affinoid or dagger algebra; if $A$ is regular, then
$\Omega^1_{A/K}$ is also locally free, hence projective.
Similarly, one verifies that for $A$ a reduced affinoid or dagger algebra,
$\Omega^1_{\calR^T_{A}/K}$ is isomorphic to the direct sum of
$\calR^T_A \otimes_A \Omega^1_{A/K}$ with the $\calR^T_A$-module freely generated
by $T$. (For instance, $\Omega^1_{\calR_K/K}$ is the free
module generated by $dt$.) In particular, if $A$ is regular, then 
$\Omega^1_{\calR^T_{A}/K}$ is also finite locally free.


\begin{defn}
If $A$ is a subring of $B$, we define the
\emph{relative module of continuous differentials}
$\Omega^1_{B/A}$ as the quotient of $\Omega^1_{B/K}$ by the submodule
generated by $da$ for all $a \in A$. 
\end{defn}

We will be particularly interested in the relative construction when $B 
= \calR^T_A$;
in this case we can write
\[
dx = \sum_{t \in T} \frac{\partial x}{\partial t} \otimes\,dt.
\]
\begin{prop} \label{prop:deriv}
Let $A$ be an affinoid or dagger algebra.
Given $x \in \calR^T_A$, write $x = \sum_i c_i t^i$ for 
$c_i \in \calR^{T \setminus \{t\}}_A$. Then:
\begin{enumerate}
\item[(a)] The kernel of $\frac{\partial}{\partial t}$ consists of those
$x$ with $c_i = 0$ for $i \neq 0$.
\item[(b)] The image of $\frac{\partial}{\partial t}$ consists of those
$x$ with $c_{-1} = 0$.
\item[(c)] The induced map $\calR^T_A/(\ker \frac{\partial}{\partial t})
\to \im \frac{\partial}{\partial t}$ is a homeomorphism.
\end{enumerate}
\end{prop}
\begin{proof}
\begin{enumerate}
\item[(a)] Obvious.
\item[(b)] 
It suffices to check that the image
includes every $x$ for which $c_i = 0$ for $i \geq -1$,
and every $x$ for which $c_i = 0$ for $i \leq -1$.
Write $T = \{t_1, \dots, t_n\}$ with $t_1 = t$,
and write $x = \sum_I a_I t^I$. Then we are given
\[
a_I p^{\lfloor s_1i_1 + \cdots + s_n i_n \rfloor} \to 0
\qquad \forall s_1, \dots, s_n \in (0, r]
\]
and need to deduce that
\[
a_I p^{\lfloor s_1' (i_1+1) + \cdots + s_n i_n - \log_p |i_1+1| \rfloor} \to 0
\qquad \forall s_1', s_2, \dots, s_n \in (0, r']
\]
for some $r' < r$.
In the first case, we get the second estimate from the first by taking
$s_1 < s_1'$ and $r'=r$; in the second case, we get it by taking $s_1 > s_1'$
and $r' < r$.
\item[(c)]
It is clear that the induced map is continuous; the proof of
(b) shows that its inverse is also continuous.
\end{enumerate}
\end{proof}

\subsection{Algebras of MW-type and Frobenius lifts}
\label{subsec:mwtype}

\begin{defn}
We say an affinoid or dagger algebra $A$ is \emph{of MW-type}
(short for ``of Monsky-Washnitzer type'') if $A^{\inte}/\pi A^{\inte}$
is a smooth $k$-algebra. In particular, $A$ is regular, and its reduction
coincides with $A^{\inte}/\pi A^{\inte}$.
\end{defn}
In the terminology of \cite{bib:mw},
a dagger algebra of MW-type is a
``formally smooth, weakly complete, weakly finitely generated algebra''.

Algebras of MW-type are known to have good lifting properties
(see for instance \cite{bib:arabia}).
For one, given any smooth $k$-algebra $R$, there exist an 
affinoid
algebra and a dagger algebra
of MW-type with reduction $R$, and both are unique up to noncanonical
isomorphism.
For another, if $A$ and $B$ are affinoid/dagger algebras of MW-type with
reductions $R$ and $S$, then any morphism $R \to S$ lifts to a morphism
$A \to B$. These lifts retain many good properties of the original map:
e.g., they are finite flat if the original map is
(see Lemma~\ref{lem:flat}).

\begin{remark} \label{rem:noethnorm}
One example of the lifting properties is Noether normalization: any
affinoid (resp.\ dagger) algebra of MW-type whose special fibre is
purely of dimension $n$ (e.g., irreducible) can be written as
a finite extension of $T_n$ (resp. $W_n$), by applying classical
Noether normalization and lifting the resulting map from $k[x_1, \dots,
x_n]$ into the reduction. 
\end{remark}

Another key example is given by the Frobenius lifts.
\begin{defn}
For $A$ a ring in any of the categories in Definition~\ref{defn:moddiff},
define a \emph{Frobenius lift} on $A$ to be a 
continuous homomorphism $\sigma: A \to A$ which is $\sigma_K$-linear over $K$
and which acts on the
reduction of $A$ via the $q$-power Frobenius map.
If we write $A_{\sigma_K}$ for $A$ with its $K$-action funneled through 
$\sigma_K$, then a 
Frobenius lift is given by a $K$-linear map $A_{\sigma_K} \to A$ lifting the
relative Frobenius on reductions. In case $A$ is an affinoid or dagger algebra
of MW-type, the existence of a Frobenius lift is thus guaranteed by the good lifting
properties discussed above. If $A$ is a Robba ring or plus part over an affinoid
or dagger algebra $B$ of MW-type, one can extend to $A$ any Frobenius lift 
$\sigma_B$ on $B$ by setting
\[
\sum_I b_I t^I \mapsto  \sum_I b_I^{\sigma_B} t^{qi}.
\]
We call this the \emph{standard extension} of $\sigma_B$.
\end{defn}

\begin{remark} \label{rem:frob}
Note that for $A$ a dagger algebra of MW-type
and $\sigma$ a Frobenius lift on $A$,
$\sigma$ is finite flat because the Frobenius map on the reduction 
of $A$ (that being a smooth $k$-scheme) is finite flat.
\end{remark}

\begin{remark} \label{rem:frobfringe}
Beware that a Frobenius lift on a dagger algebra does not act on individual
fringe algebras. For instance, the Frobenius lift on $T_n$
given by $\sum c_I t^I \mapsto \sum c_I^{\sigma_K} t^{qI}$ carries
the fringe algebra $T_{n,\rho}$ into the larger fringe algebra
$T_{n,\rho^{1/q}}$. 
\end{remark}

By restricting to dagger algebras of MW-type, we can invoke the
following result \cite[Theorem~8.3]{bib:mw} (see also
\cite[Proposition~3.6]{bib:ber2}).
\begin{prop} \label{prop:trace}
Let $f: R \to S$ be a morphism of dagger algebras of MW-type, and suppose
$S$ is finite and locally free over $R$.
Then there is an $R$-linear
map of complexes $\Trace: \Omega^{.}_{S/K} \to \Omega^{.}_{R/K}$
such that $\Trace \circ f$ is multiplication by the degree of $S$ over $R$.
\end{prop}
Here $R$-linearity means that
\[
\Trace(\omega \wedge \eta) = \omega \wedge \Trace(\eta)
\qquad
\omega \in \Omega^{.}_{R/K}, \eta \in \Omega^{.}_{S/K}.
\]

\begin{remark} \label{rem:byhand}
We believe the same result is true in any eligible category, but the
proof of \cite[Theorem~8.3]{bib:mw} uses the noetherian property of
dagger algebras, and so does not apply to Robba rings.
(The proof of \cite[Proposition~3.6]{bib:ber2} relies on
\cite[Theorem~8.3]{bib:mw}, so is of no help either.)
However, we can construct the trace ``by hand'' in two particular instances
of interest:
\begin{enumerate}
\item[(a)] In case $R = \calR_A$ and $S = \calR_B$ where $A$ and $B$ are
 reduced
dagger
algebras of MW-type
and $B/A$ is finite projective. In this case we may simply apply
the trace from $B$ to $A$ coefficientwise.
\item[(b)] In case $R = \calR_A$ and $S = \calR_A \otimes_{\calR_A^{\inte}}
\calR'$ where $\calR'$ is a finite Galois extension of
$\calR_A^{\inte}$. 
In this case, we may simply sum over $\Aut(\calR'/\calR_A^{\inte})$ to 
obtain the trace.
\end{enumerate}
\end{remark}

\subsection{$(\sigma, \nabla)$-modules}

We now introduce the coefficient objects in Monsky-Washnitzer
cohomology.
\begin{defn}
Let $R$ be a ring of one of the types described in 
Definition~\ref{defn:moddiff} equipped with a Frobenius lift $\sigma$, 
and let $A$ be a subring preserved by $\sigma$.
(If $A$ is not specified, take it to be $K$.)
Define a \emph{$(\sigma, \nabla)$-module} over $R$,
relative to $A$,
as a finite locally free module $M$ over $R$ equipped with:
\begin{enumerate}
\item[(a)] a \emph{Frobenius structure}:
an additive, $\sigma$-linear map $F: M \to M$ (that is, $F(r \bv)
= r^\sigma F(\bv)$ for $r \in R$ and $\bv \in M$) such that the induced
$R$-linear map $\sigma^* M \to M$ is an isomorphism;
\item[(b)] an \emph{integrable connection}:
an additive, $A$-linear map $\nabla: M \to M \otimes_R
\Omega^1_{R/A}$ satisfying the Leibniz rule: $\nabla(r \bv) = r
\nabla(\bv) + \bv \otimes dr$ for $r \in R$ and $\bv \in M$, and
such that, if we write $\nabla_n$ for the induced map $M \otimes_R
\Omega^n_{R/A} \to M \otimes_R \Omega^{n+1}_{R/A}$, we have
$\nabla_{n+1} \circ \nabla_n = 0$ for all $n \geq 0$;
\end{enumerate}
subject to the compatibility condition
\[
\xymatrix{
M \ar^{F}[d] \ar^(.3){\nabla}[r] & M \otimes_R \Omega^1_{R/A}
\ar^{F \otimes d\sigma}[d] \\
M \ar^(.3){\nabla}[r] & M \otimes_R \Omega^1_{R/A}
}
\]
\end{defn}
\begin{remark}
Note that the integrability condition $\nabla_{n+1} \circ \nabla_n
= 0$ is redundant if $\Omega^1_{R/A}$ is projective of rank 1, e.g.,
if $R$ is a one-dimensional Robba ring over $A$.
In any case, it suffices to check that $\nabla_1 \circ \nabla_0 = 0$.
\end{remark}

\begin{remark}
The compatibility condition means simply that the map
$\sigma^* M \to M$ induced by $F$ is horizontal for the connection.
\end{remark}

It is worth noting that over a regular
 affinoid or dagger algebra, the hypothesis
of local freeness in the definition of a $(\sigma, \nabla)$-module
can be dropped.
\begin{lemma} \label{lem:locfree}
Let $A$ be a regular affinoid or dagger algebra, and let $M$ be
a finitely generated $A$-module equipped with a connection
$\nabla: M \to M \otimes \Omega^1_{A/K}$. Then $M$ is locally free.
\end{lemma}
\begin{proof}
The proof is basically the same as for the analogous complex
analytic statement (see, e.g., \cite[Lemme~10.3.1]{bib:pham}).
It suffices to check that for each maximal ideal $\ida$ of $A$,
$M \otimes A_{\ida}$ is free, where $A_{\ida}$ is the local ring
of $A$ at $\ida$. For $r \in A_{\ida}$, let $\ord(r)$ 
denote the order of vanishing of $r$ at $\ida$.

Choose $\bv_1, \dots, \bv_n \in M$ whose images in $M/\ida M$ form
a basis of that vector space over $A/\ida$. By Nakayama's lemma,
$\bv_1, \dots, \bv_n$ generate $M \otimes A_{\ida}$, and we must only
check that the equation $\sum c_i \bv_i = 0$, with $c_i \in A_{\ida}$,
implies that the $c_i$ are all zero. Suppose that this is not the case;
choose a counterexample that minimizes $\min_i\{\ord(c_i)\}$.
This minimum cannot be zero, since the $\bv_i$ are linearly independent
in $M/\ida M$. 

Choose vector fields $\del_1, \dots, \del_d$ on $\Maxspec A$
which at $\ida$ form a basis of the tangent space there,
and let $N_j$ be the contraction of $\nabla$ with $\del_j$. Then
for each $j$,
\[
0 = \sum_i (\del_j c_i) \bv_i + c_i (N_j \bv_i).
\]
Suppose without loss of generality that $m_1 = \min_i \{m_i\}$;
then 
we can choose $j$ so that $\del_j c_1$ is nonzero and
$\ord(\del_j c_i) < m_1$.
On the other hand, the order of vanishing of the coefficient
of $\bv_1$ in $c_i (N_j \bv_i)$ is at least $m_i \geq m_1$.
We thus obtain a new linear combination of the $\bv_i$ which vanishes,
and the minimum order of vanishing of a coefficient has decreased, 
contradicting the choice of the $c_i$.

We conclude that the equation $\sum c_i \bv_i = 0$ has only the trivial
solution, so that $M \otimes A_{\ida}$ is freely generated by the $\bv_i$.
This proves the desired result.
\end{proof}

\subsection{Cohomology of $(\sigma, \nabla)$-modules}

\begin{defn}
Let $M$ be a $(\sigma, \nabla)$-module over a ring $R$ 
(again of the type described in Definition~\ref{defn:moddiff})
relative
to a subring $A$.
Thanks to the integrability condition, we have a de~Rham complex associated to
$M$:
\[
\cdots \to
M \otimes \Omega^i_{R/A} \stackrel{\nabla_i}{\to}
M \otimes \Omega^{i+1}_{R/A} \to \cdots
\]
Define the cohomology of $M$ as the cohomology of this complex:
\[
H^i(M) = \frac{\ker(M \otimes \Omega^i_{R/A} \stackrel{\nabla_i}{\to}
M \otimes \Omega^{i+1}_{R/A})}{\im(M \otimes \Omega^{i-1}_{R/A}
\stackrel{\nabla_{i-1}}{\to}
M \otimes \Omega^{i}_{R/A})}.
\]
\end{defn}

\begin{defn}
If $f: R \to S$ is an $A$-algebra homomorphism and $S$ can be equipped
with a Frobenius lift compatible with the one given on $R$, 
we define the \emph{pullback} $f^* M$ of a $(\sigma,\nabla)$-module
$M$ over $R$ as the module $M \otimes_R S$ with the natural extra
structures. If $f: S \to R$ is an injective $A$-algebra homomorphism such that
$R$ is finite \'etale over $S$, we define
the \emph{pushforward} $f_* M$ as $M$ itself viewed as an $S$-module.
(The terminology looks backward on the algebra level, but becomes
correct when we pass to geometry and reverse arrows.)
\end{defn}

We have the following relationships between cohomology under pullback
and pushforward.
\begin{prop} \label{prop:finpull}
Let $f: R \to S$ be an injective $A$-algebra homomorphism such that
$S$ is a finite projective $R$-module, and suppose $\sigma$ is a 
Frobenius lift on $S$ mapping $R$ and $A$ into themselves.
Also assume either that $R$ and $S$ are
dagger algebras of MW-type, 
or that one of the conditions of Remark~\ref{rem:byhand}
applies.
Then for any $(\sigma, \nabla)$-module $M$ over $R$, relative to $A$,
there is an injection $H^i(M) \hookrightarrow H^i(f^* M)$ for all
$i$, whose image is that of a canonical projector on $H^i(f^* M)$.
\end{prop}
\begin{proof}
By Proposition~\ref{prop:trace} or Remark~\ref{rem:byhand},
there is a trace map $\Trace: \Omega^._{S/K} \to \Omega^._{R/K}$ such that
$\Trace \circ f$ is multiplication by a positive integer, and likewise when
$K$ is replaced by $A$.
Tensoring with $M$
and taking cohomology yields maps $H^i(M) \to H^i(f^* M) \to H^i(M)$
whose composition is multiplication by a positive 
integer, yielding the desired
result.
\end{proof}

\begin{remark} \label{rem:injfrob}
Note that the above argument only uses the $A$-linearity of $f$ in
order to conclude that the maps defined between $\Omega^._{R/K}$
and $\Omega^._{S/K}$ descend to $\Omega^._{R/A}$ and $\Omega^._{S/A}$.
This holds even if we only assume $f(A) \subseteq A$; in particular
it holds for $f = \sigma$ when $A$ is a dagger algebra of MW-type
(so that $\sigma$ is finite flat).
We thus conclude that $H^i(M) \to H^i(\sigma^* M)$ is injective; since
the Frobenius structure on $M$ consists of an isomorphism of $\sigma^* M$
with $M$, we conclude that the induced ($A$-semilinear)
action of $F$ on $H^i(M)$ is also injective.
\end{remark}

\begin{prop} \label{prop:finpush}
Let $R$ be a dagger algebra, let $f: R \to S$ be a finite
\'etale extension of $R$, and let $M$ be a $(\sigma,
\nabla)$-module on $S$. Then there is a canonical isomorphism
$H^i(M) \to H^i(f_* M)$ for all $i$.
\end{prop}
\begin{proof}
In fact, there is a canonical isomorphism of the underlying
de~Rham complexes: for each $i$, one has a canonical isomorphism
\[
S \otimes_R \Omega^i_{R/K} \cong \Omega^i_{S/K}
\]
and this remains an isomorphism after tensoring with
$M$.
\end{proof}

\section{Rigid cohomology}
\label{sec:rigid}

In this chapter, we review Berthelot's theory of rigid cohomology,
and recall its connection to the Monsky-Washnitzer theory 
(see especially Remarks~\ref{rem:compare1} and~\ref{rem:compare2}).
We remind the reader of a useful analogy: Monsky-Washnitzer
cohomology is analogous to the ``naive'' algebraic de Rham cohomology
of smooth affine schemes, which is simply the cohomology of the
complex of K\"ahler differentials of the coordinate ring.
Rigid cohomology is analogous to the ``correct'' definition of
algebraic de Rham cohomology, in which one (locally) embeds a given variety
into a smooth variety, considers the formal neighborhood, and then
computes hypercohomology of the resulting de Rham complex. The choice
of the embedding does not matter because any two such embeddings
end up being homotopy-equivalent on the level of de Rham complexes.

Our treatment here is by no means definitive. For further information,
a good introduction is \cite{bib:ber1}, but it contains very few proofs.
A more comprehensive development in the case of constant coefficients
is \cite{bib:ber2}; most of the arguments there can be extended to
nonconstant coefficients relatively easily.
The notion of isocrystals is comprehensively developed in the
unpublished \cite{bib:ber5}.
Of course, all of this theory is ultimately subsumed into the
theory of arithmetic $\mathcal{D}$-modules anyway, for which
see \cite{bib:ber7} for an introduction.

\subsection{Some rigid geometry}

We begin with some preliminaries from rigid geometry,
namely
Raynaud's construction
\cite{bib:raynaud} of the ``generic fibre'' of a formal $\calO$-scheme.
The canonical reference for this subject is of course \cite[Part~C]{bib:bgr},
to which we refer for such notions as rigid analytic spaces, admissible
open subsets, and sheaf cohomology.

\begin{defn}
Let $P$ be a formal $\calO$-scheme which is separated and of 
finite type. 
We define the \emph{generic fibre} $\tilP$ to be the set of
closed formal subschemes of $P$ which are integral, finite and flat
over $\calO$. This set comes with a specialization map
$\spe: \tilP \to P$ which associates to each subscheme its support.
\end{defn}

For $P$ affine, we can write
$P= \Spf(R)$ for some complete topologically finitely generated
$\calO$-algebra $R$, and
identify $\tilP$ with $\Maxspec (R \otimes_{\calO} K)$. 
This gives $\tilP$ the structure
of a rigid analytic space; it can be shown that this structure does not
depend on the choice of $R$. In particular, the structures glue to give
$\tilP$ the structure of a rigid analytic space for any $P$. (More precisely,
$\tilP$ is a locally $G$-ringed space for an appropriate Grothendieck
topology.)

\begin{defn}
For $P$ a formal $\calO$-scheme and $X$ a subscheme of the special
fibre of $P$, we define the \emph{tube} of $X$, denoted $]X[$, to
be the inverse image $\spe^{-1}(X)$.
\end{defn}

We make the following definition following \cite[Definition~2.3]{bib:kiehl}.
\begin{defn}
A rigid analytic space is \emph{quasi-Stein} if it is the union of
an increasing sequence $U_1 \subseteq U_2 \subseteq \cdots$ of affinoid
subdomains such that $\Gamma(U_i)$ is dense in $\Gamma(U_{i+1})$ for all
$i$.
\end{defn}
The main result about quasi-Stein spaces is ``Theorem B'' 
\cite[Satz~2.4]{bib:kiehl}.
\begin{prop}
Let $\calG$ be a coherent sheaf on a quasi-Stein space $X$. Then
the cohomology groups $H^i(X, \calG)$ vanish for $i>0$.
\end{prop}

\subsection{Overconvergent $F$-isocrystals}

In this section, we recall the notion of overconvergent $F$-isocrystals
from \cite{bib:ber1}. These play the role of the coefficient objects
in rigid cohomology; they are analogous to local systems in ordinary
de Rham cohomology or lisse sheaves in \'etale cohomology. They are
also closely related to the $(\sigma, \nabla)$-modules of 
the previous chapter; see Remark~\ref{rem:compare1}.

For simplicity, we restrict to the following situation.
(One can generalize the definitions by a glueing construction.)
\begin{hypo}  \label{hypo:rigsetup}
Let $X$ be a $k$-scheme and let $\overline{X}$ be a compactification
of $X$; put $Z = \overline{X} \setminus X$. 
Suppose that there exists a closed immersion
of $\overline{X}$ into a $\calO$-formal scheme $P$ which is smooth
in a neighborhood of $X$. Suppose further that there exists
an isomorphism $\sigma$ of $P$ with itself that induces the action of
$\sigma_K$ on $K$ and the $q$-power (absolute) Frobenius on $\overline{X}$,
and fix a choice of $\sigma$.
\end{hypo}

\begin{defn}
Suppose that $P$ is affine.
Define a \emph{strict neighborhood} of $]X[$ in $]\overline{X}[$ as
an admissible open subset $V$ containing $]X[$ such that for any affinoid
$W$ contained in $]\overline{X}[$, we can find sections $f_1, \dots, f_m
\in \Gamma(\tilP)$, whose images in
$\Gamma(\overline{X})$ cut out $Z$, and some $\lambda < 1$ such that
\[
W \cap \{x \in ]\overline{X}[: \max_i\{ |f_i(x)| \} \geq \lambda \} \subseteq W \cap V.
\]
For $P$ nonaffine, we make the same definition locally.
\end{defn}

\begin{defn}
An \emph{overconvergent $F$-isocrystal} $\calE$ on $X$, relative to $K$,
is a finite locally free $\calO_{V}$-module on some strict neighborhood
$V$ of $]X[$, equipped with a connection $\nabla: \calE \to \calE 
\otimes \Omega^1_{V}$ and an isomorphism $F: \sigma^* \calE \to \calE$
of modules with connection, subject to a convergence condition which
will not be described in detail here (it asserts that the Taylor isomorphism
is defined on a strict neighborhood of $]X[_{P \times_{\calO} P}$).
This definition appears to depend on 
the choice of $P$, but in fact there are canonical equivalences
(via pullbacks, and relying on the convergence condition) 
between the categories constructed using different $P$.
\end{defn}

\begin{remark}
We reiterate that the equivalences of categories between different
choices of $P$ are homotopy equivalences on the level
of de Rham complexes. It may be useful to imagine them as
homotopy equivalences in the true topological sense; that is, there
are maps between the two rigid analytic spaces whose compositions
are close enough to being the identity maps that they can be
``continuously deformed'' into the identity map.
For an example of this, see \cite[Theorem~5.5]{bib:mw}.
\end{remark}

\begin{remark} \label{rem:compare1}
For $X$ smooth affine, the direct limit $A$ of $\Gamma(V, \calO_V)$
over all embeddings $j: ]X[ \to V$ into a strict neighborhood
is a dagger algebra of MW-type over $K$ with special fibre $X$. 
This turns out to yield an equivalence between the
category of overconvergent $F$-isocrystals on $X$, relative to $K$,
and the category of $(\sigma, \nabla)$-modules over $A$; the main point
is that in this setting, the convergence condition is an automatic
consequence of the presence of a Frobenius structure
\cite[Th\'eor\`eme~2.5.7]{bib:ber5}.
This connection extends also to the cohomology level; see
Remark~\ref{rem:compare2}.
\end{remark}

\subsection{Rigid cohomology with and without supports}
\label{subsec:rigcoh}

We now recall Berthelot's definition of rigid cohomology, first
without supports.
\begin{defn} \label{defn:rigid}
 The \emph{rigid cohomology} of $\calE$, notated
$H^i_{\rig}(X/K, \calE)$, is given by
\[
H^i_{\rig}(X/K, \calE) = H^i(]\overline{X}[, j^{\dagger} 
(\calE \otimes \Omega^{.}_{V})),
\]
where $j^\dagger$ denotes the direct limit over all inclusions
$j: ]X[ \hookrightarrow V$ into a strict neighborhood.
This construction is independent of the choice of compactification; see
\cite[2.3]{bib:ber1}.
\end{defn}
\begin{remark}
The notation is a little peculiar, because it is $\calE$ and not $X$ that
depends on the choice of $K$. It is set this way because in the case of
constant coefficients, one omits $\calE$ but still needs to
keep track of the coefficient field.
\end{remark}

\begin{remark} \label{rem:compare2}
Recall from Remark~\ref{rem:compare1} that if $X$ is smooth affine,
then an overconvergent $F$-isocrystal $\calE$ on $X$ corresponds to a
$(\sigma, \nabla)$-module $M$ over a dagger algebra $A$ of
MW-type with special fibre $X$. This comparison extends to cohomology:
the rigid cohomology $H^i_{\rig}(X/K, \calE)$ is canonically isomorphic
to $H^i(M)$.
See \cite[Proposition~1.12]{bib:ber1} for a detailed discussion
(in the constant coefficient case, but the general case is similar).
\end{remark}

We next turn to cohomology with supports.
\begin{defn} \label{def:supports}
Let $\iota: ]Z[ \cap V \hookrightarrow V$ denote the obvious inclusion;
note that $]Z[$ is an admissible open subset of $]\overline{X}[$.
For $\calF$ a sheaf of abelian groups on $V$, put
\[
\underline{\Gamma}(\calF) = \ker(\calF \to \iota_* \iota^* \calF);
\]
this functor is left exact and admits right derived functors.
For $\calF$ coherent, by Theorem B we have $\RR^q \iota_* \iota^* \calF = 0$ for 
$q>0$, so the total derived complex $\RR^q \underline{\Gamma}(\calF)$
is isomorphic (in the derived category) to the two-term complex
\[
0 \to \calF \to \iota_* \iota^* \calF \to 0.
\]
The \emph{rigid cohomology with compact supports} of $\calE$,
notated $H^i_{c,\rig}(X/K, \calE)$, is given by
\[
H^i_{c,\rig}(X/K,\calE) = H^i(V, \RR \underline{\Gamma}
(\calE \otimes \Omega^{.}_{V})).
\]
In addition to being independent of the compactification,
this construction is also independent of the choice of $V$; see
\cite[3.3]{bib:ber1} and subsequent remarks.
\end{defn}

\begin{remark}
As one might expect, there is a natural ``forget supports'' map
\[
H^
i_{c,\rig}(X/K,\calE) \cong H^i_{\rig}(X/K,\calE),
\] and this
map is an isomorphism when $X$ is proper. We will not use
either of these facts in what follows; however, they are crucial when
considering the Weil conjectures, as in \cite{bib:me10}.
\end{remark}

\subsection{Additional properties}

We will use a few other formal properties of rigid cohomology,
particularly in Chapter~\ref{sec:final} as we prove our main results.

One has in rigid cohomology a notion of cohomology with support
along a closed subscheme, as described in \cite[2.1]{bib:ber2}.
Briefly, if $U$ is an open subscheme of $X$ and $Z = X \setminus U$,
and $\calE$ is an overconvergent $F$-isocrystal on $X$,
there is a natural surjective map of sheaves from
the sheaf $j^\dagger (\calE \otimes \Omega^.)$ to its ``restriction to $U$''
(really, the limit of the direct images of its restrictions 
to strict neighborhoods of $]U[$).
The cohomology of the kernel of this map is denoted 
$H^i_{Z,\rig}(X/K, \calE)$. By construction, it fits into a long
exact sequence with the cohomologies of $X$ and $U$; more careful
consideration shows that one actually has a long exact sequence
\begin{equation} \label{eq:excis1}
\cdots \to H^i_{T/K,\rig}(X, \calE) \to H^i_{Z/K,\rig}(X, \calE)
\to H^i_{Z \setminus T,\rig}(X \setminus T/K, \calE) \to \cdots
\end{equation}
whenever $T$ is a closed subscheme of $Z$.
(See \cite[Proposition~2.5]{bib:ber2} for the constant coefficient case;
as usual, the general case is similar.)

In order for excision to be useful, one needs to have 
a Gysin map comparing $H^i_{Z,\rig}(X/K,\calE)$
with $H^i_{\rig}(Z/K, \calE)$, at least when $X$ and $Z$ are smooth.
As it stands, one only has this with some additional restrictions,
but this will be enough for our purposes.
\begin{hypo} \label{hypo:gysin}
Suppose 
$Z \hookrightarrow X$ is a smooth pair, with $\dim Z \leq n$, admitting a
lifting $\mathcal{Z} \hookrightarrow \mathcal{X}$ over $\calO$,
such that $\mathcal{Z}, \mathcal{X}$ are smooth over $\calO$ and
there exists an \'etale
morphism $\mathcal{X} \to \Spec \calO[x_1, \dots, x_d]$ for some integer $d$,
such that $\mathcal{Z}$ maps to the vanishing locus of $d-n$ of the $x_i$.
\end{hypo}
Under Hypothesis~\ref{hypo:gysin}, Tsuzuki
\cite[Theorem~4.1.1]{bib:tsu5} constructs a
Gysin isomorphism
\[
H^i_{Z,
\rig}(X/K, \calE) \cong H^{i-2c}_{\rig}(Z/K, \left. \calE \right|_Z),
\]
where $c = \dim X - \dim Z$.

In cohomology with compact supports, one has a similar excision sequence
going the other way:
\begin{equation} \label{eq:excis2}
\cdots \to H^i_{c,\rig}(U, \calE) \to
H^i_{c,\rig}(X, \calE) \to H^i_{c,\rig}(Z, \calE) \to \cdots;
\end{equation}
see \cite[Proposition~2.5.1]{bib:tsu5}.
As one expects, all three terms in this sequence are of the same
type, so there is no need for a Gysin construction.

Finally, we will frequently compute by pushing forward
along a finite \'etale morphism. This operation behaves as expected:
in particular, one has the following relationship
(see \cite[Corollaries~2.6.5 and~2.6.6]{bib:tsu5}).
\begin{prop} \label{prop:finpush2}
Let $f: X \to Y$ be a finite \'etale morphism of separated $k$-schemes
of finite type, and let
$\calE$ be an overconvergent $F$-isocrystal on $X$.
Then for each $i$, there are canonical isomorphisms
$H^i_{\rig}(X/K, \calE) \cong H^i_{\rig}(Y/K, f_* \calE)$ and
$H^i_{c,\rig}(X/K, \calE) \cong H^i_{c,\rig}(Y/K, f_* \calE)$.
\end{prop}

\section{Relative local monodromy}
\label{sec:rel}

In this section, we concentrate on $(\sigma, \nabla)$-modules over
(one-dimensional) Robba rings and their cohomology. We
recall the notion of unipotence and the statement of
the statement of the $p$-adic local 
monodromy theorem (Crew's conjecture) for $(\sigma, \nabla)$-modules
over $\calR_K$. We then formulate an analogue for a Robba ring over a dagger
algebra, and apply this analogue to study the cohomology of
a $(\sigma, \nabla)$-module over such a Robba ring.

\subsection{Unipotence and $p$-adic local monodromy}

\begin{defn}
Let $A$ be a reduced dagger algebra
and let $M$ be a free $(\sigma, \nabla)$-module over the Robba ring $\calR_A$,
relative to $A$.
Define a \emph{unipotent basis} of $M$ to be a basis
$\bv_1, \dots, \bv_n$ for which
$\nabla \bv_i \in (\calR_A \bv_1 + \cdots + \calR_A \bv_{i-1})
\otimes \frac{dt}{t}$ for all $i$. We say $M$ is \emph{unipotent}
if it admits a unipotent basis. (Note that the definition does not
refer explicitly to the Frobenius structure.)
\end{defn}
Note that the definition of unipotence does
not depend on the implicit choice of the series parameter in the
definition of $\calR_A$;
that is, it is invariant under pullback by any continuous automorphism
of $\calR_A$.

Over a field, the definitive statement about the property of unipotence
is the $p$-adic local monodromy theorem (``Crew's conjecture'').
This theorem has been proved independently by
Andr\'e \cite{bib:andre}, Mebkhout \cite{bib:mebkhout},
and the author \cite{bib:me7}.
\begin{theorem} \label{thm:monodromy1}
Let $M$ be a free $(\sigma, \nabla)$-module over $\calR_K$. Then
there exist an integer $m \geq 0$, a
finite unramified extension $K_1$ of $K_0 = K^{\sigma^{-m}}$,
a finite Galois (\'etale) extension $\calR'$ of $\calR_{K_1}^{\inte}$,
and a continuous $K_1$-algebra isomorphism of
\[
\calR'' = \calR_{K_1} \otimes_{\calR_{K_1}^{\inte}} \calR'
\]
with 
$\calR_{K_1}$ such that $M \otimes \calR''$ is unipotent.
\end{theorem}
Note that the conclusion 
does not depend on the choice of the isomorphism
$\calR'' \cong \calR_{K_1}$, only on the fact that such an isomorphism
exists.
\begin{proof}
We explain here how to deduce the theorem from
\cite[Theorem~6.12]{bib:me7}. First recall that the category of
$(\sigma, \nabla)$-modules over $\calR_K$ is independent of the choice of
$\sigma$;
see for instance \cite[Definition~6.1.2]{bib:me13}. We may thus assume
that $\sigma$ is a power of a $p$-power Frobenius lift. (This step is needed
because \cite{bib:me7} always assumes Frobenius lifts are powers of
$p$-power Frobenius lifts. This can be avoided by instead invoking
\cite[Theorem~7.2.5]{bib:me12}.)

Now by \cite[Theorem~6.12]{bib:me7}, $M$ becomes unipotent after tensoring
up to the extension of $\calR_K$ corresponding to a ``nearly finite
separable extension'' of $k((t))$, i.e., a finite separable extension $L$ of
$k_0((t))$ with $k_0 = k^{q^{-m}}$ for some $m \geq 0$.
There is no harm in enlarging $L$, so we may replace 
$L$ by its Galois closure over $k_0((t))$.
By increasing $m$ (and correspondingly enlarging $L$), 
we can ensure that the residue field of $L$ is separable
over $k_0$, and hence coincides with the integral closure of $k_0$ in $L$;
call this field $k_1$. Let $K_1$ be the unramified extension of 
$K_0 = K^{\sigma^{-m}}$ with residue field $k_1$. Since
$\calR_{K_1}^{\inte}$ is a henselian discrete valuation ring
\cite[Lemma~3.9]{bib:me7}, we can uniquely lift $L$ to a finite \'etale
extension $\calR'$ of $\calR_{K_1}^{\inte}$. Since the residue field
of $L$ coincides with $k_1$, we can find a continuous isomorphism
$L \cong k_1((u))$ by the Cohen structure theorem; 
by lifting $u$, we obtain an isomorphism
of
\[
\calR'' = \calR_{K_1} \otimes_{\calR_{K_1}^{\inte}} \calR'
\]
with 
$\calR_{K_1}$ such that $M \otimes \calR''$ is unipotent, as desired.
\end{proof}

The purpose of this section is to salvage
Theorem~\ref{thm:monodromy1} in case $K$ is replaced by an integral
dagger algebra of MW-type. For starters, we can apply
Theorem~\ref{thm:monodromy1} with $K$ replaced by the completion
of $\Frac A$ for the affinoid topology. (One might call the
outcome of this application ``generic quasi-unipotence''.) 
However, one would like to perform a less drastic enlargement, replacing
$A$ only by a localization. Happily, this is indeed possible.
\begin{theorem} \label{thm:monodromy2}
Let $A$ be an integral affinoid or
dagger algebra of MW-type equipped with a Frobenius lift,
and let $M$ be a free $(\sigma, \nabla)$-module over $\calR_A$,
relative to $A$.
Then there exist a localization $B$ of $A$, an
integer $m \geq 0$, a
finite \'etale extension $B_1$ of $B_0 = B^{\sigma^{-m}}$,
a finite Galois (\'etale) extension $\calR'$ of $\calR_{B_1}^{\inte}$,
and a continuous $B_1$-algebra isomorphism of
\[
\calR'' = \calR_{B_1} \otimes_{\calR_{B_1}^{\inte}} \calR'
\]
with 
$\calR_{B_1}$ such that $M \otimes \calR''$ is unipotent.
\end{theorem}
The proof of Theorem~\ref{thm:monodromy2} occupies the remainder of this
chapter.

\subsection{More on unipotence}
\label{subsec:setup}

We begin our approach to Theorem~\ref{thm:monodromy2} with some preliminaries.
First, we widen the class of rings we are working with.
\begin{defn}
For $A$ a reduced affinoid algebra and $r>0$, define
$R_{A,r}$ be the ring of formal series $\sum_{i \in \ZZ} a_i t^i$
over $A$ for which there exists $c>0$ such that
\[
\lim_{i \to \pm \infty} v_A(a_i) + ri - c|i| = \infty,
\]
and put
\[
w_{A,r}\left(\sum a_i t^i \right) = \min_i \{v_A(a_i) + ri\}.
\]
For $A$ a reduced dagger algebra, define $R_{A,r}$ as the union
of $R_{B,r}$ over all fringe algebras $B$ of $A$.
Define $\Omega^1_{R_{A,r}}$ as the free $R_{A,r}$-module generated
by $dt$, with
the obvious $A$-linear derivation $d: R_{A,r} \to \Omega^1_{R_{A,r}}$.
\end{defn}

\begin{defn}
Define a \emph{module with connection} over $R_{A,r}$ as a free
$R_{A,r}$-module $M$ of finite rank, equipped with an additive
$A$-linear map $\nabla: M \to M\otimes \Omega^1_{R_{A,r}}$ satisfying
the Leibniz rule. Define the operator $D: M \to M$ so that
$\nabla \bv = D\bv \otimes \frac{dt}{t}$.
\end{defn}

The notion of unipotent basis makes sense for modules with connection;
it will help to also have a more refined notion.
\begin{defn}
Let $A$ be a reduced affinoid or dagger algebra.
A \emph{strongly unipotent basis} of a $(\sigma, \nabla)$-module
over $\calR_A$ relative to $A$, 
or a module with connection over $R_{A,r}$ for some
$r>0$,
is a basis $\bv_1, \dots, \bv_n$ such that
\[
D \bv_i \in A \bv_1 + \cdots + A \bv_{i-1} \qquad (i=1, \dots, m).
\]
\end{defn}
Note that this definition,
unlike the definition of a unipotent basis, is not stable under
continuous automorphisms of $\calR_A$ or $R_{A,r}$.

\begin{lemma}
Let $A$ be a reduced affinoid or dagger algebra.
Let $M$ be a $(\sigma, \nabla)$-module over $\calR_A$ relative to $A$, or a 
module with connection over $R_{A,r}$ for some $r>0$.
Suppose $\bv_1, \dots, \bv_n$ is a strongly unipotent basis
of $M$. Then for $\bw \in M$, $\bw$ belongs to the $A$-span of the
$\bv_i$ if and only if $D\bw$ does.
\end{lemma}
\begin{proof}
One inclusion is straightforward, so suppose that
$D\bw$ is in the span of the $\bv_i$.
If we write $\bw = \sum c_i \bv_i$, then the coefficient of
$\bv_n$ in $D\bw$ is precisely $t \frac{dc_n}{dt}$. This can only
be in $A$ if it is zero. Hence $c_n \in A$; the same argument
applied to $\bw - c_n \bv_n$ shows that $c_{n-1} \in A$, and so forth.
\end{proof}
\begin{cor} \label{cor:strongspan}
Let $A$ be a reduced affinoid or dagger algebra.
Let $M$ be a $(\sigma, \nabla)$-module over $\calR_A$ relative
to $A$, or a 
module with connection over $R_{A,r}$ for some $r>0$.
Suppose $\bv_1, \dots, \bv_n$ and $\bw_1, \dots, \bw_n$ are 
strongly unipotent bases of $M$. Then
the $A$-span of the $\bv_i$ is the same as the $A$-span of the $\bw_i$.
\end{cor}

\begin{prop} \label{prop:unipbasis}
Let $A$ be a reduced affinoid or dagger algebra.
Let $M$ be a free unipotent $(\sigma, \nabla)$-module over $\calR_A$
relative to $A$, or a module with connection over $R_{A,r}$ for some $r>0$.
Then $M$ admits a strongly unipotent basis.
\end{prop}
\begin{proof}
Let $\bw_1, \dots, \bw_n$ be any unipotent basis of $M$.
We construct the $\bv_i$ so that $\bv_1, \dots, \bv_i$ has the same
span as $\bw_1, \dots, \bw_i$ for $i=1, \dots, n$; this will
ensure that $\bv_1, \dots,\bv_n$ form a basis of $M$.

To begin with, set $\bv_1 = \bw_1$. Given $\bv_1, \dots, \bv_i$, write
\[
D \bw_{i+1} = \sum_{j=1}^i a_{ij} \bv_j.
\]
Suppose that $a_{ij} \in A$ for $j=l+1, \dots, i$.
Write $a_{il} = b_{il} + \sum_{m \neq 0} c_{ilm} t^m$ with 
$b_{il}, c_{ilm} \in A$.
By Proposition~\ref{prop:deriv} (or its easy analogue for $R_{A,r}$), 
we can find $e_{il} \in \calR_A$
such that $t \frac{de_{il}}{dt} = a_{il} - b_{il}$. We then observe that
\[
D(\bw_{i+1} - e_{il} \bv_j) \in \calR_A \bv_1 + \cdots +
\calR_A \bv_{l-1} + A\bv_l + \cdots + A\bv_i.
\]
Thus by repeating this procedure with $l = i,i-1,\dots,1$ in succession,
we can modify $\bw_{i+1}$, without changing the span of
$\bv_1,\dots,\bv_i,\bw_{i+1}$, until we ultimately obtain a new vector
which we call $\bv_{i+1}$, with the property that
\[
D \bv_{i+1} \in A \bv_1 + \cdots + A \bv_i.
\]
Hence this iteration produces the desired result.
\end{proof}

We will also need one lemma from elementary number theory.
\begin{lemma} \label{lem:bounddenom}
Let $m,l,e$ be integers with $l,e$ positive. Then the product
\[
\prod_{i=1}^l \frac{m+x+i}{i}
\]
in $\QQ_p[x]/(x^e)$ becomes an element of $\ZZ_p[x]/(x^e)$ after
multiplication by $p^a$ for some positive integer $a \leq 
(e-1) \lceil \log_p (|m| + l)\rceil$.
\end{lemma}
\begin{proof}
The coefficient of $x^j$ in the product can be written as the
integer $\prod_{i=1}^l (m+i)/i$ times a sum of products of
$j$ fractions, each of the form $c/(c+m)$ for some $c \in \{1, \dots, l\}$.
From this description, the claim is evident.
\end{proof}

\subsection{Generic unipotence}

In this section, we show that in some sense,
unipotence of a module with connection
over the completed fraction field implies unipotence
after a localization.

\begin{lemma} \label{lem:iterate}
Let $M$ be a unipotent module with connection over $R_{K,r}$
for some $r>0$.
Let $e$ be the nilpotency index of the matrix via which
$D$ acts on a strongly unipotent basis of $M$.
Define the functions $f_l: M \to M$ for
$l=0,1,\dots$ by setting $f_0(\bw) = D^{e-1} \bw$ and
\[
f_l(\bw) = (1-D^2/l^2)^e f_{l-1}(\bw).
\]
Then for any $\bw \in M$, the sequence $\{f_l(\bw)\}_{l=0}^\infty$
converges (under $w_{K,r}$)
in $M$ to an element $f(\bw)$ such that $\nabla f(\bw) = 0$.
\end{lemma}
\begin{proof}
Let $\bv_1, \dots, \bv_n$ be a strongly nilpotent basis of $M$,
and define the upper triangular nilpotent matrix $X$ over $K$ by
\[
D\bv_j = \sum_{i} X_{ij} \bv_i;
\]
by the choice of $e$, $X^e = 0$ but $X^{e-1} \neq 0$.
Write $\bw = \sum_i a_i \bv_i$ and $a_i = \sum_{m \in \ZZ} a_{im} t^m$,
and put $h = \max\{0, \max_{i,j} \{-v_K(X_{ij}) \}\}$.

Write $f_l(\bw) = \sum_i b_{il} \bv_i$ and write 
$b_{il} = \sum_{m \in \ZZ}
b_{ilm} t^m$. We then have
\begin{equation} \label{eq:iterate}
b_{jlm} = \sum_i \left( (mI+X)^{e-1} \frac{(I-(mI+X)^2)^e \cdots (l^2 I - 
(mI+X)^2)^e}{(l!)^{2e}}
\right)_{ij} a_{im}.
\end{equation}
For $m=0$, this equation yields $b_{jlm} = \sum_i (X^{e-1})_{ij} a_{im}$,
which does not depend on $l$. For $0 < |m| \leq l$, we have
$b_{jlm} = 0$ because $(m^2 I - (|m|I+X)^2)^e = 0$.
For $|m| > l$, we can interpret the parenthesized expression in
\eqref{eq:iterate} as the product in $\QQ[x]/(X^e)$ of $(m+X)^{e-1}$ times
$2e$ expressions, each of the form $\prod_{i=1}^l (\pm(m+X)+i)/i$.
By Lemma~\ref{lem:bounddenom}, we thus have
\[
w_{K,r}(b_{i(l+1)} - b_{il})
\geq \min_{1 \leq i \leq n, |m| > l} \{ v_K(a_{im}) + mr
- 2e(e-1) \lceil \log_p (|m| + l+1) \rceil - (e-1)h\} 
\]
for $0 < r \leq s$.

Choose $c,d>0$ such that $v_K(a_{im}) + mr \geq c|m| - d$ for all $m$.
We then have
\begin{align*}
w_{K,r}(b_{i(l+1)} - b_{il}) &\geq \min_{x \geq l}
\{ cx - 2e(e-1) \log_p (x + l+1)\} - 2e(e-1) - d - (e-1)h\\
&\geq \min_{x \geq l} \{cx - 2e(e-1) \log_p (2x+1)\} - 2e(e-1) - d - (e-1)h.
\end{align*}
The expression on the right tends to infinity with $l$;
moreover, for any $0<c'<c$,
we can choose $d'>0$ so that the right side is at least
$c'l - d'$ for all $l$.

We conclude that $f_l(\bw)$ converges under $w_{K,r}$
to an element of $M$ in the $K$-span of the $\bv_i$. Moreover,
the limit is in the image of $D^{e-1}$ on this span, so it must
be killed by one more application of $D$. This yields the desired
result.
\end{proof}

Next, we use the iteration from the previous lemma to produce 
horizontal elements of a $(\sigma, \nabla)$-module
over a dagger algebra.
\begin{lemma} \label{lem:iterdagger}
Let $A$ be a reduced integral dagger algebra, and let $L$ be the
completion of $\Frac A$ for the affinoid topology.
Let $M$ be a module with connection over $R_{A,r}$
for some $r>0$
such that $M \otimes R_{L,r}$ is unipotent. Then
for any $\bw \in M$, the element $f(\bw)$ of $M \otimes R_{L,r}$
constructed by Lemma~\ref{lem:iterate} lies in $M$.
\end{lemma}
\begin{proof}
Apply Lemma~\ref{lem:iterate} to $M \otimes R_{L,r}$ and retain
all notation therein.
Choose a basis $\be_1, \dots, \be_n$ of $M$, and define the matrix $N$
by
\[
D\be_j = \sum_i N_{ij} \be_j.
\]
Choose a fringe algebra $B$ over $A$ such that everything in sight
defined over $R_{A,r}$ is actually defined over $R_{B,r}$ 
(i.e., the matrix $N$, and the coefficients of the $\be_i$
in the expansion of $\bw$).

Write $f_l(\bw) = \sum_i g_{il} \be_i$. For some $c',d'> 0$,
we then have the estimates
\begin{align*}
  w_{L,r}(g_{i(l+1)} - g_{il}) &\geq c'l - d' \\
w_{B,r}(g_{i(l+1)} - g_{il}) &\geq -2el \left( \frac{1}{p-1}
+ \max\{0, -w_{B,r}(N)\} \right) - d'.
\end{align*}
Namely, the first estimate is obtained from the corresponding
estimate for $b_{i(l+1)}-b_{il}$ in the proof of Lemma~\ref{lem:iterate};
the change of basis from the $\bv_i$ to the $\be_i$ only shifts
the estimate by a constant.
The second estimate is a straightforward
consequence of the definition of $f_l$.

By Proposition~\ref{prop:hadamard}, for any rational $\epsilon \in (0,1]$,
we can find a fringe algebra $C$ of $A$ containing $B$ such that
\[
|x|_{\sup,C} \leq |x|_{\sup,A}^{1-\epsilon} |x|_{\sup,B}^{\epsilon}
\]
for all $x \in B$. For $y \in R_{B,r}$, this implies
\[
w_{C,r}(y) \geq (1-\epsilon) w_{L,r}(y) + \epsilon w_{B,r}(y).
\]
In particular, by taking $\epsilon$ sufficiently small, we can ensure that
\[
w_{C,r}(g_{i(l+1)} - g_{il}) \geq c''l -d''
\]
for some $c'',d'' > 0$.
We thus conclude that the sequences $\{g_{il}\}_{l=0}^\infty$
converge under $w_{C,r}$, so their limits in $R_{L,r}$ are actually
in $R_{A,r}$. (It must be checked that
each limit $\sum a_m t^m$ satisfies the condition
$v_D(a_m) + mr - c|m| \to \infty$ for some $c>0$ and some fringe algebra $D$,
but this follows by another application of Proposition~\ref{prop:hadamard}
as above.)
Hence $f(\bw) \in M$, as desired.
\end{proof}

Finally, we use what we have to obtain a result in the vein
of Theorem~\ref{thm:monodromy2}.
\begin{prop} \label{prop:mono1}
Let $A$ be a reduced integral dagger algebra, and let $L$ be the
completion of $\Frac A$ for the affinoid topology.
Let $M$ be a module with connection over $R_{A,r}$ for some $r>0$,
such that $M \otimes R_{L,r}$ is unipotent.
Then $M \otimes R_{B,r}$ is unipotent for some localization $B$ of $A$.
\end{prop}
\begin{proof}
We proceed by induction on the rank of $M$.
Let $M_0$ be the span of a strongly unipotent basis of $M \otimes R_{L,r}$,
and let $e$ be the nilpotency index of the action of $D$ on $M_0$.
Then as $\bw$ varies over $M \otimes R_{L,r}$, $f(\bw)$ varies
over the image of $D^{e-1}$ on $M_0$.
Since $M \otimes_A L$ is dense in $M \otimes R_{L,r}$,
even if $\bw$ varies only over $M \otimes_A L$, $f(\bw)$
still covers the entire image of $D^{e-1}$ on $M_0$.  In particular,
$f(\bw)$ is nonzero for some $\bw \in M$.

Choose $\bv \in M$ in the image of $f$ with $\bv \neq 0$; 
then $\nabla \bv = 0$.
Since the free submodule of $M \otimes R_{L,r}$ generated
by $\bv$ is a direct
summand (because $\bv$ belongs to the $L$-span of a strongly unipotent
basis by Corollary~\ref{cor:strongspan}),
for some localization $B$ of $A$,
the free submodule of $M \otimes R_{B,r}$ spanned by $\bv$ is a
direct summand.
(What we are using here is that 
if a row vector over $R_{A,r}$ is the first row
of an invertible square matrix over $R_{L,r}$, then the same is true
over $R_{B,r}$ for some localization $B$ of $A$. But that is clear:
any suitably close approximation of the original square matrix will
still be invertible.) We may thus quotient by the submodule
generated by $\bv$ and
apply the induction hypothesis to deduce the desired result.
\end{proof}

\subsection{Generic quasi-unipotence}

We now proceed to the proof of Theorem~\ref{thm:monodromy2}.
First, let us reconcile the auxiliary language in which 
Proposition~\ref{prop:mono1} is couched with the situation of
Theorem~\ref{thm:monodromy2}.
\begin{prop} \label{prop:relunip}
Let $A$ be a reduced integral dagger algebra, and let $L$ be the
completion of $\Frac A$ for the affinoid topology.
Let $M$ be a free $(\sigma, \nabla)$-module over $\calR_A$,
relative to $A$, such that $M \otimes \calR_L$ is unipotent.
Then $M \otimes \calR_B$ is unipotent for some localization $B$ of $A$.
\end{prop}
\begin{proof}
Choose a basis $\be_1, \dots, \be_n$ of $M$.
Then for $r>0$ sufficiently small, the $(\calR_A \cap R_{A,r})$-span 
$N$ of the
$\be_i$ will be a module with connection which becomes unipotent
after tensoring with $R_{L,r}$. We thus conclude that
$N \otimes R_{B,r}$ admits a strongly unipotent basis 
$\bv_1, \dots, \bv_n$ for some localization $B$ of $A$.
On the other hand, $N \otimes (\calR_L \cap R_{L,r})$ also admits a 
strongly unipotent basis; by Corollary~\ref{cor:strongspan},
these two bases have the same $L$-span within $N \otimes R_{L,r}$. 
Consequently
the $\bv_i$ form a basis of
\begin{align*}
(N \otimes \calR_L) \cap (N \otimes R_{B,r}) &= 
N \otimes (\calR_L \cap R_{B,r}) \\
& = N \otimes (\calR_B \cap R_{B,r}).
\end{align*}
This yields the desired result.
\end{proof}

\begin{remark}
In some cases, it may be possible to show that the conclusion of
Proposition~\ref{prop:relunip} holds without the localization.
For instance, this can be shown in case $M$ is a
``log-$(\sigma, \nabla)$-module'' over $\calR^+$ (in which
the module of differentials is enlarged to include $dt/t$);
this is a variant of ``Dwork's trick'' \cite[Lemma~6.3]{bib:dej4}.
\end{remark}

Theorem~\ref{thm:monodromy2} can now be deduced as follows.
\begin{proof}[Proof of Theorem~\ref{thm:monodromy2}]
We will treat only the case where $A$ is
a dagger algebra; for the affinoid case, the argument is the same
except that any fringe algebra of $A$ should be replaced by $A$ itself.

Let $L$ be the completion of $\Frac A$ for the affinoid topology.
Apply Theorem~\ref{thm:monodromy1} to obtain an integer $m \geq 0$,
a finite unramified extension $L_1$ of $L_0 = L^{\sigma^{-m}}$,
a finite \'etale extension $\calR'_L$ of $\calR^{\inte}_{L_1}$,
and a continuous $L_1$-algebra isomorphism of
\[
\calR''_L = \calR_{L_1} \otimes_{\calR^{\inte}_{L_1}} \calR'_L
\]
with $\calR_{L_1}$ such that $M \otimes \calR''_L$ is unipotent.
Recall that $M \otimes \calR''_L$ is actually unipotent for any
choice of the last isomorphism.

Let $X$ be the special fibre of $A_0 = A^{\sigma^{-m}}$, and let $Y$ be 
the normalization of $X$ in the residue field of $L_1$. Since $L_1/L_0$
is unramified, its residue field extension is separable: that is, the
natural finite morphism $f: Y \to X$ is generically \'etale. The residue field
of $\calR^{\inte}_{L_1}$ can be canonically identified with $K(Y)((t))$,
where $K(Y)$ is the function field of $Y$ (i.e., the residue field of $L_1$);
the residue field of $\calR'_L$ is a finite separable extension of
$K(Y)((t))$. 

We can choose an open dense subset $X'$ of $X$ such that the following
conditions hold.
\begin{itemize}
\item[(a)]
The restriction of $f$ to $Y' = f^{-1}(X')$ is finite \'etale.
\item[(b)]
For $R$ the coordinate ring of $Y'$,
the integral closure of $R[[t]]$ in the residue field of $\calR'_L$
is isomorphic over $R$ to (another copy of) $R[[t]]$.
\item[(c)]
The complement $X \setminus X'$ is the zero locus of some regular
function on $X$.
\end{itemize}
Choose a localization $A'_0$ of $A_0$ with special fibre $X'$,
and let $A'_1$ be the integral closure of $A'_0$ in $L_1$. Then
$A'_1$ is a dagger algebra of MW-type with special fibre $Y'$,
and the isomorphism in (b) lifts (noncanonically) to an isomorphism
of the integral closure of $\calR_{A'_1}$ in $\calR''_L$ with
$\calR_{A'_1}$.

We may thus apply Proposition~\ref{prop:relunip} to deduce that
$M$ becomes unipotent after tensoring with the integral closure of
$\calR_{A''_1}$ in $\calR''_L$, for some localization $A''_1$ of $A'_1$.
Choose a localization $B_0$ of $A'_0$ whoes integral closure $B_1$
in $L_1$ contains $A''_1$, and put $B = B_0^{\sigma^m}$. 
Let $\calR''$ denote the integral closure
of $\calR_{B_1}$ in $\calR''$; we now have a continuous $B_1$-linear
isomorphism of $\calR''$ with $\calR_{B_1}$, under which
$M \otimes \calR''$ becomes unipotent. This is the desired result.
\end{proof}

\section{Cohomology of curves}
\label{sec:cohcurves}

In this chapter, we establish the finite dimensionality of the
cohomology of an overconvergent $F$-isocrystal on an affine curve.
Strictly speaking, this chapter can
be skipped, as there are several other techniques available for 
establishing this fact (see Remark~\ref{rem:cohcurve}).
However, we decided to include this argument because it is most in 
the spirit of the rest of the paper (in particular, it can be used
in the relative setting, although we will not do so), and
because its strategy may be
of use for computational purposes.

\subsection{Near-isomorphisms}

Given a short exact sequence
\[
0 \to A \to B \to C \to 0
\]
of $K$-vector spaces, $B$ is finite dimensional if and only if $A$ and
$C$ both are. Thus within the abelian category of
$K$-vector spaces, the full subcategory of finite dimensional vector spaces is
a Serre subcategory, and we can form the quotient category.

\begin{defn}
We will denote properties of morphisms that hold in the
quotient category by the prefix ``near''. That is, a $K$-linear transformation
$f: A \to B$ is a \emph{near-injection} if $\ker(f)$
is finite dimensional, a \emph{near-surjection} if $\coker(f)$ is
finite dimensional, and a \emph{near-isomorphism} if it is both
near-injective and near-surjective.
\end{defn}

From the Serre construction, the following facts are evident:
\begin{itemize}
\item
If $f: A \to B$ and $g: B \to C$ are $K$-linear transformations such that
$g \circ f$ is a near-injection (resp.\ is a near-surjection),
then so is $f$ (resp.\ so is $g$).
\item
If  $f: A \to B$ and $g: B \to C$ are $K$-linear transformations
and any two of $f,g, g\circ f$ are near-isomorphisms, then so is the third.
\item
Let
\[
\xymatrix{
A \ar[r] \ar[d] & B \ar[d] \\
C \ar[r] & D
}
\]
be a commuting square of $K$-linear transformations. If any three of
the maps are near-isomorphisms, then so is the fourth.
\end{itemize}


\subsection{Constant coefficients}

\begin{hypo} \label{hypo:curve}
Throughout this section unless otherwise specified,
let $A$ be a dagger algebra of MW-type
whose special fibre $X = \Spec R$ is a (smooth affine) irreducible
curve over $k$ such that, for $\overline{X}$ a smooth compactification
of $X$, the points of the reduced scheme 
$\overline{X} \setminus X$ are geometrically reduced.
(This is automatic if $k$ is perfect.)
\end{hypo}

In this section, we show that the cohomology of $A$,
 as a $(\sigma,\nabla)$-module over itself, is finite dimensional.
This result was first shown by Monsky \cite{bib:monsky}; our argument
is a form of that argument recast in a setting due to Crew \cite{bib:crew2},
which will pave the way for an argument involving coefficients.
(We will return to this setting in Chapter~\ref{sec:push}.)

First recall a consequence of \cite[Theorem~6.2]{bib:mw}.
\begin{lemma} \label{lem:frac}
Any finitely generated $A$-module which is annihilated by some element of $A$
is finite dimensional over $K$. In particular,
any homomorphism of finitely generated
$A$-modules which becomes an isomorphism over
$\Frac(A)$ is a near-isomorphism of $K$-vector spaces.
\end{lemma}
This lemma can also be proved by reducing to the case $A = K \langle t
\rangle^\dagger$ by Noether normalization, then using Weierstrass
preparation to show that any ideal of $K \langle t \rangle^\dagger$
is generated by a polynomial in $t$.

Following Crew \cite[7.2--7.3]{bib:crew2}, we associate to each point of
$\overline{X} \setminus X$ a ``local algebra'' isomorphic to a Robba ring.
First (again following Crew) we do this for points on $X$.
Given $x \in X$ with residue field $k'$, let $\idm_x$ be the preimage in
$A^{\inte}$ of the maximal ideal of $R$ at $x$. Then by the
Cohen structure theorem, the $\idm_x$-adic
completion of $A^{\inte}$ is homeomorphic to $\calO' \llbracket t \rrbracket
= \calR^{+,\inte}_{K'}$, for $K'$ the unramified
extension of $k$ with residue field $k'$ and $\calO'$ the ring of integers
of $K'$.

Now we consider the missing points.
\begin{defn}
Let $x_1, \dots, x_m$ be the (closed) points of $\overline{X} \setminus X$,
and let $K_i$ be the unramified extension of $K$ whose residue field $k_i$
is the residue field of $x_i$.
Choose an open affine $U$ of $\overline{X}$
containing all of the $x_i$,
and let $B$ be a dagger algebra with special fibre $U$.
Let $\calR_i$ be a copy of the Robba ring with residue field
$K_i$ and series parameter $t_i$, 
and choose an embedding $B^{\inte} \hookrightarrow \calR_i^{+,\inte}$
as in the previous paragraph.
This extends to an embedding $C^{\inte} \hookrightarrow \calR_i^{\inte}$,
for $C$ the localization of $B$ with special fibre $U - \{x_1,\dots,x_m\}$.
Choose an embedding of $A$ into $C$, and let $\rho_i: A^{\inte}
\to \calR_i^{\inte}$ be the resulting embedding; note that
any Frobenius lift on $A$ induces a compatible Frobenius
lift on $C$, hence also on $\calR_i$. Let $\rho:
A^{\inte} \to \oplus_i \calR_i^{\inte}$ denote the direct sum of the
$\rho_i$.
\end{defn}

\begin{remark}
We mention in passing Crew's description of the above construction, which
is more geometric.
Choose a formally smooth lifting of $\overline{X}$, 
and let $\calR_i$ be the union, over
strict neighborhoods $V$ of the tube of $X$, of the ring of rigid
analytic functions on the intersection
of $V$ with the tube of $x_i$.
\end{remark}

The embedding $\rho$ is a ``faithful'' description of the dagger algebra
$A$ in the following sense.
(For a much stronger result, see \cite[Theorem~7.5]{bib:crew2}.)
\begin{prop} \label{prop:closed}
The image of $\rho$ is closed for the
$\pi$-adic topology.
\end{prop}
The proposition is a consequence of the following lemma
together with the weak completeness of $A$.
\begin{lemma}
Let $R = k[y_1,\dots, y_n]/\ida$ be the coordinate ring of a smooth affine curve
$X$ over $k$, and let $\overline{X}$ be a smooth compactification of $X$.
For $r \in R$, let $d(r)$ denote the maximum pole degree of $r$ at
any point of $\overline{X} \setminus X$. Then for some $c>0$, every $r
\in R$ can be lifted to an element of $k[y_1,\dots,y_n]$ of degree at most
$cd(r)$.
\end{lemma}
\begin{proof}
It is enough to produce a single set of generators of $R$ for which the
result holds. Let $x_1, \dots, x_m$ be the points of $\overline{X} \setminus
X$, and for $r \in R$, let $p(r)$ be the $m$-tuple whose $i$-th component
is $\max\{0, -\deg_{x_i}(r)\}$. By Riemann-Roch, $p(r)$ is the zero tuple
if and only if $r$ is in the integral closure of $k$.

The submonoid of $\ZZ_{\geq 0}^m$ generated by the image of $p$ is 
finitely generated (again by Riemann-Roch, which implies that the submonoid
contains all but finitely elements
of $\ZZ_{\geq 0}^m$);
choose $z_1, \dots, z_l \in R$ whose images under $p$ generate this monoid.
Then there is a straightforward division algorithm to write
$r \in R$ in terms of the $z_i$: 
manufacture a product of $z_j$'s with the same leading term as $r$ 
in the pole expansion
at each $x_i$, subtract off this monomial, and repeat. 
This gives the desired estimate.
\end{proof}

\begin{remark}
Note that the images of the individual $\rho_i$ need not be closed if $m>1$;
for instance,
consider the example $A = K \langle x,x^{-1} \rangle^\dagger$ with
embeddings into $\calR_K$ given by $x \mapsto t$ and $x \mapsto t^{-1}$.
On the other hand, if we $p$-adically complete $A^{\inte}$ and $\calR_i^{\inte}$,
then the image of the map induced by $\rho_i$ is closed.
\end{remark}

By Riemann-Roch, the maps
\[
R \to \oplus_i (k_i ((\overline{t}_i)) / k_i \llbracket \overline{t}_i \rrbracket),
\qquad
\Omega^1_{R/k} 
\to \oplus_i (\Omega^1_{k_i ((\overline{t}_i))/k_i} / \Omega^1_{k_i 
\llbracket \overline{t}_i \rrbracket/k_i})
\]
induced by the reduction of $\rho$ are near-isomorphisms of $k$-vector spaces.
Our next lemmas extend this property to $\rho$ itself.
\begin{lemma} \label{lem:near2}
For any $a,b \in A$ with $b \neq 0$,  the map
\[
\rho: \frac{a}{b} A \to \oplus_i \calR_i/\calR_i^+
\]
is a near-isomorphism of $K$-vector spaces. In particular, we can choose
$a$ so that the map is injective when $b=1$, and we can choose $b$ so that
the map is surjective when $a=1$.
\end{lemma}
\begin{proof}
The near-injectivity and near-surjectivity for any given $a,b$ are
equivalent to the same statements for $a,b=1$, because the inclusions
$aA \hookrightarrow (a/b)A$ and $aA \hookrightarrow A$ are near-isomorphisms
by Lemma~\ref{lem:frac}.
It thus suffices to prove the claims of the final sentence.

The injectivity is easy: for any $a \in A^{\inte}$ whose reduction has
a pole at each $x_i$, the reduction of the map $\rho: aA^{\inte} \to
\oplus_i \calR_i^{\inte}/\calR_i^{+,\inte}$ is injective by Riemann-Roch.
Thus $\rho$ injects $aA$ into 
$\oplus_i \calR_i/\calR_i^+$.

As for surjectivity, again by Riemann-Roch we can choose $b \in A^{\inte}$
with nonzero reduction so that the reduction of
$\rho: \frac{1}{b} A^{\inte} \to
\oplus_i \calR_i^{\inte}/\calR_i^{+,\inte}$ is surjective.
Thus the image of $\frac{1}{b} A$ in $\oplus_i \calR_i/\calR_i^+$
is dense; since the image is also closed by Proposition~\ref{prop:closed},
it is full.
\end{proof}

\begin{lemma} \label{lem:near3}
The map
\[
d\rho: \Omega^1_{A/K} \to \oplus_i \Omega^1_{\calR_i/K}/\Omega^1_{\calR_i^+/K}
\]
is a near-isomorphism of $K$-vector spaces.
\end{lemma}
\begin{proof}
Again, the near-injectivity follows from the corresponding statement
on reductions.
As for near-surjectivity, 
let $u$ be any element of $A$ which is not integral over $K$.
By Lemma~\ref{lem:near2}, we can choose $w \in A$ such that $\frac{1}{w} A \to 
\oplus_i \calR_i/\calR_i^+$ is surjective; we can thus find $v \in A$
such that $v/w$ maps to the class of $dt_i/du$ in $\calR_i/\calR_i^+$
for all $i$. This means the map
\[
\frac{v}{w} \Omega^1_{A/K} \to \oplus_i \Omega^1_{\calR_i/K}
/ \Omega^1_{\calR_i^+/K}
\]
induced by $d\rho$ is surjective. As in Lemma~\ref{lem:near2},
this means that $d\rho$ is near-surjective, as desired.
\end{proof}

\begin{prop} \label{prop:const}
Let $A$ be a dagger algebra of MW-type whose special fibre $X$ is 
one-dimensional over $k$. Then $d: A \to \Omega^1_{A/K}$ is a near-isomorphism,
that is, $H^0_{\rig}(X/K)$ and $H^1_{\rig}(X/K)$ are finite dimensional.
\end{prop}
In particular, in the statement of the proposition, we do not assume that
$X$ is irreducible or
that $\overline{X} \setminus X$ is geometrically reduced. However,
we may make these assumptions without loss of generality, the former
because formation of $\Omega^1$ commutes with direct sums, the latter
because we can replace $k$ by a finite extension thanks to
Proposition~\ref{prop:finpull}.
Thus we may continue to assume Hypothesis~\ref{hypo:curve}.
\begin{proof}
By Lemma~\ref{lem:near2} and Lemma~\ref{lem:near3},
the horizontal arrows in the diagram
\[
\xymatrix{
A \ar[r] \ar^\nabla[d] & \oplus_i \calR_i/\calR_i^+ \ar^\nabla[d] \\
\Omega^1_{A/K} \ar[r] & \oplus_i \Omega^1_{\calR_i/K}/\Omega^1_{\calR_i^+/K}
}
\]
are near-isomorphisms. The right vertical arrow is also a
near-isomorphism: in each factor, the kernel is trivial and the cokernel is
generated by $dt_i/t_i$. Since three of the arrows are near-isomorphisms,
the fourth must be as well.
\end{proof}
In fact, with a bit more care, one can extract from this argument that
$H^1_{\rig}(X/K)$ has the expected dimension, namely twice the genus of a
smooth compactification $\overline{X}$, plus the degree of 
$\overline{X} \setminus X$, minus one.

\subsection{Interlude: matrix factorizations}

To handle general $(\sigma, \nabla)$-modules, we will need a factorization
lemma in the vein of those from \cite[Chapter~6]{bib:me7}, but which
accounts for the embedding $\rho$. We develop that lemma in this section.

First recall the statement of \cite[Proposition~6.5]{bib:me7}.
\begin{lemma} \label{lem:factor1}
Let $U$ be an invertible $n \times n$ matrix
over $\calR_K$. Then $U$ can be written as a product
$VW$, with $V$ an invertible matrix over $\calR_K^{\inte}[\frac 1p]$
and $W$ an invertible matrix over $\calR_K^+$.
\end{lemma}

With a bit more care, we can shift the denominators from $V$ to $W$.
First, we recall some terminology from linear algebra.
\begin{defn}
Recall that an \emph{elementary matrix} over a ring $R$ is defined as
a matrix over $R$ obtained from the identity matrix by performing one of
the following operations:
\begin{enumerate}
\item[(a)] multiplying a row by a unit in $R$;
\item[(b)] switching two rows;
\item[(c)] adding a multiple of one row to another.
\end{enumerate}
Multiplying a matrix $U$ by an elementary matrix on the right effects
the corresponding operation on $U$.
\end{defn}

\begin{lemma} \label{lem:factor2}
Let $U$ be an invertible $n \times n$ matrix over 
$\calR_K$. Then $U$ can be written as a product
$VW$, with $V$ an invertible matrix over $\calR_K^{\inte}$
and $W$ an invertible matrix over $\calR_K^+$.
\end{lemma}
\begin{proof}
By Lemma~\ref{lem:factor1}, we may reduce to the case where
$U$ itself is invertible over $\calR_K^{\inte}[\frac 1p]$.
Namely, a general $U$ admits the decomposition $V_1 W_1$ from 
Lemma~\ref{lem:factor1}.
If we can then decompose $V_1 = V_2 W_2$
with $V_2$ invertible over $\calR_K^{\inte}$
and $W_2$ invertible over $\calR_K^+$, we may take $V = V_2$
and $W = W_2W_1$ and we are done.

We thus assume that $U$ is invertible over $\calR_K^{\inte}[\frac 1p]$.
We may also assume
that $U$ has integral entries, though its inverse may not.

We induct on $v_K(\det(U))$; if it is zero, then $U^{-1}$ has
integral entries and we are done.
Otherwise, the reduction of $U$ is a matrix over $k((t))$
 with less than full rank.
By the usual theory of elementary divisors over $k \llbracket t \rrbracket$, 
we can find a
sequence of elementary row operations on the reduction which will
produce a zero row. Namely, perform the Euclidean algorithm on entries
in the first column until at most one is nonzero; if there
remains a nonzero entry, strike out that row. Then perform the Euclidean
algorithm on entries in the second column in the remaining rows, and so on.

The coefficients of the corresponding elementary matrices can be lifted
to $\calO \llbracket t \rrbracket$. The product of these
lifts, together with an elementary matrix that divides the offending row by
$\pi$, yields an
invertible matrix over $\calR_{K}^+$ whose product on the right with
$U$ has determinant of strictly smaller valuation than that of $U$.
Thus the induction hypothesis can be used to complete the decomposition,
yielding the desired result.
\end{proof}

We now introduce Hypothesis~\ref{hypo:curve} and give a lemma
of a similar flavor.
\begin{lemma} \label{lem:multnear}
For any sequence $U_1, \dots, U_m$ with $U_i$ an invertible $n \times n$
matrix over $\calR_i$, there exists a matrix $V$ over $\Frac A$
such that $U_i \rho_i(V)$ is an invertible matrix over $\calR_i^+$
for $i=1, \dots, n$.
\end{lemma}
\begin{proof}
By Lemma~\ref{lem:factor2}, we can find an invertible matrix $Q_i$
over $\calR_i^+$ such that $Q_iU_i$ is invertible over $\calR_i^{\inte}$.
By Riemann-Roch, we can simultaneously
approximate the reductions of the $Q_iU_i$ arbitrarily
well by the reductions of matrices over $(\Frac A)^{\inte}$.
In particular, we can choose an invertible matrix $P$ over $(\Frac A)^{\inte}$
such that for some $r>0$,
$w_l(Q_iU_i \rho_i(P)) > 0$ for $0 < l \leq r$ and $i=1,\dots,m$. 
We can then find an invertible matrix $R_i$ over $\calR_i^{+,\inte}$
such that $v_K(R_i Q_i U_i \rho_i(P) -I) > 0$.
For simplicity, put $U'_i = R_i Q_i U_i \rho_i(P)$ and put
$h = \min_i \{v_{K_i}(U'_i - I)\} > 0$.

By Lemma~\ref{lem:near2}, we can find $b \in A^{\inte}$ such that
$\rho: b^{-1}A^{\inte} \to \oplus_i \calR_i^{\inte}/\calR_i^{+,\inte}$ 
is surjective.
We construct sequences $\{V_j\}_{j=1}^\infty, \{W_{i,j}\}_{j=1}^\infty$
of matrices, 
where $V_j$ has entries in $b^{-1}A^{\inte}$ and $W_{i,j}$ has entries
in $\calR_i^{+,\inte}$, such that
\[
v_K(U'_i \rho_i(V_j) - W_{i,j}) \geq j h.
\]
Namely, we start with $V_1 = W_{i,1} = I$. Given $V_j$ and $W_{i,j}$
for some $j$, write $U'_i \rho_i(V_j) - W_{i,j} = X_{i,j}$.
Then choose matrices $Y_j$ over $\pi^j b^{-1}A^{\inte}$ and $Z_{i,j}$ over
$\calR_i^{+,\inte}$ so that $X_{i,j} = \rho_i(Y_j) + Z_{i,j}$,
$v_K(Y_j) \geq v_{K_i}(X_{i,j})$, and $v_{K_i}(Z_{i,j}) \geq v_{K_i}(X_{i,j})$.
Finally, set $V_{j+1} = V_j - Y_j$ and $W_{i,j+1} = W_{i,j} + Z_{i,j}$.
Then
\begin{align*}
v_{K_i}(U'_i \rho_i(V_{j+1}) - W_{i,j+1}) &= v_{K_i}(U'_i \rho_i(V_j - Y_j) - 
W_{i,j} - Z_{i,j}) \\
&= v_{K_i}((I- U'_i) \rho(Y_j)) \\
&\geq (j+1)h,
\end{align*}
and the iteration continues.

From the construction, we see that $v_{K_i}(W_{i,j+1} - W_{i,j}) \geq jh$.
Consequently each sequence $\{W_{i,j}\}_{j=1}^\infty$ 
converges to an invertible matrix $W_i$ over $\calR_i^{+,\inte}$.
The entries in the direct sum of the matrices $W_i^{-1} U'_i$ 
are in the closure of the image of $\rho$ on $b^{-1} A$.
By Proposition~\ref{prop:closed}, there is a matrix $V'$
over $b^{-1}A$ with $U'_i \rho_i(V') = W_i$ for all $i$.
We may thus take $V = PV'$.
\end{proof}

\subsection{General coefficients}

We now analyze the cohomology of $(\sigma, \nabla)$-modules on a curve,
focusing on the crucial case of unipotent local monodromy.
\begin{prop} \label{prop:unipcohom}
Let $M$ be a $(\sigma, \nabla)$-module over $A$.
Assume that $M \otimes \calR_i$ is unipotent for each $i$.
Then there is a commuting diagram
\begin{equation} \label{eq:diagram}
\xymatrix{ M \ar[r] \ar_{\nabla}[d] &
\oplus_i (\calR_i/\calR_i^+)^n \ar[d] \\ M \otimes
\Omega^1_{A/K} \ar[r] & \oplus_i (\Omega^1_{\calR_i/K_i} / t_i^{-1}
\Omega^1_{\calR_i^+/K_i})^n }
\end{equation}
in which the horizontal arrows and the right vertical arrow are
near-isomorphisms of $K$-vector spaces. 
Consequently, $\nabla: M \to M \otimes \Omega^1_{A/K}$
is a near-isomorphism.
\end{prop}
\begin{proof}
Since $M$ is locally free over $A$, we can 
choose a basis $\be_1, \dots, \be_n$ of $M$ over $A[a^{-1}]$ for some $a$.
(Note that $A[a^{-1}]$ is the localization in the category of
rings, not dagger algebras.)
Define the matrices $N_i$ over $\calR_i$
by $\nabla \be_l = \sum_{j} (N_i)_{jl} \be_j \otimes 
\frac{dt_i}{t_i}$. By the unipotence hypothesis
and Proposition~\ref{prop:unipbasis}, there exist matrices $U_i$
over $\calR_i$ such that $X_i = U_i^{-1} N_i U_i + U_i^{-1} t_i 
\frac{d}{dt_i} U_i$ is a nilpotent matrix over $K_i$.

By Lemma~\ref{lem:multnear}, there exists a matrix $V$ over $\Frac A$
such that $W_i = U_i^{-1} \rho_i(V)$ and its inverse have entries in
$\calR_i^+$ for $i=1, \dots, m$.
Put $\bv_l = \sum_j U_{jl} \be_j$ for $l=1,\dots,n$.
We now can describe the horizontal arrows in \eqref{eq:diagram} as follows.
Form the free module $S$ over $\oplus_i \calR_i$ generated by
$\bw_1, \dots, \bw_n$, map $M$ to this module by sending
$\sum_l c_l \bv_l$ to $\sum_l \rho(c_l) \bw_l$, and quotient as needed.

Put 
\[
N'_i = \rho(V)^{-1} N_i \rho(V) + \rho(V)^{-1} t_i \frac{d}{dt_i} \rho(V)
= W_i^{-1} X_i W_i + W_i^{-1} t \frac{d}{dt} W_i,
\]
which has entries in $\calR^+$. Then
the right vertical arrow in \eqref{eq:diagram} is induced by the map
$S \to S \otimes \oplus_i \Omega^1_{\calR_i/K_i}$ given by
\[
\sum_l c_l \bw_l \mapsto
\sum_l \bw_l \otimes dc_l + \sum_{j,l} (N'_i)_{jl} c_l \bw_j \otimes
\frac{dt}{t}.
\]
The right vertical arrow is actually an isomorphism, by a straightforward
calculation (or skip ahead to Proposition~\ref{prop:pluscohom}).
The top horizontal arrow can be factored as an inclusion
$M \to \frac{1}{u} A^n$ (for some $u \in A$) followed by $\rho$, which are
near-isomorphisms by Lemma~\ref{lem:frac} and Lemma~\ref{lem:near2},
respectively;
hence it is a near-isomorphism. Similarly, the bottom horizontal arrow
is a near-isomorphism by Lemma~\ref{lem:frac} and Lemma~\ref{lem:near2}.

Thus three of the arrows in \eqref{eq:diagram} are near-isomorphisms, so
the left vertical arrow is as well.
\end{proof}

\begin{theorem} \label{thm:fincurve}
Let $X$ be a smooth affine curve over $k$ and let $\calE$ be an overconvergent
$F$-isocrystal on $X$. Then $H^i(X/K, \calE)$ is finite dimensional over $K$
for all $i$.
\end{theorem}
\begin{proof}
Let $A$ be a dagger algebra of MW-type 
with special fibre $X$, let $\sigma: A \to A$
be a Frobenius lift, and let $M$ be a $(\sigma, \nabla)$-module on $A$
corresponding to $\calE$.
As in Proposition~\ref{prop:const}, by pulling back along a suitable
cover, we may reduce to the case where Hypothesis~\ref{hypo:curve}
holds; we may also reduce to the case where $M \otimes \calR_i$ is unipotent
for each $i$. But then Proposition~\ref{prop:unipcohom} asserts that
$H^0(M)$ and $H^1(M)$ are finite dimensional over $K$. This
yields the desired result.
\end{proof}

\begin{remark} \label{rem:cohcurve}
As mentioned earlier, 
we know of three other ways to prove Theorem~\ref{thm:fincurve}
(at least for $k$ perfect, but one can easily reduce to this case).
One proof is given by Crew \cite{bib:crew2},
who establishes finiteness of cohomology on a curve
using techniques of $p$-adic functional analysis (the key argument being
that a locally convex space
which is both Fr\'echet and dual-of-Fr\'echet is finite dimensional).
A second approach is to ``compactify coefficients'' as in \cite{bib:me8} and then
apply Kiehl's finiteness theorems \cite{bib:kiehl}, or put another way,
invoke a comparison between rigid and log-crystalline cohomology.
A third approach, closely related to the second, is to ``algebrize''
the connection and compare rigid cohomology to algebraic 
de~Rham cohomology \cite[Corollaire~5.0-12]{bib:cm4}.
\end{remark}

\begin{remark}
As noted at the beginning of this chapter,
one can actually deduce the existence of generic pushforwards
by relativizing the proof of Theorem~\ref{thm:fincurve}.
However, the argument is somewhat
messier than the one we use in the next chapter, and gives less precise
information. We thus leave it as an exercise for the interested
reader.
\end{remark}

\section{Pushforwards in rigid cohomology}
\label{sec:push}

In this section, we establish the existence of ``generic'' pushforwards,
with and without supports,
for an overconvergent $F$-isocrystal on a family 
of affine curves. A similar result for families of projective
curves has been given by Tsuzuki \cite{bib:tsu7}.

\subsection{Relative local cohomology}

The main relevance of a relative version of the $p$-adic local monodromy
theorem in this paper is that it controls the local cohomology of
a $(\sigma, \nabla)$-module over the Robba ring over a dagger algebra.
In this section, we illustrate how this control is exerted, starting
with a computation in the unipotent case.
\begin{lemma} \label{lem:freegen}
Let $A$ be a reduced dagger algebra
equipped with a Frobenius lift.
Let $M$ be a free unipotent $(\sigma, \nabla)$-module over
$\calR_A$, relative to $A$.
Then $H^0(M)$ and $H^1(M)$ are finitely generated $A$-modules.
\end{lemma}
\begin{proof}
By Proposition~\ref{prop:unipbasis},
there exists a strongly unipotent basis  $\bv_1, \dots, \bv_n$ of $M$.
In terms of this basis, it is easy to read off the kernel and cokernel
of $\nabla$.
\begin{itemize}
\item
Any element of the kernel of $\nabla$ must lie in the
$A$-span of the $\bv_i$. Hence $\ker(\nabla) = 
H^0(M)$ is a submodule
of a finitely generated $A$-module, and so is finitely generated.
\item
Any element of $M \otimes \Omega^1_{\calR_A/A}$ can be formally written as
\[
\sum_{i=1}^n \sum_{l \in \ZZ} f_{il} t^l \bv_i \otimes \frac{dt}{t}
\]
for some $f_{il} \in A$.
If $f_{i0} = 0$ for $i=1, \dots, n$, then this element is in the image of
$\nabla$. Hence $\coker(\nabla) = H^1(M)$ is a quotient
of a finitely generated $A$-module, and so is finitely generated.
\end{itemize}
\end{proof}

\begin{remark} \label{rem:freetop}
From the proof, we can easily read off two topological facts about a
unipotent $(\sigma, \nabla)$-module $M$ over $\calR_A$, for
$A$ a reduced dagger algebra.
\begin{enumerate}
\item[(a)] The image of $\nabla$ is closed.
\item[(b)] The continuous bijection $M/\ker(\nabla) \to \im(\nabla)$ is a
homeomorphism.
\end{enumerate}
These will later form the basis for more general statements of the same form.
\end{remark}

In the manner of Lemma~\ref{lem:freegen}, we can also verify the observation
that was made without proof in the course of proving
Theorem~\ref{thm:fincurve}.
\begin{prop} \label{prop:pluscohom}
Let $M$ be a unipotent $(\sigma, \nabla)$-module over $\calR_K$.
Let $\bv_1, \dots, \bv_n$ be a strongly unipotent basis of $M$,
and let $M_1$ be the
$\calR_K^+$-span of the $\bv_i$. Then the induced map
\[
\nabla: M/M_1 \to (M \otimes dt) /(M_1 \otimes \frac{dt}{t})
\]
is an isomorphism.
\end{prop}
\begin{proof}
The fact that the map is an isomorphism means two things, both of which
are evident from the special form of $\nabla$.
\begin{itemize}
\item If $\bv \in M$ satisfies $\nabla \bv \in M_1 \otimes \frac{dt}{t}$,
then $\bv \in M_1$. This follows because 
we may assume the coefficient of $t^j$ in $\bv$
is zero for $j \geq 0$, but then any power of $t$ which appears in $\bv$
also appears in $\nabla \bv$.
\item For every element $\omega$ of $M \otimes \frac{dt}{t}$, there
exists $\bv \in M$ such that $\nabla \bv - \omega \in M_1 \otimes dt/t$.
This follows because once we subtract off the coefficient of $dt/t$
in $\omega$, the result maps to zero in the cokernel of $\nabla:
M \to M \otimes dt$ as seen in the proof of Lemma~\ref{lem:freegen}.
\end{itemize}
\end{proof}

Now we consider relative local cohomology in the non-unipotent case.
\begin{prop} \label{prop:relloc}
Let $A$ be an integral dagger algebra of MW-type
equipped with a Frobenius lift.
Let $M$ be a free $(\sigma, \nabla)$-module over $\calR_A$,
relative to $A$,
and let $B$ be a localization satisfying the conclusion of
Theorem~\ref{thm:monodromy2}. Then $H^0(M \otimes \calR_B)$
and $H^1(M \otimes \calR_B)$ are finitely generated $B$-modules.
\end{prop}
\begin{proof}
An extension $\calR''$ as in Theorem~\ref{thm:monodromy2} is the
composition of two extensions as in Remark~\ref{rem:byhand}
(the first because $B^{\sigma^{-m}}$ is flat over $B$ by 
Remark~\ref{rem:frob}),
so $H^i(M \otimes \calR_B)$ injects into $H^i(M \otimes \calR'')$ for
$i=0,1$. Since $H^i(M \otimes \calR'')$ is finitely generated
by Lemma~\ref{lem:freegen}, we have the desired result.
\end{proof}

\begin{remark}
Over $K$, one can also deduce Proposition~\ref{prop:relloc} from
the $p$-adic index theorem of Christol-Mebkhout 
\cite[Th\'eor\`eme~7.5-2]{bib:cm3}.
We do not know whether that theorem admits a relative version from
which one can deduce Proposition~\ref{prop:relloc} in general.
\end{remark}

\subsection{Base change in local cohomology}

We next examine how the local pushforwards we have constructed behave with 
respect to flat base change.
\begin{prop} \label{prop:locbasechange}
Let $A$ be an integral dagger algebra of MW-type
equipped with a Frobenius lift.
Let $M$ be a free $(\sigma, \nabla)$-module over $\calR_A$,
relative to $A$,
and let $B$ be a localization satisfying the conclusion of
Theorem~\ref{thm:monodromy2}.
Let $C$ be an integral dagger algebra of MW-type which is flat over $B$.
Then the base change maps
$H^i(M \otimes \calR_B) \otimes_B C \cong H^i(M \otimes
\calR_C)$ (induced by the multiplication map $\calR_B \otimes C \to \calR_C$)
are isomorphisms. Moreover, if we can take $B=A$ and
$\calR'' = \calR_A$ (i.e.,
$M$ is unipotent), the
flatness hypothesis on $C$ can be dropped.
\end{prop}
\begin{proof}
We first handle the case where $M$ is unipotent, i.e., in case we can take
$B=A$
and $\calR'' = \calR_A$. Let
$\bv_1, \dots, \bv_n$ be a unipotent basis, and let $M_j$ be the span
of $\bv_1, \dots, \bv_j$. Then $\nabla$ carries $M_j$ to $M_j \otimes
\Omega^1_{\calR_A/A}$ and $M_j/M_{j+1}$ is isomorphic, as a module
with connection, to $\calR_A$ itself. Thus $H^i(M_j/M_{j+1})$
is a free $A$-module of rank 1 for $i=0,1$, and the map
\[
H^i(M_j/M_{j+1}) \otimes_A C \to H^i((M_j/M_{j+1}) \otimes \calR_C)
\]
is clearly an isomorphism. We deduce by induction on $j$ and the snake lemma
that $H^i(M_j) \otimes_A C \to H^i(M_j \otimes \calR_C)$ is an isomorphism,
and the case $j=n$ gives what we wanted.

We now treat the general case.
Set notation as in Theorem~\ref{thm:monodromy2}, and
put $\calR''_C = \calR'' \otimes_{\calR_B} \calR_C$; we then have a
commuting diagram
\begin{equation} \label{eq:basechange}
\xymatrix{
H^i(M \otimes \calR'') \otimes_B C
 \ar[r] \ar[d] & H^i(M \otimes \calR'') \otimes_B C\ar[d] \\
H^i(M \otimes \calR''_C) \ar[r] & H^i(M \otimes \calR''_C)
}
\end{equation}
in which the horizontal arrows are projectors induced by the trace map.
Since the left vertical arrow is an isomorphism, it induces an isomorphism
on the images of the projectors. The image of the 
lower horizontal arrow is of course
$H^i(M \otimes \calR_C)$. As for the upper horizontal arrow, since $C$
is flat over $B$, tensoring with it commutes with taking images (since
the image of $f: X \to Y$ is the kernel of $Y \to \coker(f)$, and tensoring
with $C$ commutes with taking kernels and cokernels).
Since the image of $\Trace$ on $H^i(M \otimes \calR'')$ is $H^i(M
\otimes \calR_B)$,
the image of the upper horizontal arrow is thus $H^i(M \otimes
\calR_B) \otimes_B C$.
We conclude that \eqref{eq:basechange} yields the desired
base change isomorphism.
\end{proof}

\begin{remark} \label{rem:flatfrac}
The same argument works with $C$ replaced
by the completion of $\Frac B$ for the affinoid topology.
\end{remark}

\subsection{Generic pushforwards}

We now exhibit ``generic'' higher direct images in a simple geometric
setting, using the formalism of \cite[8.1]{bib:crew2}.

\begin{hypo}
Let $A$ be a dagger algebra of MW-type with special fibre $X$.
Let $f: A \to A \langle x \rangle^\dagger$ be the canonical inclusion,
and embed $A \langle x \rangle^\dagger$ into $\calR_A$ by mapping
$\sum c_i x^i$ to $\sum c_i t^{-i}$.
\end{hypo}

\begin{defn}
Let $M$ be a $(\sigma, \nabla)$-module over $A \langle x \rangle^\dagger$,
and let $\nabla_v$ denote the composition
\[
M \stackrel{\nabla}{\to} M \otimes \Omega^1_{A \langle x \rangle^\dagger/K}
\to M \otimes \Omega^1_{A \langle x \rangle^\dagger/A} = M\otimes dx.
\]
Let
\begin{align*}
\nabla_v^{\loc} &: M \otimes \calR_A \to M \otimes \Omega^1_{\calR_A/A} \\
\nabla_v^{\qu} &: (M \otimes \calR_A)/M \to
(M \otimes \Omega^1_{\calR_A/A})/(M \otimes \Omega^1_{A
\langle x \rangle^\dagger/A})
\end{align*}
be the maps induced by $\nabla_v$.
We then define the following $A$-modules:
\begin{gather*}
R^0 f_* M = \ker(\nabla_v), \qquad
R^1 f_* M = \coker(\nabla_v) \\
R^0_{\loc} f_* M = \ker(\nabla_v^{\loc}), \qquad
R^1_{\loc} f_* M = \coker(\nabla_v^{\loc}) \\
R^1 f_! M = \ker(\nabla_v^{\qu}), \qquad
R^2 f_! M = \coker(\nabla_v^{\qu}).
\end{gather*}
The snake lemma applied to the diagram with exact rows
\begin{equation} \label{eq:presnake}
\xymatrix{
0 \ar[r] & M \ar^{\nabla_v}[d] \ar[r] & M \otimes \calR_A
\ar^{\nabla_{\loc}}[d] \ar[r] & (M \otimes \calR_A)/M \ar^{\nabla_{\qu}}[d]
\ar[r] & 0 \\
0 \ar[r] & M \otimes \Omega^1_{A \langle x \rangle^\dagger/A}
\ar[r] & M \otimes \Omega^1_{\calR_A/A} 
\ar[r] & (M \otimes \Omega^1_{\calR_A/A})/(M \otimes \Omega^1_{A\langle x \rangle^\dagger/A}) \ar[r] & 0
}
\end{equation}
gives an exact sequence
\begin{equation} \label{eq:snake}
0 \to R^0 f_* M \to R^0_{\loc} f_* M \to R^1 f_! M \to
R^1 f_* M \to R^1_{\loc} f_* M \to R^2 f_! M \to 0
\end{equation}
using which we can define the ``primitive middle cohomology'':
\[
R^1_p f_* M = \im(R^1 f_! M \to R^1 f_* M).
\]
\end{defn}

The Frobenius on $M$ acts on all of these modules.
Moreover, if we write $\nabla_h$ for the composition
\[
M \stackrel{\nabla}{\to} M \otimes \Omega^1_{A \langle x \rangle^\dagger/K}
\to M \otimes (A \langle x\rangle^\dagger \otimes_A \Omega^1_{A/K}) =
M \otimes_A \Omega^1_{A/K},
\]
then each module $N$ listed above admits a $K$-linear connection
$\nabla_h: N \to N \otimes_A \Omega^1_{A/K}$. Note that this is true
even if $N$ is not finitely generated; that is because
$\Omega^1_{A/K}$ is projective over $A$, so computing the kernel/cokernel
of any of $\nabla_v, \nabla_v^{\loc}, \nabla_v^{\qu}$ commutes with tensoring
with $\Omega^1_{A/K}$.

We now prove our main theorem on higher direct images. Note the
essential use of the nonrelative case (Theorem~\ref{thm:fincurve}).
\begin{theorem} \label{thm:fingen}
Let $A$ be an integral dagger algebra of MW-type, and
let $f: A \to A \langle x \rangle^\dagger$ denote the canonical inclusion.
Let $M$ be a $(\sigma,\nabla)$-module over $A \langle x \rangle^\dagger$.
Then after replacing $A$ by a suitable localization, the modules
$R^i f_* M, R^i_{\loc} f_* M, R^i f_! M, R^1_p f_* M$ 
become $(\sigma, \nabla)$-modules
over $A$.
\end{theorem}
\begin{proof}
Choose a nonzero $g \in A \langle x \rangle^\dagger$ such that $M$ becomes
free over $A[g^{-1}]$. Write $g = \sum_{i=0}^\infty g_i x^i$ and choose the
largest $i$ for which $|g_i|$ is maximized. By Lemma~\ref{lem:nilpotent},
the series expansion
of $(1 - (1-g/(g_i x^i)))^{-1}$ converges in $\calR_{A\langle g_i^{-1} \rangle^\dagger}$;
hence for any localization
$B$ of $A \langle g_i^{-1} \rangle^\dagger$, $g$ is invertible in
$\calR_B^{\inte}[\frac 1p]$. In particular, $M \otimes \calR_B$ is free
for such $B$, and so we can choose $B$ satisfying the conclusion of
Theorem~\ref{thm:monodromy2}.

We first verify that over $B$,
$R^i f_* M$, $R^i_{\loc} f_* M$, $R^i f_! M$, $R^1_p f_* M$ 
become finitely generated $B$-modules.
We have that $R^0_{\loc} f_* M$ and $R^1_{\loc} f_* M$ are
finitely generated by Proposition~\ref{prop:relloc}; by
the exactness of \eqref{eq:snake}, it suffices to check that
$R^1_p f_* M$ is finitely generated.

Let $L$ be the completion of $\Frac A$ for the affinoid topology
and put $M_L = M \otimes L \langle x \rangle^\dagger$. By
Theorem~\ref{thm:fincurve} and Proposition~\ref{prop:relloc},
respectively, $R^1 f_* M_L$ and $R^0_{\loc} f_* M_L$ are finite
dimensional over $L$; by the snake lemma, so then is $R^1 f_! M_L$.
Now view $(M \otimes \calR_A)/ M$ and $(M_L \otimes \calR_L)/M_L$ as 
$A$-submodules of direct
products of copies of $A$ and $L$, respectively, indexed by pairs consisting
of a positive power of $t$ and a basis element of $M$. 
In particular, we may view $(M \otimes \calR_A)/ M$
as an $A$-submodule of $(M_L \otimes \calR_L)/M_L$; in this point of view,
$R^1 f_! M$ is none other than the intersection of $R^1 f_! M_L$ with
$(M \otimes \calR_A)/M$ within $(M_L \otimes \calR_L)/M_L$. Since $A$
is noetherian, we may apply Lemma~\ref{lem:domain} below to
deduce that $R^1 f_! M$ is finitely generated, as then is $R^1_p f_* M$.

For $N$ any one of the modules under consideration, we now know that $N$
is finitely generated. Since $N$ is equipped with a connection given by
$\nabla_h$, $N$ is also locally free over $B$ by Lemma~\ref{lem:locfree}.
We also have an action of $F$ on $N$, which is injective by 
Remark~\ref{rem:injfrob}. 
The quotient of $N$ by the submodule generated by the image of $F$ also
comes equipped with a connection, so it is locally free. However,
the rank of this quotient must be zero, so the quotient itself must be zero.
Hence the image of $F$ generates $N$; that is, the induced $B$-linear
map $\sigma_A^* N \to N$ is surjective. Now the kernel of this map
is equipped with a connection, so is locally free by
Lemma~\ref{lem:locfree} and necessarily of rank zero,
hence trivial. Thus $\sigma_A^* N \to N$ is an isomorphism.

We conclude that all of the modules in question become
$(\sigma, \nabla)$-modules over $B$, as desired.
\end{proof}

\begin{lemma} \label{lem:domain}
Let $R$ be a noetherian domain, let $L$ be a field containing $R$,
and let $S$ be an arbitrary set. Let $N$ be a finite dimensional $L$-subspace
of the direct product $L^S$. Then $N \cap R^S$ is a finitely generated 
$R$-module.
\end{lemma}
\begin{proof}
We induct on the dimension of the $L$-span of $N \cap R^S$. 
If this dimension is zero, then
$N \cap R^S = 0$ and there is nothing to check. Otherwise, we may 
pick $s \in S$
such that, writing $\pi_s$ for the projection of $R^S$ onto the $s$-component,
we have $\pi_s(N \cap R^S) \neq 0$. Now $\pi_s(N \cap R^S)$ is an ideal of
$R$, which by the noetherian hypothesis must be finitely generated. Pick
elements $\bv_1, \dots, \bv_n \in N \cap R^S$ whose images under $\pi_s$
generate $\pi_s(N \cap R^S)$. Then each element of $N \cap R^S$ can be
written as an $R$-linear combination of $\bv_1, \dots, \bv_n$ plus an element
of $(N \cap \pi_s^{-1}(0)) \cap R^S$, so it suffices to check that the
latter is finitely generated. But the latter has the same form as the
original intersection except with $N$ replaced by $N \cap \pi_s^{-1}(0)$, and
the dimension of the $L$-span of $(N \cap \pi_s^{-1}(0)) \cap R^S$ has gone
down by one since $\pi_s(N \cap R^S) \neq 0$. Thus the induction hypothesis
yields the desired result.
\end{proof}

\begin{remark}
By passing down an \'etale map to this geometric situation (using
Proposition~\ref{prop:belyi}), one can use these methods to construct
``generic'' higher direct images for a family of affine curves.
(By contrast, Tsuzuki \cite{bib:tsu5} constructs generic higher direct
images for families of projective curves, using more global techniques.)
On the other hand, 
pushforwards with compact supports are more difficult to construct, and
we do not believe our technique is sophisticated
enough to handle them. In particular, in the ``non-generic'' setting,
the objects we call $R^i f_! M$ may not
agree with the direct image with compact supports in the category of arithmetic
$\calD$-modules.
\end{remark}

\begin{remark}
Although our argument essentially uses the presence of the Frobenius structure,
in order to invoke the $p$-adic local monodromy theorem, it is conceivable
that this dependence could be lifted. 
To do so, it may be necessary to extend the $p$-adic index theory of
Christol and Mebkhout \cite{bib:cm4} to cover modules with connection
over the Robba ring over a dagger algebra.
\end{remark}

\subsection{A Leray spectral sequence}
\label{subsec:leray}

The value of the pushforward construction
is that the cohomology of a $(\sigma, \nabla)$-module
can be controlled using the cohomology of its pushforwards, as follows.

\begin{prop} \label{prop:leray}
Let $A$ be a dagger algebra, and
let $M$ be a $(\sigma, \nabla)$-module over $A \langle x \rangle^\dagger$.
Let $P$ and $Q$ be the kernel and cokernel, respectively, of
\[
\nabla_v: M \stackrel{\nabla}{\to} M \otimes_{A \langle x \rangle^\dagger} 
\Omega^1_{A \langle x \rangle^\dagger /K}
\to M \otimes_{A \langle x \rangle^\dagger} \Omega^1_{A \langle x 
\rangle^\dagger/A}.
\]
Then there is a long exact sequence
\[
\cdots \to H^i(P) \to H^i(M) \to H^{i-1}(Q) \to H^{i+1}(P) \to \cdots.
\]
\end{prop}
\begin{proof}
The claim is a formal consequence of the Leray spectral sequence
associated to the double complex
\[
E_0^{pq} = M \otimes_A \Omega^p_{A/K} \otimes_{K \langle x
\rangle^\dagger} \Omega^q_{K \langle x\rangle^\dagger/K}
\]
whose total cohomology computes $H^i(M)$.
(We follow \cite[Section~14]{bib:bt} 
in our terminology concerning spectral sequences.)
Using the flatness of $\Omega^i_{A/K}$ over $A$, we have
\begin{gather*}
E_1^{p0} = P \otimes \Omega^p_{A/K},
E_1^{p1} = Q \otimes \Omega^p_{A/K} \\
E_2^{p0} = H^p(P), E_2^{p1} = H^p(Q),
\end{gather*}
the differential $d_2$ maps $E_2^{pq}$ to $E_2^{(p+2)(q-1)}$, and
the sequence degenerates at $E_3$ (that is, $d_3=d_4 = \cdots =0$).
Thus in the desired long exact sequence, the maps $H^{i-1}(Q)
\to H^{i+1}(P)$ are given by $d_2$; the remaining maps come from
the fact that $E_3$ abuts to $H^i(M)$.
\end{proof}

\subsection{Pushforwards and base change}

To establish the K\"unneth formula, we need to prove that our relative
pushforwards commute with flat base change. The argument requires a
bit of $p$-adic functional analysis; since this will be introduced in 
Chapter~\ref{sec:push2}, the proof will contain a forward reference, but
it can be verified that no vicious circle has been created.

\begin{lemma} \label{lem:closedim}
Let $A$ be an integral dagger algebra of MW-type, and
let $M$ be a $(\sigma, \nabla)$-module over $A \langle x \rangle^\dagger$,
relative to $A$.
Then there exists a localization $A'$ of $A$
for which the image of $\nabla$ on $M \otimes_{A \langle x \rangle^\dagger}
A' \langle x \rangle^\dagger$ is closed for the fringe topology.
\end{lemma}
\begin{proof}
This is a corollary of Proposition~\ref{prop:strict}. Note that the proposition
includes the conclusion that $\nabla$ is strict; although we do not
need that conclusion here, its presence is a key aspect of the proof.
\end{proof}

We can now establish a flat base change theorem.
\begin{prop} \label{prop:basechange}
Let $A$ and $B$ be integral dagger algebras of MW-type,
with $B$ flat over $A$.
Let $f: A \to A \langle x \rangle^\dagger$ and $g: B \to 
B \langle x \rangle^\dagger$ be the canonical inclusions.
Let $M$ be a $(\sigma,\nabla)$-module over $A \langle x \rangle^\dagger$,
and put $M_B = M \otimes_{A \langle x \rangle^\dagger} B \langle x 
\rangle^\dagger$.
Suppose that the conclusion of Theorem~\ref{thm:fingen} holds
for $M$
without further localization. Then the canonical base change morphisms
\begin{align*}
(R^i f_* M) \otimes_A B &\to R^i g_* M_B \\
(R^i_{\loc} f_* M) \otimes_A B &\to R^i_{\loc} g_* M_B \\
(R^i f_! M) \otimes_A B &\to R^i g_! M_B \\
(R^1_p f_* M) \otimes_A B &\to R^1_p g_* M_B
\end{align*}
are isomorphisms. In particular, the conclusion of Theorem~\ref{thm:fingen}
holds for $M_B$ without further localization.
\end{prop}
\begin{proof}
For $R^i_{\loc}$, this follows from Proposition~\ref{prop:locbasechange}.
As in Theorem~\ref{thm:fingen}, we have that
$R^1_p g_* M_B = \im(R^1 g_! M_B \to R^1 g_* M_B)$ is finitely generated,
which shows (by the exactness of \eqref{eq:snake})
that all of the pushforwards are finitely generated without
further localization.

It is clear that $(R^0 f_* M) \otimes_A B \to R^0 g_* M_B$ and
$(R^1 f_! M) \otimes_A B \to R^1 g_! M_B$ are injective, since
the map of complexes $M \otimes 
\Omega^._{A \langle x \rangle^\dagger/A} \otimes_A B
\to M_B \otimes \Omega^._{B \langle x \rangle^\dagger/B}$ is injective and 
$R^0 f_* M$ and $R^1 f_! M$ are defined by taking certain kernels of $\nabla$.

We next show that $(R^1 f_* M) \otimes_A B \to R^1 g_* M_B$ is surjective; since
source and target are $\nabla$-modules, 
surjectivity may be checked after
localizing $B$ (and this localization preserves flatness over $A$).
In particular, we may reduce to the case where the
conclusion of Lemma~\ref{lem:closedim} holds for
$N = M_B \otimes \Omega^1_{B \langle x \rangle^\dagger/B}$ without further
localization.
For brevity, put $N_1 =
M \otimes \Omega^1_{A \langle x \rangle^\dagger/A} \otimes_A B$,
let $N_2$ be the image of $\nabla$ on $M_B$.
Then $N_1$ is dense in $N$, and $N_2$ is closed in $N$ by 
Lemma~\ref{lem:closedim}. Since $B$ is noetherian and $R^1 g_* M_B = N/N_2$ is
finitely generated over $B$, 
the image of $N_1$ in $N/N_2$ is also finitely
generated. Thus $N_1 + N_2$ is  closed in $N$; however, it is also dense
because $N_1$ already is. Consequently $N_1 + N_2 = N$, and we deduce that
$(R^1 f_* M) \otimes_A B \to R^1 g_* M_B$ is surjective.

To conclude, write down the snake lemma exact sequences \eqref{eq:snake}
corresponding to $M$ and $M_B$. Then the base change maps give a map
between the former, tensored with $B$, and the latter, and we are
in the situation of the following Lemma~\ref{lem:snake}. It follows
that all of the base change maps are isomorphisms.
\end{proof}

\begin{lemma} \label{lem:snake}
Suppose in the category of abelian groups, the diagram
\[
\xymatrix{
0 \ar[r] & A \ar[r] \ar^a[d] & B \ar[r] \ar^b[d] & C \ar[r] \ar^c[d]
& D \ar[r] \ar^d[d] & E \ar[r] \ar^e[d] & F \ar[r] \ar^f[d] & 0 \\
0 \ar[r] & A' \ar[r] & B' \ar[r] & C' \ar[r] & D' \ar[r] & E' \ar[r]
& F' \ar[r] & 0
}
\]
commutes and has exact rows. Also suppose that $c$ is injective,
$b$ and $e$ are isomorphisms, and $d$ is surjective. Then all six
vertical arrows are isomorphisms.
\end{lemma}
\begin{proof}
Apply the five lemma first to $c$ and $d$, then to $a$ and $f$.
\end{proof}

\begin{remark} \label{rem:flatfrac2}
As in Remark~\ref{rem:flatfrac}, we can also take $B$ to be
the completion of $\Frac A$ under the affinoid topology. This
implies that the pushforwards all have the ``correct'' ranks, namely
the ranks over the generic fibre.
\end{remark}

\subsection{Pushforwards and external products}

Using flat base change, we can now verify that our pushforward construction
commutes with external products. This will be crucial for verifying the
K\"unneth decomposition.

\begin{lemma} \label{lem:extpush}
Let $A$ be an integral dagger algebra of MW-type,
and let $f: A \to A \langle x \rangle^\dagger$ be the natural inclusion.
Let $M$ be a $(\sigma, \nabla)$-module over $A \langle x \rangle^\dagger$,
and let $N$ be a $(\sigma, \nabla)$-module over $A$. Then the natural maps
\begin{align*}
(R^i f_* M) \otimes_A N &\to R^i f_* (M \otimes_A N) \\
(R^i_{\loc} f_* M) \otimes_A N &\to R^i_{\loc} f_* (M \otimes_A N) \\
(R^i f_! M) \otimes_A N &\to R^i f_! (M \otimes_A N) \\
(R^1_p f_* M) \otimes_A N &\to R^1_p f_* (M \otimes_A N)
\end{align*}
are isomorphisms.
\end{lemma}
\begin{proof}
All of the pushforwards are computed as kernels and cokernels of various
$A$-module homomorphisms. Since $N$ is flat over $A$, these computations
commute with tensoring with $N$.
\end{proof}

Given modules $M$ and $N$ over dagger algebras $A$ and $B$, respectively,
define the external product of $M$ and $N$ as the
$(A \hattimes B)$-module
\[
M \boxtimes_{A,B} N = M \otimes_A (A \hattimes B) \otimes_B N.
\]

\begin{prop} \label{prop:pushext}
Let $A$ and $B$ be integral dagger algebras of MW-type.
Let $M$ be a $(\sigma, \nabla)$-module over $A \langle x \rangle^\dagger$ 
and let $N$ be a $(\sigma, \nabla)$-module over $B$.
Let $f: A \to A \langle x \rangle^\dagger$
and $g: A \hattimes B \to (A \hattimes B) \langle x \rangle^\dagger$
be the natural inclusions.
Suppose that the conclusion of Theorem~\ref{thm:fingen} holds for $M$
without further localization. Then the natural maps
\[
(R^i f_* M) \boxtimes_{A,B} N \to
R^i g_* (M \boxtimes_{A \langle x \rangle^\dagger,B} N)
\]
are isomorphisms.
\end{prop}
Again, the truth of Theorem~\ref{thm:fingen} for 
$M \boxtimes_{A \langle x \rangle^\dagger,B} N$ without further
localization is part of the conclusion.
\begin{proof}
Note that $A \hattimes B$ is flat over $A$ by Lemma~\ref{lem:flat}
(see Remark~\ref{rem:flatprod}).
Hence we can apply Proposition~\ref{prop:basechange} to deduce that
the natural maps
\[
(R^i f_* M) \boxtimes_{A,B} N
\to (R^i g_* (M \otimes_{A \langle x \rangle^\dagger} (A \hattimes B)\langle
x \rangle^\dagger)) \otimes_{B} N
\]
are isomorphisms.
By Lemma~\ref{lem:extpush}, the natural maps
\[
(R^i g_* (M \otimes_{A \langle x \rangle^\dagger} (A \hattimes B)\langle
x \rangle^\dagger)) \otimes_{B} N
\to R^i g_* (M \boxtimes_{A \langle x \rangle^\dagger,B} N)
\]
are also isomorphisms. Putting these isomorphisms together yields
the desired result.
\end{proof}

\section{Cohomology with compact supports}
\label{sec:push2}

In this chapter, we review Berthelot's construction of cohomology
with compact supports, then make this construction explicit in the case
of affine space. The resulting computation has a flavor similar to that
of the Monsky-Washnitzer cohomology we have been using so far; it also
bears a passing resemblance to Dwork's ``dual theory''.
We then verify Poincar\'e duality on affine space; this requires
some of the machinery of $p$-adic functional analysis, for which we
draw from Crew \cite{bib:crew2}.

\subsection{The cohomology of affine space}
\label{subsec:cohomaff}

We now restrict consideration to a simple geometric situation.
\begin{hypo} \label{hypo:affine}
Put $X = \AAA^n_k$ and $\overline{X} = \PP^n_k$, and let
$P$ be the formal completion of $\PP^n_{\calO}$
along $\PP^n_k$. We can identify $\tilP$ with the set of 
tuples $(z_0,\dots,z_n)$ over $\calO^{\alg}$ (the integral
closure of $\calO$ in the algebraic closure $K^{\alg}$ of $K$)
with $\max\{|z_i|\}=1$ modulo the action of $(\calO^{\alg})^*$
and of $\Gal(K^{\alg}/K)$. We use the Frobenius lift $\sigma$
sending $(z_0, \dots, z_n)$ to $(z_0^q, \dots, z_n^q)$.
\end{hypo}

\begin{defn}
For $i=0, \dots, n$,
let $U_i$ be the subspace of $\tilP$ on which $|z_i| = 1$. Then
each $U_i$ is an affinoid subdomain of $\tilP$, and
\[
\Gamma(U_i) = K \left\langle \frac{z_0}{z_i}, \dots, \frac{z_{i-1}}{z_i},
\frac{z_{i+1}}{z_i}, \dots, \frac{z_n}{z_i} \right\rangle.
\]
For $S \subseteq \{0, \dots, n\}$, put $U_S = \cap_{i\in S} U_i$.
\end{defn}

\begin{defn}
For $\rho>1$ in the norm group of $K^{\alg}$,
let $V_\rho$ be the subspace of $\tilP$ on which $|z_0| \geq \rho^{-1}$.
Then $V_\rho$ is a strict neighborhood of $]X[$, and is an affinoid:
namely, if we put $x_i = z_i/z_0$, then $V_\rho
= \Maxspec T_{n,\rho}$.
\end{defn}

We now want to compute $\RR\underline{\Gamma}(\calF)$ 
(as defined in Definition~\ref{def:supports}) for $\calF$ the
coherent sheaf associated to a finite locally free $T_{n,\rho}$-module
$E$, viewed as a one-term complex in degree 0.
Our model for the computation is \cite[Exemple~3.4]{bib:ber1},
but our situation is complicated by the fact that
$V_\rho \setminus U_0$ is not quasi-Stein if $n > 1$. However, it admits a nice
cover by quasi-Stein varieties. (This fact is alluded to in the proof
of \cite[Proposition~1.1]{bib:ber6}.)

\begin{defn}
For $T \subseteq \{1, \dots, n\}$, let $W_T$ be the subspace of $V_\rho$
on which $|x_i| > 1$ for $i \in T$.
\end{defn}
\begin{lemma} \label{lem:quasistein}
The rigid analytic space $W_{T}$ is quasi-Stein for each $T$.
\end{lemma}
\begin{proof}
For $\eta \in (1,\rho) \cap \Gamma^*$,
let $X_\eta$ be the subset of $W_T$ on which
$\eta \leq |x_i| \leq \rho$ for $i \in T$.
Then $X_\eta$ is an affinoid, $X_\eta \subseteq X_{\tau}$ if $\eta > \tau$,
and $W_T$ is the union of the $X_\eta$. The density of
$\Gamma(X_\tau)$ in $\Gamma(X_\eta)$ is
straightforward, as in fact the 
Laurent polynomial ring $K[x_1,\dots,x_n][x_i^{-1}: i \in T]$ is dense
in $\Gamma(X_\eta)$. Thus $W_T$ is quasi-Stein, as claimed.
\end{proof}

We do not have names for the coordinate rings of the $W_T$, so 
we introduce them now.
\begin{defn}
Put $\calS^{(n)}_r = \calR^{\{t_1,\dots,t_n\}}_{K,r}$.
For $r = \log_p(\rho)$, we may
identify $T_{n,\rho}$ as a subring of $\calS^{(n)}_r$ by equating
$\sum_I c_I x^I$ with $\sum_I c_I t^{-I}$.
For $T \subset \{1, \dots, n\}$,
let $\calS^{(n),T}_r$ denote the subring of the ring $\calS^{(n)}_r$
of series in which the variable $t_i$ may occur with a positive power
if and only if $i \in T$.
For short, write $i'$ for the set $\{1, \dots, n\} \setminus \{i\}$.
Observe that $\Gamma(W_T) = \calS^{(n),T}_r$ and
that $\calS^{(n),\emptyset}_r = T_{n,\rho}$.
\end{defn}

Recall that the complex $\RR \underline{\Gamma}(\calF)$ is isomorphic, in the
derived category, to the two-term complex
\begin{equation} \label{eq:twoterm}
0 \to \calF \to \iota_* \iota^* \calF \to 0.
\end{equation}
We now compute the (hyper)cohomology of this complex.
\begin{prop} \label{prop:berthcomplex}
Let $\calF$ be the coherent sheaf associated to the finite locally free
$T_{n,\rho}$-module $E$. Then the cohomology of \eqref{eq:twoterm}
is equal to 
\[
\frac{E \otimes \calS^{(n)}_r}{\sum_i 
E \otimes \calS^{(n),i'}_r}
\]
in degree $n$ and zero elsewhere.
\end{prop}
Note in particular that this space does not depend on $r$.
\begin{proof}
We use the admissible covering $U_0, U_1 \cap V_\rho, \dots, U_n \cap V_\rho$
of $V_\rho$ to resolve the two terms of \eqref{eq:twoterm} 
with \v{C}ech complexes. For each nonempty subset $T$ of
$\{1, \dots, n\}$, we have
\[
(V_\rho \setminus U_0) \cap \bigcap_{j \in T} U_j = W_T,
\]
which by Lemma~\ref{lem:quasistein} is quasi-Stein; also, each
$U_S$ is affinoid and hence quasi-Stein.
Hence (by Theorem B again) we may replace the \v{C}ech complexes
with their global sections. In other words,
the complex \eqref{eq:twoterm} is isomorphic in the derived category
to the simple complex associated to the double complex with rows
\begin{equation} \label{eq:dblcomplex}
\bigoplus_i \left(\oplus_{S} \Gamma(\calF, U_{S})\right) \to 
\bigoplus_i \left( \oplus_{T}
\Gamma(\calF, W_{T}) \right),
\end{equation}
where $S$ runs over all subsets of $\{0,1,\dots, n\}$ of size
$i$ and $T$ runs over all subsets of $\{1, \dots, n\}$ of size $i$.

We now compute the cohomology of \eqref{eq:dblcomplex} using the
Leray spectral sequence (which already arose in Section~\ref{subsec:leray}).
Half of the computation of the $E_1$ term is easy: 
we have $E_1^{00} = E$ and $E_1^{0q} = 0$ for $q>0$.
To compute the rest, we use the fact that
the $j$-th cohomology of the \v{C}ech complex with $i$-th term
\[
\bigoplus_i \left( \oplus_T \calS^{(n),T}_r \right)
\]
is given by
\[
H^j = \begin{cases} T_{n,\rho} & j= 0 \\
\calS^{(n)}_r / \sum_i \calS^{(n),i'}_r & j=n-1\\
0 & \mbox{otherwise.}
      \end{cases}
\]
This is straightforward to prove by induction on $n$, if notationally
cumbersome; we leave details to the reader.
The upshot is that $E_1^{10} = E$, $E_1^{1(n-1)} = (E \otimes \calS^{(n)}_r) / 
\sum_i (E \otimes \calS^{n,i'}_r)$
and $E_1^{1q} = 0$ for $q \neq 0,n-1$.

Recall that $E_2^{pq}$ is the $p$-th cohomology of the complex
$E_1^{\cdot q}$. For $q=0$, the complex is just $E \to E$ via
the identity map, so $E_2^{p0} = 0$ for all $p$. The only
other term left is $E_2^{1(n-1)} = E_1^{1(n-1)}$. Thus the spectral sequence
degenerates at $E_2$, yielding the desired result.
\end{proof}

We can now compute the cohomology with compact supports of an overconvergent
$F$-isocrystal on $\AAA^n$.
\begin{prop} \label{prop:affcohsupp}
Let $\calE$ be an overconvergent $F$-isocrystal on $\AAA^n$ which
corresponds to a $(\sigma, \nabla)$-module $M$ over $T_{n,\rho}$.
Then
\[
H^{n+j}_{c,\rig}(\AAA^n/K, 
\calE) = H^{j}((M \otimes \Omega^. \otimes \calS^{(n)}_r) 
/ \oplus_i (M \otimes \Omega^. \otimes \calS^{(n),i'}_r)).
\]
\end{prop}
\begin{proof}
The cohomology is again computed by a Leray spectral sequence.
This time, by Proposition~\ref{prop:berthcomplex}, we have 
\[
E_1^{pn} = (M \otimes \Omega^p \otimes \calS^{(n)}_r) 
/ \oplus_i (M \otimes \Omega^p \otimes \calS^{(n),i'}_r)
\]
and $E_1^{pq} = 0$ for $q \neq n$.
Thus the spectral sequence degenerates at $E_2$ and yields the
desired result.
\end{proof}
\begin{cor} \label{cor:afftop}
The space $H^{2n}_{c,\rig}(\AAA^n/K)$ is one-dimensional over $K$.
\end{cor}
\begin{proof}
It is clear from the description in Proposition~\ref{prop:affcohsupp}
that $H^{2n}_{c,\rig}(\AAA^n/K)$ is generated by the class of
$dt_1/t_1 \wedge \cdots \wedge dt_n/t_n$.
\end{proof}

\subsection{The Poincar\'e pairing}
\label{subsec:poinpair}

We next recall the formalism of the Poincar\'e pairing, and explicate this
formalism over affine space. Our discussion is drawn from \cite{bib:ber6}.

Retain Hypothesis~\ref{hypo:rigsetup}, and assume now that $X$ has
dimension $n$.
Let $j$ denote the inclusion of $]X[$ into a strict neighborhood $V$.
As in \cite[Lemme~2.1]{bib:ber6}, for any sheaves $\calF, \calG$
of $K$-vector spaces on $V$, one has an isomorphism
\[
\Hom(j^\dagger  \calF, \calG) \cong \underline{\Gamma}(\Hom_K(\calF,\calG)).
\]
If $\calF, \calG$ are overconvergent $F$-isocrystals, we thus obtain
a pairing in cohomology
\begin{equation} \label{eq:cupprod}
H^j_{\rig}(X/K,\calF) \otimes_K H^l_{c,\rig}(X/K, \calG)
\to H^{j+l}_{c,\rig}(X/K, \calF \otimes \calG).
\end{equation}
In particular, if $\calG = \calF^\dual$ is the dual of $\calF$, we have
a natural 
morphism from $\calF \otimes \calF^\dual$ to the constant $F$-isocrystal.
This yields a map
\[
H^{2n}_{c,\rig}(X/K, \calF \otimes \calF^\dual) \to H^{2n}_{c,\rig}(X/K).
\]
One can now construct a trace morphism from $H^{2n}_{c,\rig}(X/K)$
to $K$.
In \cite[1.2]{bib:ber6}, this is done by a reduction to the corresponding
construction in coherent cohomology; one then relates the construction
to crystalline cohomology \cite[1.5]{bib:ber6} to check that
the trace is an isomorphism when $X$ is geometrically irreducible. We
will check the latter instead by reducing to the case of affine space
(Lemma~\ref{lem:onedim});
in any case,
we already have this for $X = \AAA^n$ by Corollary~\ref{cor:afftop},
which will suffice for the moment.


\subsection{Review of $p$-adic functional analysis}
\label{subsec:functional}

To make more progress, we need some
basic facts from $p$-adic functional analysis; we state
the relevant results here and refer to \cite{bib:crew2} 
(or \cite{bib:schneider}) for further discussion.
For brevity, in this section we use the word ``space'' to mean
``topological vector space over $K$''.

\begin{defn}
Define a subset $S$ of a space $V$ to be
\begin{itemize}
\item \emph{convex} if for any $\bv, \bw, \bx \in S$ and any $a,b,c \in \calO$
with $a+b+c=1$, $a\bv + b \bw + c\bx \in S$;
\item \emph{balanced} if $0 \in S$ (normally only said of convex sets);
\item \emph{absorbing} if for any $\bv \in V$, we have
$\bv \in \lambda S$ for some $\lambda \in K$;
\item \emph{bounded} if for any open neighborhood $U$ of zero, we have
$S \subseteq \lambda U$ for some $\lambda \in K$;
\item \emph{linearly compact} if $S$ is convex, and every convex filter
on $S$ (i.e., every filter admitting a base of convex sets)
has an accumulation point.
\end{itemize}
\end{defn}

\begin{defn}
Define a space $V$ to be 
\begin{itemize}
\item
\emph{locally convex} if its zero element admits
a neighborhood basis consisting of
convex subsets;
\item
\emph{Fr\'echet} if it is Hausdorff and complete,
and its topology is induced by a countable collection of seminorms
(and hence is locally convex);
\item
\emph{LF} (for ``limit-of-Fr\'echet'') if it is separated and locally convex, and
is a direct limit of Fr\'echet spaces;
\item
\emph{barreled} if any closed, convex, balanced, absorbing subset
is a neighborhood of zero;
\item
\emph{Montel} if it is barreled, and every closed convex bounded subset
is linearly compact.
\end{itemize}
\end{defn}

\begin{defn}
The \emph{strong dual} of a space $V$ is the space
$V^\dual$ of continuous linear maps $f: V \to K$, topologized by decreeing
a basis of open sets of the form $\{ f: V \to K | f(B) \subseteq U\}$ for
$B \subseteq V$ bounded and $U \subseteq K$ open.
\end{defn}

\begin{defn} \label{defn:strict}
A continuous linear map $f: V \to W$ of spaces is
\emph{strict} if the induced continuous bijection $V/\ker(f) \to \im(f)$ is a
homeomorphism. 
As noted in \cite[3.7]{bib:crew2}, given a finite collection of continuous
linear maps $f_i: V_i \to W_i$ between locally convex spaces, the map
$\oplus f_i: \oplus V_i \to \oplus W_i$ is strict (resp.\ strict with closed image)
if and only if each $f_i$
is strict (resp.\ strict with closed image).
\end{defn}

The notion of strictness is somewhat delicate: for example, a composition
$g \circ f$ of strict maps with closed image is only guaranteed to be strict
with closed image if $f$ is surjective or $g$ is injective. Otherwise
it may not
even be strict, as in the example (suggested by the referee)
\[
K \langle x \rangle \to K \langle x,y \rangle \to K \langle x,y \rangle/(x-py)
\cong K \langle y \rangle.
\]
A number of useful criteria for strictness appear in \cite[\S 3]{bib:crew2};
here is one of them \cite[Corollary~3.12]{bib:crew2}.
\begin{prop} \label{prop:strdual}
Let $f: V \to W$ be a continuous linear map between Fr\'echet-Montel
spaces. Then $f$ is strict if and only if $f^\dual: W^\dual \to V^\dual$
is strict.
\end{prop}

We also need the following form of the Hahn-Banach theorem; see for example
\cite[Corollary~9.4]{bib:schneider} (or \cite[Th\'eor\`eme~V.2]{bib:monna}). 
Note that we need here the fact
that the valuation on $K$ is discrete, so that $K$ is spherically complete.
Note also that we will only use this result in case $W$ is closed.
\begin{prop} \label{prop:hb}
Let $V$ be a locally convex space, let $W$ be a subspace of $V$,
and let $f: W \to K$
a continuous linear map. Then there exists a continuous linear map
$g: V \to K$ extending $f$.
\end{prop}

\subsection{Topology of $(\sigma, \nabla)$-modules}
\label{subsec:topology}

All of the rings we have been considering (dagger algebra, fringe
algebras, Robba rings) carry topologies, so it is natural to fit them
into the framework of the previous section. Let us see how this is done.

We start with some simple observations.
\begin{itemize}
\item
A dagger algebra with its fringe topology is an LF-space. This can be
verified by noting that $K \langle x_1, \dots, x_n \rangle^\dagger$ is an
LF-space (as in \cite[Proposition~5.5]{bib:crew2}), 
then noting that a separated quotient of an LF-space is again an LF-space.
\item
For any dagger algebra $A$ and any finite set $T$,
the rings $\calR^T_{A,r}$ and $\calR^{T,+}_A$ are Fr\'echet; hence
$\calR^T_A$ is an LF-space.
\item
For any dagger algebra $A$, the continuous
map $A \langle x \rangle^\dagger \to \calR_A$
carrying $x$ to $t^{-1}$ is strict with closed image,
as is the map $\calR_A^+ \to \calR_A$.
Similarly, the map $K \langle x_1, \dots, x_n \rangle^\dagger
\to K \langle x_1, \dots, x_n, x_1^{-1} \rangle^\dagger$ is strict with closed image.
\item
For any affinoid algebra $A$ and any $f \in A^{\inte}$ whose
reduction is not a zero divisor, the localization $A \to A \langle f^{-1} \rangle$
is an injective isometry, and hence is strict.
\item
If $A \to B$ is a strict morphism of reduced affinoid or dagger algebras, then
the map $\calR_A \to \calR_B$ is also strict.
\item
Any morphism $R \to S$ satisfying one of the conditions of 
Proposition~\ref{prop:trace} or Remark~\ref{rem:byhand} is strict with closed image:
the topological vector space $S$ splits as the direct sum of $R$ with the kernel of
the trace map.
\end{itemize}

\begin{prop}
Let $A$ be an integral dagger algebra of MW-type, and suppose $f \in A$ satisfies
$|f|_{\sup,\widehat{A}} = 1$. Then the map
$A \to A \langle f^{-1} \rangle^\dagger$ is strict.
\end{prop}
\begin{proof}
Since it is safe to pass from $K$ to a finite extension at any point in the argument,
we will do so without further comment.
Let $X$ be the special fibre of $A$, so that $X$ is smooth over $k$.
We are trying to show that the passage to the localization of $A$ corresponding
to the passage from $X$ to a particular basic open subset $U$ is strict; it suffices
to show the same thing after passing to an even smaller $V$, 
or after passing from $X$
to a cover by open affines.

Note that the claim is clear in case there is a strict injection
$A \hookrightarrow \calR^+_B$ carrying $f$ to $t$, as we already know that
$\calR^+_B \to \calR_B$ is a strict injection. This injection exists if and only
if the completion of $X$ along $X \setminus U$ can be written as a fibration
over $X \setminus U$; this is true locally if $X \setminus U$ is smooth, so
we are done in that case.

In general, we can now cut out any smooth subvariety of $X$ by embedding it into
a smooth divisor (we may have to pass to a cover of $X$ in order for the
complement of that divisor to become a basic open). 
Now stratify $X \setminus U$ as follows. First enlarge $k$ if
needed to ensure that each component of $X \setminus U$ is generically smooth.
Then discard the smooth locus of each component and repeat until the process
is forced to terminate by the noetherianness of $X$. We can now cut out the
elements of the stratification one by one, from smallest to largest, so as
to ultimately cut out all of $X \setminus U$ (and possibly more). This yields
the claim.
\end{proof}

We turn next to the topology of $(\sigma, \nabla)$-modules and their
cohomology. The basic observation we need is the following, which follows
directly from Proposition~\ref{prop:deriv}.
\begin{prop} \label{prop:strict1}
For any reduced dagger algebra $A$,
the map $d: \calR_A \to \Omega^1_{\calR_A/A}$ is strict with closed image.
\end{prop}
From this simple fact we gain quite a lot of mileage.
\begin{prop} \label{prop:strict2}
Let $A$ be an integral dagger algebra of MW-type,
and let $M$ be a $(\sigma, \nabla)$-module over $\calR_A$, relative
to $A$. Then there exists a localization $A'$ of $A$ such that
$\nabla: M \otimes_{\calR_A} \calR_{A'} \to M \otimes_{\calR_A} 
\Omega^1_{\calR_{A'}/A'}$
is strict with closed image, and remains so after further localization.
\end{prop}
\begin{proof}
Set notation as in Theorem~\ref{thm:monodromy2}.
Then $M \otimes_{\calR_A} \calR'' \to M \otimes_{\calR_A} \Omega^1_{\calR''/A}$
is strict by Proposition~\ref{prop:strict1} plus a direct calculation
(as in Lemma~\ref{lem:freegen}). By taking traces and invoking the criterion
of Definition~\ref{defn:strict}, we obtain the desired strictness. (We get stuck
at $A'$ because we cannot take a trace to pass from $A'$ down to $A$.)
\end{proof}

The following result has already been invoked via
Lemma~\ref{lem:closedim}, but neither that result nor any of its
consequences have been invoked in this chapter, so no vicious circle has been
created.
\begin{prop} \label{prop:strict}
Let $A$ be an integral dagger algebra of MW-type,
and let $M$ be a $(\sigma, \nabla)$-module over $A \langle x \rangle^\dagger$,
relative to $A$.
Then there exists a localization $A'$ of $A$ such that
$\nabla_v: M \otimes_{A \langle x \rangle^\dagger} A' \langle x \rangle^\dagger
\to M \otimes_{A \langle x \rangle^\dagger} 
\Omega^1_{A' \langle x \rangle^\dagger/A'}$ is strict with closed image,
and remains so after further localization.
\end{prop}
\begin{proof}
Replace $A$ with a localization for which Theorem~\ref{thm:fingen}
and Proposition~\ref{prop:strict2} hold for $M$.
We first check that $\im(\nabla_v)$ is closed.
Since the horizontal rows in \eqref{eq:presnake} are strict exact, 
$\im(\nabla_v)$ is homeomorphic to the image of the composite map 
$M \to M \otimes \Omega^1_{\calR_A/A}$, and it suffices to check that the
latter is closed. Since $\nabla_v^{\loc}$ is strict with closed image, 
it also suffices to
check that $M + \ker(\nabla_v^{\loc})$ is closed in $M \otimes \calR_A$,
or that the image of $\ker(\nabla_v^{\loc})$ in
$(M \otimes \calR_A)/M$ is closed. But the latter image is 
finitely generated over $A$, hence always closed.

Since $\im(\nabla_v)$ is closed, $R^1 f_* M$ is separated;
after possibly localizing $A$ further, we can ensure that $R^1 f_* M$ is
finite free over $A$.
Choosing free generators and
lifting to $M \otimes_{A \langle x \rangle^\dagger} \Omega^1_{A \langle x
\rangle^\dagger/A}$, we see that $\im(\nabla_v)$ admits a topological 
complement and so
is a direct summand. By
\cite[Proposition~3.5]{bib:crew2} (applicable since the source and
target of $\nabla_v$ are LF-spaces), $\nabla_v$ is strict.
\end{proof}

\begin{prop} \label{prop:strict3}
Let $M$ be a $(\sigma, \nabla)$-module over a dagger algebra $A$,
and suppose that the spaces $H^i(M)$ are finite-dimensional and separated for
all $i$. Then all of the differentials $\nabla_i: M \otimes_A \Omega^i_{A/K}
\to M \otimes_A \Omega^{i+1}_{A/K}$ are strict with closed image.
\end{prop}
\begin{proof}
The map $f: M \otimes_A \Omega^i_{A/K} \to \ker(\nabla_{i+1})$ has cokernel
$H^{i+1}(M)$, which by hypothesis is finite-dimensional and separated, so
the map is strict with closed image by \cite[Corollary~3.6]{bib:crew2}.
The map $\nabla$ factors as the composition of $f$ with the
inclusion of a closed subspace, so it is also strict with closed image.
\end{proof}

\subsection{Poincar\'e duality on affine space}

Having given an explicit description of the cohomology with compact supports
of an overconvergent $F$-isocrystal on $\AAA^n$, we can now explicitly
verify the nondegeneracy of the  Poincar\'e pairing. We continue to appeal to
\cite{bib:crew2} for topological details. Note that 
Hypothesis~\ref{hypo:affine} and all subsequent definitions and notations
are in force throughout this section.

Recall that $\calR^{\{t_1,\dots,t_n\},+}_K$ is a Fr\'echet space
and observe that it is also Montel, by the same proof as
\cite[Corollary~5.3]{bib:crew2}.
Then note that the natural map from
$t_1\cdots t_n \calR^{\{t_1,\dots,t_n\},+}_K$ to
$\calS^{(n)}_r/\sum_i \calS^{(n),i'}_r$ is a homeomorphism,
so the latter is also a Fr\'echet-Montel space.

The following basic result is proved in the same manner as
\cite[Proposition~5.5]{bib:crew2} (which is essentially the one-dimensional
case).
\begin{lemma} \label{lem:dual}
The pairing
\[
K \langle x_1, \dots, x_n \rangle^\dagger
\otimes_K \calS^{(n)}_r/\sum_i \calS^{(n),i'}_r \to K
\]
taking $a \otimes b$ to the coefficient of $t_1 \cdots t_n$ in $ab$
is perfect: that is, the induced map of each space to the strong dual of the
other is an isomorphism of topological vector spaces.
\end{lemma}

\begin{defn}
Given a $(\sigma, \nabla)$-module $M$ on $K \langle x_1, \dots, x_n 
\rangle^\dagger$, put
\[
C^{i} = \frac{M \otimes \Omega^i_{\calS^{(n)}_r/K}}{\sum_{j=1}^n
M \otimes \Omega^i_{\calS^{(n),j'}_r/K}},
\]
and put
\[
D^i = M^\dual \otimes \Omega^i_{K\langle x_1, \dots, x_n \rangle^\dagger/K},
\]
topologized using the fringe
topology of $K\langle x_1, \dots, x_n \rangle^\dagger$.
\end{defn}

We quickly recall the pairing $[\cdot,\cdot]: C^i \otimes D^{n-i} \to K$
on the chain level that induces the Poincar\'e pairing described in
Section~\ref{subsec:poinpair}.
Given elements of $C^i$ and $D^{n-i}$ represented by
$\bv \in M \otimes \Omega^i_{\calS^{(n)}_r/K}$ and 
$\bw \in M^\dual \otimes \Omega^{n-i}_{K 
\langle x_1, \dots, x_n \rangle^\dagger/K}$,
pair $\bv$ and $\bw$ and apply the natural pairing $M \otimes M^\dual
\to K$ to obtain an element of $\Omega^n_{\calS^{(n)}_r/K}$. Then take
the residue (the
coefficient of $(dt_1 \wedge \cdots \wedge dt_n)/(t_1\cdots t_n)$)
and call the result $[\bv, \bw]$.

The pairing $[\cdot,\cdot]$
expresses each of $C^i$ and $D^{n-i}$ as the strong dual of the
other by Lemma~\ref{lem:dual}. It also
has the property that $[\bv, \nabla \bw] = -[\nabla \bv, \bw]$ (by
the Leibniz rule),
and so induces a pairing on the homologies of the complexes $C^.$
and $D^.$. If $M$ corresponds to an overconvergent $F$-isocrystal
$\calE$ on $\AAA^n$, that induced pairing corresponds to
the Poincar\'e pairing $H^{n+i}_{c,\rig}(\AAA^n, \calE) \otimes_K
H^{n-i}_{\rig}(\AAA^n, \calE^\dual) \to K$.

\begin{prop} \label{prop:affpoincare}
Let $\calE$ be an overconvergent $F$-isocrystal on $\AAA^n$,
and assume that $H^i_{\rig}(\AAA^n, \calE^\dual)$ is finite dimensional and
separated for each $i$.
Then the Poincar\'e pairing 
$H^{n+i}_{c,\rig}(\AAA^n, \calE) \otimes_K
H^{n-i}_{\rig}(\AAA^n, \calE^\dual) \to K$
is nondegenerate (i.e., each factor injects into the dual of the other).
\end{prop}
\begin{proof}
By Proposition~\ref{prop:strict3}, $\nabla: D^{n-i} \to D^{n-i+1}$ is 
strict for each $i$; by Proposition~\ref{prop:strdual}
(applicable because each $D^i$ is Fr\'echet-Montel),
$\nabla: C^i \to C^{i+1}$ is also strict for each $i$.

We first check that the Poincar\'e pairing is nondegenerate on the
right. For any $\bv \in D^{n-i} \setminus \im(D^{n-i-1} \to D^{n-i})$, we may
use the fact that $\im(D^{n-i-1} \to D^{n-i})$ is closed 
(by Proposition~\ref{prop:strict3}) to define a continuous
linear map $\im(D^{n-i-1} \to D^{n-i}) + K\bv \to K$ sending
$\im(D^{n-i-1} \to D^{n-i})$ to zero and $\bv$ to 1. By
Proposition~\ref{prop:hb},
this map extends to a continuous linear
map $f: D^{n-i} \to K$, which corresponds by duality
to some $\bw \in C^{i}$.
Now $[\nabla \bw, \bx] = -[\bw, \nabla \bx] = 0$ for all $\bx \in D^{n-i-1}$,
so by perfectness of the residue pairing, $\nabla \bw = 0$. Thus $\bw$
represents an element of $H^{n+i}_{c,\rig}(\AAA^n, \calE)$, which pairs
with the given element of $H^{n-i}_{\rig}(\AAA^n, \calE^\dual)$ to a nonzero
element of $K$. Hence the Poincar\'e pairing is nondegenerate on the right.

We next check that the Poincar\'e pairing is nondegenerate on the left.
Given $\bw \in C^{i}$ such that $[\bw, \bv] = 0$ for any $\bv \in
D^{n-i}$ such that $\nabla \bv = 0$, we obtain from $\bw$ a continuous
linear map $f: D^{n-i} / \ker(D^{n-i} \to D^{n-i+1}) \to K$. Since $\nabla:
D^{n-i} \to D^{n-i+1}$ is strict, the map $f$ induces a continuous linear map
$g: \im(D^{n-i} \to D^{n-i+1}) \to K$. By Proposition~\ref{prop:hb},
$g$ extends to a continuous linear map $h: D^{n-i+1} \to K$, which corresponds
by duality to an element $\bx$ of $C^{i-1}$. Now for any $\bv \in D^{n-i}$,
\[
[\bw + \nabla(\bx), \bv] = [\bw, \bv] - [\bx, \nabla(\bv)]
= f(\bv) - g(\nabla(\bv)) = 0.
\]
Thus $\bw + \nabla(\bx) = 0$ by perfectness of the residue pairing,
so $\bw$ vanishes in cohomology.
Therefore the Poincar\'e pairing is nondegenerate on the left,
completing the proof.
\end{proof}

\subsection{Relative Poincar\'e duality}
\label{subsec:relpoin}

We mention in passing that Poincar\'e duality also holds
in the relative setting. Although we will not use this fact hereafter,
it provides
a useful tool for controlling the shape of the pushforwards; see
for instance \cite{bib:me10}. 

\begin{prop} \label{prop:relpoincare}
Let $A$ be an integral dagger algebra of MW-type,
and let $f: A \to A \langle x \rangle^\dagger$ be the canonical inclusion.
Let $M$ be a $(\sigma, \nabla)$-module over $A \langle x \rangle^\dagger$,
and suppose $M$ and $M^\dual$ satisfy the conclusion
of Theorem~\ref{thm:fingen} without further localization. 
Then there are canonical perfect pairings
\begin{align*}
R^0 f_* M \otimes_A R^{2} f_! M^\dual &\to A  \\
R^i_{\loc} f_* M \otimes_A R^{1-i}_{\loc} f_* M^\dual &\to A \\
R^1 f_! M \otimes_A R^1 f_* M^\dual &\to A \\
R^1_p f_* M \otimes_A R^1_p f_* M^\dual &\to A.
\end{align*}
\end{prop}
Again, these pairings are induced by the natural pairing 
$M \otimes_A M^\dual \to A$ and the residue map $\Omega^1_{\calR_A/A} \to A$
sending $\omega$ to its coefficient of $dt/t$. The fourth pairing exists
because of a compatibility between the others: given $\alpha \in
R^0_{\loc} f_* M$ mapping to $\beta \in R^1 f_!M$
and $\gamma \in R^1 f_* M^\dual$ mapping to $\delta \in 
R^1_{\loc} f_* M^\dual$,
then $[\alpha, \delta] = [\beta, \gamma]$.
\begin{proof}
By flat base change to the generic point (Remark~\ref{rem:flatfrac2}),
Poincar\'e duality on $\AAA^1$ (Proposition~\ref{prop:affpoincare}, whose
finite-dimensionality and separatedness hypotheses are satisfied 
as in the proof of
\cite[Theorem~9.5]{bib:crew2}), 
and the finite dimensionality of
$R^0 f_* M$ and $R^2 f_! M^\dual$ (Theorem~\ref{thm:fingen}), 
the map
\[
R^0 f_* M \to (R^2 f_! M^\dual)^\dual
\]
becomes an isomorphism upon tensoring with the completion of
$\Frac A$ for the affinoid topology. However, the kernel and cokernel
of this map are finitely generated
$A$-modules with connection, so by Lemma~\ref{lem:locfree}
they are locally free, necessarily of rank zero. Thus in fact they both vanish,
and the map is an isomorphism, that is, the first pairing above is perfect.
The same argument applies to $R^1 f_* M \to (R^1 f_! M^\dual)^\dual$;
the perfectness of the other two pairings follows from the five lemma.
\end{proof}

\section{Finite dimensionality: the main results}
\label{sec:final}

In this chapter, we prove the results enunciated in Chapter~1.
While up to now we have made all our calculations of rigid cohomology
on smooth affine schemes, where the Monsky-Washnitzer description can
be used, at this point we wish to prove more general results. We will
thus use the formalism of rigid cohomology in full strength, reducing
everything to consideration of smooth affine schemes.
See Chapter~\ref{sec:rigid} for more details.

\subsection{Cohomology without supports: the smooth case}

In this section, we prove the restriction of Theorem~\ref{thm:finite}
to the case where $X$ is smooth.
\begin{prop} \label{prop:finite}
Let $\calE$ be an overconvergent $F$-isocrystal on a smooth scheme
$X$ over a field $k$ of characteristic $p>0$. Then
the rigid cohomology spaces $H^i_{\rig}(X/K, \calE)$ are finite
dimensional for all $i$.
\end{prop}
The gross outline of the proof follows that given 
by Berthelot \cite{bib:ber2} in the constant coefficient
case, that given by Tsuzuki \cite[Theorem~6.1.1]{bib:tsu5} in the
unit-root case, and the proposed proof of Shiho
\cite{bib:shiho} in the general case. However, we simplify matters
using the following geometric fact, for whose proof see \cite{bib:me11}.
(A proof in case $k$ is infinite and perfect can be found in \cite{bib:me9}.)
\begin{prop} \label{prop:belyi}
Let $X$ be a separated, geometrically reduced,
finite type $k$-scheme of dimension $n$, for $k$ a field of
characteristic $p>0$, and let $x$ be a smooth geometric point of $X$. Then
$X$ contains an open dense affine subscheme containing $x$ and
admitting a finite \'etale morphism to affine $n$-space.
\end{prop}
This result allows us to ``go down'' to simpler spaces (namely affine space).
By contrast, in \cite{bib:ber2}, \cite{bib:tsu5}, and \cite{bib:shiho}
(and Section~\ref{subsec:descent}, where we invoke cohomological descent),
one simplifies spaces by ``going up'' using
de~Jong's alterations theorem.

\begin{proof}[Proof of Proposition~\ref{prop:finite}]
We prove the following two assertions by induction on $n$:
\begin{enumerate}
\item[(a)$_n$]
$H^i_{\rig}(X/K, \calE)$ is finite dimensional for any $X$ smooth over $k$ of
dimension at most $n$, and any overconvergent $F$-isocrystal $\calE$
on $X$;
\item[(b)$_n$]
$H^i_{Z,\rig}(X/K, \calE)$ is finite dimensional for any closed immersion
$Z \hookrightarrow X$ of a geometrically reduced
scheme $Z$ of finite type over $k$ of dimension
at most $n$ into a scheme $X$ smooth over $k$, and for any overconvergent
$F$-isocrystal $\calE$ on $X$.
\end{enumerate}
For starters, (a)$_n$ is trivial for $n=0$ and (b)$_n$ is vacuous
for $n=-1$. Now (b)$_{n-1}$ and (a)$_n$ together imply (b)$_n$ as
follows. (The argument is copied essentially verbatim from
\cite[Theorem~3.1.9]{bib:shiho}.) First suppose that
Hypothesis~\ref{hypo:gysin} is satisfied; that is,
$Z \hookrightarrow X$ is a smooth pair, with $\dim Z \leq n$, admitting a
lifting $\mathcal{Z} \hookrightarrow \mathcal{X}$ over $\calO$,
such that $\mathcal{Z}, \mathcal{X}$ are smooth over $\calO$ and
there exists an \'etale 
morphism $\mathcal{X} \to \Spec \calO[x_1, \dots, x_d]$ for some integer $d$,
such that $\mathcal{Z}$ maps to the vanishing locus of $d-n$ of the $x_i$.
Then we may apply Tsuzuki's construction of the Gysin isomorphism
\cite[Theorem~4.1.1]{bib:tsu5} to conclude that $H^i_{Z,
\rig}(X/K, \calE) \cong H^{i-2c}_{\rig}(Z/K, \left. \calE \right|_Z)$,
where $c = \dim X - \dim Z$; by (a)$_{n}$, the latter is finite dimensional.

To handle the general case, we use the excision exact sequence 
\eqref{eq:excis1}, which we recall for convenience:
\[
\cdots \to H^i_{T/K,\rig}(X, \calE) \to H^i_{Z/K,\rig}(X, \calE)
\to H^i_{Z \setminus T,\rig}(X \setminus T/K, \calE) \to \cdots,
\]
Suppose $Z$ is any geometrically reduced
subscheme of a smooth $k$-scheme $X$ with
$\dim Z \leq n$. Then we can find
some closed subscheme $T$ with $\dim T \leq n-1$ such that
Hypothesis~\ref{hypo:gysin} holds for the pair
$Z \setminus T \hookrightarrow X \setminus T$. Thus
$H^i_{Z \setminus T,\rig}(X \setminus T/K,\calE)$ is finite dimensional;
by (b)$_{n-1}$, $H^i_{T,\rig}(X/K, \calE)$ is also finite dimensional.
Hence by \eqref{eq:excis1}, $H^i_{Z,\rig}(X/K, \calE)$ is finite dimensional;
we thus conclude that (b)$_{n-1}$ and (a)$_n$ imply (b)$_n$.

Now we prove that (a)$_{n-1}$ and (b)$_{n-1}$ imply (a)$_n$.
To prove finite dimensionality of $H^i_{\rig}(X/K, \calE)$
for an overconvergent $F$-isocrystal
$\calE$ on $X$,
we may make the following replacements without loss of generality:
\begin{enumerate}
\item[(i)]
replace $K$ by a finite extension;
\item[(ii)]
replace $(X, \calE)$ by $(U, j^* \calE)$ for any
open immersion $j: U \hookrightarrow X$ with dense image;
\item[(iii)]
replace $(X, \calE)$, if $X$ is affine,
by $(Y, f_* \calE)$ for any finite \'etale
morphism $f: X \to Y$.
\end{enumerate}
The reduction (i) follows from Proposition~\ref{prop:finpull}
(or by finite flat base change, as in the next section).
The reduction (ii) follows from (b)$_{n-1}$ and the excision
sequence \eqref{eq:excis1}
with $T$ the reduced subscheme of
$X \setminus U$ and $Z=X$ (and applying (i) if necessary so
that $T$ becomes geometrically reduced).
The reduction (iii) follows from
Proposition~\ref{prop:finpush}, as the dimensions of the cohomology
spaces do not change.

For starters, (ii) means we can replace $X$ with an open dense
affine subset $V$. By Proposition~\ref{prop:belyi},
there exists
an open immersion $j: U \hookrightarrow V$ with dense image such
that $U$ admits a finite \'etale map $f: U \to \AAA^n$. By (ii) and (iii),
to prove finiteness for $(X, \calE)$, it suffices to prove
finiteness for $(\AAA^n, f_* j^* \calE)$. Now write $\AAA^n$ as a
product $\AAA^1 \times \AAA^{n-1}$. By
Theorem~\ref{thm:fingen},
the kernel $\calF_0$ and cokernel $\calF_1$ of the vertical connection
$\nabla_v$
become overconvergent $F$-isocrystals on some open subset $W$ of $\AAA^{n-1}$.

By (a)$_{n-1}$, $H^i_{\rig}(W/K, \calF_0)$ and
$H^i_{\rig}(W/K,\calF_1)$ are finite dimensional for all $i$; by
Proposition~\ref{prop:leray},
we conclude that $H^i_{\rig}(\AAA^1 \times W/K, f_* j^* \calE)$ 
is finite dimensional for all $i$.
Thus $H^i_{\rig}(X/K, \calE)$ is finite dimensional for all $i$, and the
induction is complete.
\end{proof}

\subsection{Cohomology without supports: the general case}
\label{subsec:descent}

We now proceed to the proof of Theorem~\ref{thm:finite} in general,
relaxing the hypothesis that $X$ is smooth.
We will obtain the theorem from Proposition~\ref{prop:finite} using the
technique of cohomological descent, developed in rigid cohomology by
Chiarellotto and Tsuzuki \cite{bib:ct}, \cite{bib:tsu6}. 

\begin{proof}[Proof of Theorem~\ref{thm:finite}]
First recall a consequence of the finite flat base change theorem
\cite[Theorem~11.8.1]{bib:ct}
(or see \cite[Proposition~1.8]{bib:ber2} for the
constant coefficient case): if $K'$ is a finite extension of $K$
with residue field $k'$,
then for any separated scheme $X$ of finite
type over $k$ and any overconvergent $F$-isocrystal $\calE$ on $X$,
we have
 \[
H^i_{\rig}(X/K, \calE) \otimes_K K' \cong H^i_{\rig}(X \times_k
 k'/K', \calE).
\]

We now imitate the argument of \cite[Theorem~5.1.1]{bib:tsu6} using
nonconstant coefficients. Let $X$ be a separated scheme of finite type 
over $k$, and let
$\calE$ be an overconvergent $F$-isocrystal over $X$. The alterations
theorem of de~Jong \cite[Theorem~4.1]{bib:dej0} and the construction of
\cite[Section~4]{bib:tsu6} produce a proper hypercovering $f_.: Y_.
\to X$ over $\Spec k$ and an increasing sequence $k_0 \subseteq k_1 \subseteq
\cdots$ of finite purely inseparable extensions of $k$, such that
$Y_n$ is smooth over $k_n$ for each $n$.
Choose an increasing sequence $K_0 \subseteq
K_1 \subseteq \cdots$ of finite extensions of $K$ such that $K_n$
has residue field $k_n$ for each $n$. 
Then $Y_n \times_k  k_n$ may not be smooth over $k_n$, but the canonical
immersion $Y_n \hookrightarrow Y_n \times_k k_n$ over $\Spec k_n$
is a nilpotent immersion, so it induces an isomorphism
$H^i_{\rig}(Y_n \times_k k_n/K_n, f_n^* \calE) \cong
H^i_{\rig}(Y_n/K_n, f_n^* \calE)$.
By Proposition~\ref{prop:finite}, these spaces are finite dimensional over
$K_n$, as then is $H^i_{\rig}(Y_n/K, f_n^* \calE)$ over $K$ by the finite flat
base change theorem described above. Now
\cite[Theorem~4.5.1]{bib:tsu6} produces a spectral sequence in which
\[
E_1^{qr} = H^r_{\rig}(Y_q/K, f_q^* \calE) \Rightarrow H^{q+r}_{\rig}(X/K, \calE).
\]
Since all terms of $E_1$ are finite dimensional over $K$, the
spaces $H^i_{\rig}(X/K, \calE)$ are as well.
\end{proof}

\begin{remark}
By excision \eqref{eq:excis1}, it also follows that $H^i_{Z,\rig}(X, \calE)$
is finite dimensional for any closed subscheme $Z$ of $X$, without 
any smoothness hypothesis; one also obtains the base change property.
\end{remark}

\begin{remark}
We believe one can also carry out the above reduction of Theorem~\ref{thm:finite}
to Proposition~\ref{prop:finite} using the techniques used by Grosse-Kl\"onne
\cite{bib:gk2} in the constant coefficient case. (Namely, one should locally
embed the nonsmooth scheme into a smooth scheme and carefully
control the cohomology of some ensuing dagger spaces.) This belief amounts
to the claim that one can repeat the entire argument of \cite{bib:gk2}
with nonconstant coefficients (coherent modules with connection on dagger
spaces), using Proposition~\ref{prop:finite} in lieu of \cite{bib:ber2}.
We have not verified this claim.
\end{remark}

\subsection{Poincar\'e duality and cohomology with compact supports}

We now consider rigid
cohomology with compact supports. In this case, the relevant
excision sequence, for $Z$ a closed subscheme of $X$ and $U = X \setminus Z$,
is \eqref{eq:excis2}:
\[
\cdots \to H^i_{c,\rig}(U, \calE) \to
H^i_{c,\rig}(X, \calE) \to H^i_{c,\rig}(Z, \calE) \to \cdots
\]
Since all three types of terms in the sequence have the same form,
the proof of finite
dimensionality is somewhat simpler in this case, and in particular
there is no need to separately handle the smooth case, use a Gysin
isomorphism or invoke any descent results. All we need to do is first
check Poincar\'e duality, then use Theorem~\ref{thm:finite}.

As in the case of cohomology without supports, we have a base
change isomorphism
 \[
H^i_{c,\rig}(X/K, \calE) \otimes_K K' \cong H^i_{c,\rig}(X \times_k
 k'/K', \calE)
\]
for $K'$ a finite extension of $K$ with residue field $k'$,
and the proof is analogous; again, see \cite[Proposition~1.8]{bib:ber2}.

Before proving anything about Poincar\'e duality, we need to clarify
the nature of its target space. We do this by starting with some
basic computations; these duplicate
\cite[Propositions~1.1 and~1.5]{bib:ber6}, but we prefer to reproduce
them in our context (in particular, without comparison to crystalline
cohomology).
\begin{lemma} \label{lem:vanishes}
Let $X$ be a separated, finite type $k$-scheme of dimension $n$. 
Then $H^{d}_{c,\rig}(X/K) = 0$ for $d > 2n$.
\end{lemma}
\begin{proof}
We proceed by induction on $n$. 
We may ensure that $X$ is geometrically
reduced by making a base change.
By excision \eqref{eq:excis2}, we may
replace $X$ by an open dense subset $U$; in particular, we choose
$U$ by virtue of Proposition~\ref{prop:belyi} so that
there exists a finite \'etale map $f: U \to \AAA^n$.
Let $\calE$ be the pushforward along $f$ of the trivial
$F$-isocrystal on $U$; then $H^d_{c,\rig}(U/K) \cong 
H^d_{c,\rig}(\AAA^n/K, \calE)$ by Proposition~\ref{prop:finpush2}.
However, the latter manifestly vanishes from the description
given in Section~\ref{subsec:cohomaff}. This yields the desired result.
\end{proof}

Again, let $X$ be a separated, finite type $k$-scheme of pure dimension $n$.
Using Lemma~\ref{lem:vanishes}, one may construct a canonical trace morphism
\[
\Tr: H^{2n}_{c,\rig}(X/K) \to H^0_{c,\rig}(\Spec k/K) \cong K
\] 
as in \cite[1.2]{bib:ber6}.
That construction completes the description of the Poincar\'e pairing begun in
Section~\ref{subsec:poinpair}, and now we may begin to consider
its perfectness.

\begin{lemma} \label{lem:onedim}
Let $X$ be a separated, finite type, geometrically irreducible
$k$-scheme of dimension $n$. Then $H^{2n}_{c,\rig}(X/K)$ is one-dimensional
over $K$.
\end{lemma}
Note that this assertion does not include the fact that the trace map
on $X$ is nonzero; that will follow from the general proof of Poincar\'e
duality.
\begin{proof}
Again, we can reduce to the case where $X$ is geometrically reduced
by a base change.
By excision \eqref{eq:excis2} and Lemma~\ref{lem:vanishes},
we may again replace $X$ by an open dense subset admitting a finite \'etale
map to $\AAA^n$.
Viewing $\AAA^n$ as a family of $\AAA^1$'s with base $\AAA^{n-1}$,
we may shrink the base so that pushforwards exist as in Theorem~\ref{thm:fingen},
and the conclusion of Proposition~\ref{prop:strict2} holds.
Shrink the base and push forward again to convert into another
$\AAA^{n-1}$, then view $\AAA^{n-1}$ as $\AAA^1 \times \AAA^{n-2}$, shrink
$\AAA^{n-2}$ until we can push forward again, et cetera.

At the end of this process, we obtain an open dense subset $U$ of $X$
admitting a finite \'etale map $f: U \to \AAA^n$. 
Retain notation as in Lemma~\ref{lem:vanishes}, and
topologize the spaces
$H^i(\AAA^n/K, \calE)$ by identifying
them with the cohomology of the corresponding $(\sigma, \nabla)$-module $M$
over $K \langle x_1, \dots, x_n \rangle^\dagger$. In the argument (and notation) of
Proposition~\ref{prop:leray}, the differentials $d_0$
are strict by Proposition~\ref{prop:strict},
so the cohomology of $H^i(M)$ is correctly computed as a topological
space by $E_1$.
Suppose that the hypothesis of Proposition~\ref{prop:strict3} is satisfied
for $P$ and $Q$; then by Proposition~\ref{prop:strict3}, the
differentials $d_1$ are strict, so the cohomology of $H^i(M)$ is correctly
computed by $E_2$. 
That is, the exact sequence of Proposition~\ref{prop:leray}
is strict, so the $H^i(M)$ are finite dimensional and separated.
Applying the same argument to $P$ and $Q$, we recurse down the dimensions
to ultimately conclude that the finite dimensionality
and separatedness conditions in 
Proposition~\ref{prop:affpoincare} are satisfied.

We now may invoke Poincar\'e duality on $\AAA^n$ (Proposition~\ref{prop:affpoincare})
to ensure that
$H^{2n}_{c,\rig}(\AAA^n/K,\calE)$ is canonically isomorphic to the
dual of $H^0_{\rig}(\AAA^n/K,\calE^\dual)$, which in turn by
Proposition~\ref{prop:finpush} is canonically isomorphic to
$H^0_{\rig}(U/K)$. (Here we are using the fact that $\calE^\dual \cong \calE$
via the trace from $U$ to $\AAA^n$.)
But the latter is canonically isomorphic to
$H^0_{\rig}(\Spec k/K)$ by \cite[Theorem~7.1]{bib:mw}, whence the desired 
result.
\end{proof}

For details on the functoriality of the Poincar\'e pairing, see
\cite[Section~6.2]{bib:tsu5}. One fact we need is that for $X$ smooth,
the excision sequences
\eqref{eq:excis1} and \eqref{eq:excis2} fit into a commutative diagram
\begin{equation} \label{eq:exdiagram}
\xymatrix{
\cdots \ar[r] & H^{2d-i}_{c,\rig}(Z/K, \calE^\dual) \ar[d]
\ar[r] & H^{2d-i}_{c, \rig}(X/K, \calE^\dual) \ar[d]
\ar[r] & H^{2d-i}_{c, \rig}(U/K, \calE^\dual) \ar[d]
\ar[r] & \cdots \\
  \cdots \ar[r] & H^i_{Z,\rig}(X/K, \calE)^\dual  \ar[r] & H^i_{\rig}(X/K, \calE)^\dual
\ar[r]  & H^i_{\rig}(U/K, \calE)^\dual \ar[r] & \dots
}
\end{equation}
in which the vertical arrows are induced by the pairing \eqref{eq:pairing}.
Note that the spaces being dualized in the bottom row are already known to be
finite dimensional by Theorem~\ref{thm:finite}, so there is no ambiguity about
what sort of duals to take.

We exploit this diagram via an inductive procedure similar to
that of \cite[Theorem~6.2.5]{bib:tsu5}. 

\begin{proof}[Proof of Theorem~\ref{thm:poincare}]
We prove the following assertions by simultaneous induction on $n$,
where (a)$_0$ is straightforward.
\begin{enumerate}
\item[(a)$_n$] 
The pairing \eqref{eq:pairing} is perfect whenever $Z = X$ and
$\dim X \leq n$.
\item[(b)$_n$]
The pairing \eqref{eq:pairing} is perfect whenever $\dim Z \leq n$.
\end{enumerate}
We first observe that (b)$_{n-1}$ and (a)$_n$ imply (b)$_n$, by the same argument
as in the proof of Theorem~\ref{thm:finite}, only this time using the
commutative diagram \eqref{eq:exdiagram}
in place of the excision sequence \eqref{eq:excis1}.

We next establish that (b)$_{n-1}$ implies (a)$_{n}$. 
Given $X$ and $\calE$, we may choose an
open dense subscheme $U$ of $X$ admits a finite \'etale morphism
$f: U \to \AAA^n$.
Taking $Z = X \setminus U$, considering the
excision diagram \eqref{eq:exdiagram}, and applying the induction hypothesis
to the left vertical arrow, we see that to prove perfectness of the pairing
for $X$, it suffices to check it for $U$. To prove that, it in turn suffices
by Proposition~\ref{prop:finpush2}
to check perfectness for the isocrystal $f_* \calE$ and its dual on
$\AAA^n$. 

As in the proof of Lemma~\ref{lem:onedim}, we can arrange the map $f$
so that the spaces $H^i_{\rig}(\AAA^n/K, f_* \calE)$ are known to be
finite dimensional and
separated. We may thus apply Proposition~\ref{prop:affpoincare}; 
 the nondegeneracy of the Poincar\'e pairing implies
first that the spaces $H^{2n-i}_{c,\rig}(\AAA^n/K, f_* \calE^\dual)$ are finite
dimensional, then that the pairing is in fact perfect. 
Hence (b)$_{n-1}$ implies (a)$_{n}$, completing the proof.
\end{proof}

We now prove finite dimensionality of rigid cohomology with compact supports.
\begin{proof}[Proof of Theorem~\ref{thm:finite2}]
Again, we induct on the dimension of $X$, the case of dimension -1 being
vacuous. After replacing $k$ by a finite extension and enlarging
$K$ appropriately, we may replace $X$ by its reduced subscheme
and assume that the latter is irreducible and geometrically reduced.
Then the complement of the singular locus $Z$ of $X$ is an open
dense subscheme of $X$.
The cohomology spaces
$H^i_{c,\rig}(Z/K, \calE)$ are finite dimensional by the induction hypothesis,
while the cohomology spaces $H^i_{c,\rig}(X \setminus Z/K, \calE)$ are finite
dimensional by Theorems~\ref{thm:finite} and Theorem~\ref{thm:poincare}.
Thus $H^i_{c,\rig}(X/K, \calE)$ is finite dimensional by the exactness of
the excision sequence
\eqref{eq:excis2}.
\end{proof}

\subsection{The K\"unneth decomposition}

This time, we defer to \cite[Section~6.3]{bib:tsu5} for the construction
and functoriality of the K\"unneth morphisms \eqref{eq:kun1},
\eqref{eq:kun2} and concentrate on showing that these are isomorphisms,
following \cite[Th\'eor\`eme~3.2]{bib:ber6} and 
\cite[Theorem~6.3.6]{bib:tsu5}.

\begin{proof}[Proof of Theorem~\ref{thm:kunneth}]
We first establish that \eqref{eq:kun2} is an isomorphism,
by proving the following assertions by simultaneous induction on $n$.
\begin{enumerate}
\item[(a)$_n$] 
The morphism \eqref{eq:kun2} is an
isomorphism whenever the $X_i$ are smooth, $Z_i = X_i$ ($i=1,2$) and
$\dim X \leq n$.
\item[(b)$_n$]
The morphism \eqref{eq:kun2} is an isomorphism
whenever the $X_i$ are smooth and $\dim Z \leq n$.
\end{enumerate}
To begin with, (b)$_{-1}$ is vacuous and (a)$_0$ is trivial.
To show that (b)$_{n-1}$ and (a)$_{n}$ imply (b)$_n$,
it is easiest to first prove (b)$_n$ assuming $Z_2 = X_2$.
In that case, one can ``tensor'' the excision exact sequence on $X_1$
with the cohomology of $X_2$ and map to $X$, then argue in the
resulting diagram as in Proposition~\ref{prop:finite}.
This is a bit painful to write out in our notation; it is better
put in the language of derived functors and triangulated categories
(as in \cite[3.1]{bib:ber6})
by saying that we are tensoring the distinguished triangle
\[
\RR\Gamma_{T_1,\rig}(X_1/K, \calE_1)
\to
\RR\Gamma_{Z_1,\rig}(X_1/K,\calE_1)
\to
\RR\Gamma_{Z_1\setminus T_1,\rig}(X_1 \setminus T_1/K, \calE_1)
\stackrel{+1}{\to}
\]
with $\RR\Gamma_{Z_2,\rig}(X_2/K,\calE_2)$, and so on.

To deduce that (a)$_{n-1}$ and
(b)$_{n-1}$ imply (a)$_n$, use Proposition~\ref{prop:belyi}
to choose open dense subschemes $U_1, U_2$ of $X_1, X_2$ 
which admit finite \'etale morphisms to affine spaces.
By applying excision (as in the previous paragraph), (b)$_{n-1}$
and Proposition~\ref{prop:finpush},
we may reduce to the case where $X_1$ is an affine
space, $X_2 = V \times \AAA^1$, and $R^i f_* \calE_2$ is an overconvergent
$F$-isocrystal on $V$ (for $f$ the projection of $V \times \AAA^1$ onto $V$).
Let $g$ be the projection of $X_1 \times X_2$ onto $X_1 \times V$,
Let $\pr_1: X_1 \times V \to X_1$ and $\pr_2: X_1 \times V \to V$
denote the obvious projections,
and write $\calF \boxtimes \calG$ for the external product
$\pr_1^* \calF \otimes \pr_2^* \calG$.
We then have a commuting diagram
\[
\xymatrix{
\vdots \ar[d] & \vdots \ar[d] & \vdots \ar[d] \\
\bigoplus_{j+l=i} 
H^j(X_1, \calE_1) \otimes
H^l(V, R^0 f_* \calE_2) \ar[r] \ar[d] &
H^i(X_1 \times V,
\calE_1 \boxtimes (R^0 f_* \calE)) \ar[r] \ar[d] &
H^i(X_1 \times V, R^0 g_* \calE) \ar[d] \\
\bigoplus_{j+l=i} 
H^j(X_1, \calE_1) \otimes
H^l(X_2, \calE_2) \ar[r] \ar[d] &
H^i(X_1 \times X_2, \calE) \ar[r] \ar[d] &
H^i(X_1 \times X_2, \calE) \ar[d] \\
\bigoplus_{j+l=i} 
H^j(X_1, \calE_1) \otimes
H^{l-1}(V, R^1 f_* \calE_2) \ar[r] \ar[d] &
H^{i-1}(X_1 \times V, \calE_1 \boxtimes (R^1 f_* \calE)) \ar[r] \ar[d] &
H^{i-1}(X_1 \times V, R^1 g_* \calE) \ar[d] \\
\vdots & \vdots & \vdots
}
\]
as follows. 
(We dropped ``$K$'' and ``$\rig$'' from the notation to reduce clutter.)
The left and right vertical columns are long exact sequences
obtained from
Proposition~\ref{prop:leray}. The middle column and its arrows to the
right column are constructed from the isomorphisms
of Proposition~\ref{prop:pushext} (the middle horizontal arrow being the
identity map). The arrows between the left and middle columns
are the K\"unneth morphisms; of these, the top and bottom arrows are
isomorphisms by (a)$_{n-1}$. Since this is true all the way down the
long exact sequences, we conclude that the middle arrow is also an
isomorphism by the five lemma. This yields (a)$_n$, completing the induction
and proving that \eqref{eq:kun2} is an isomorphism.

We now establish that \eqref{eq:kun1} is always an isomorphism
by induction on the dimension of $X$.
Given $X$, let $Z$ be the singular locus of $X$ and put
$U = X \setminus Z$. By excision as above and the induction hypothesis,
it suffices to prove the K\"unneth isomorphism on $U$, but
this follows from
Poincar\'e duality (Theorem~\ref{thm:poincare})
together with the fact that \eqref{eq:kun2} is always an isomorphism.
\end{proof}

\end{document}